\documentclass[a4paper,10pt,reqno]{amsart}
 \usepackage[utf8]{inputenc}
\usepackage{amsaddr}

\usepackage{lmodern}  
\usepackage[T1]{fontenc}
\usepackage{tikz}
\usetikzlibrary{calc,trees,positioning,arrows,chains,shapes.geometric,%
    decorations.pathreplacing,decorations.pathmorphing,shapes,%
    matrix,shapes.symbols}
\tikzset{
>=stealth',
  punktchain/.style={
    rectangle, 
    rounded corners, 
    draw=black, very thick,
    text width=25	em, 
    minimum height=3em, 
    text centered, 
    on chain},
  line/.style={draw, thick, <-},
  element/.style={
    tape,
    top color=white,
    bottom color=blue!50!black!60!,
    minimum width=8em,
    draw=blue!40!black!90, very thick,
    text width=10em, 
    minimum height=3.5em, 
    text centered, 
    on chain},
  every join/.style={->, thick,shorten >=1pt},
  decoration={brace},
  tuborg/.style={decorate},
  tubnode/.style={midway, right=2pt},
}
\usepackage{color,soul}

\newcommand{\blue}[1]{\textcolor{blue}{#1}}

\usepackage{bbm}

\usepackage{sidecap}
\usepackage[normalem]{ulem}
\usepackage[margin=.9in]{geometry}
\usepackage{amssymb, cancel, amsmath, amstext}
\usepackage[dvips]{epsfig}
\usepackage[small]{caption}
\usepackage{graphicx}
\usepackage[all]{xy}
\usepackage{mathtools}

\usepackage[colorlinks=true]{hyperref}
\hypersetup{urlcolor=blue, citecolor=red}

\usepackage{caption}
\usepackage{subcaption}
\usepackage{colortbl}
\usepackage{mathrsfs}
\usepackage{multirow}
 \usepackage[shortlabels]{enumitem}
\usepackage{marginnote}
\usepackage{geometry}
\newtheorem{Lemma}{Lemma}[section]
\newtheorem{Theorem}[Lemma]{Theorem}
\newtheorem{Proposition}[Lemma]{Proposition}

\newtheorem{Corollary}[Lemma]{Corollary}
\newtheorem{Remark}[Lemma]{Remark}
\newtheorem{Definition}[Lemma]{Definition}

\newtheorem{Observation}[Lemma]{Observation}
\newenvironment{Proof}%
 {\begin{trivlist} \item[]{\bf Proof. }}%
 {\hspace*{\fill}$\rule{.4\baselineskip}{.4\baselineskip}$\end{trivlist}}

\newcommand{\gm}{g_{-1}}
\newcommand{\gp}{g_{+1}}

\newcommand{\ep}{\varepsilon}

\newcommand{\om}{\omega}
\newcommand{\kp}{\kappa}
\newcommand{\tkp}{\tilde{\kappa}}

\newcommand{\tv}{\widetilde{v}}
\newcommand{\hv}{\widehat{v}}

\newcommand{\n}{_{\text{near}}}
\newcommand{\f}{_{ \text{far}}}
\newcommand{\m}{_{ -1}}
\newcommand{\p}{_{ +1}}
\newcommand{\ur}{u_{\text{rolls}}}

\newcommand{\T}{\mathbb{T}}
\newcommand{\R}{\mathbb{R}}

\newcommand{\N}{\mathbb{N}}

\author{Rafael Monteiro$^*$}

\address{Mathematics for Advanced Materials Research - Open Innovation Lab (MathAM-OIL), Sendai, Japan}
\author{Natsuhiko Yoshinaga}

\address{Mathematics for Advanced Materials Research - Open Innovation Lab (MathAM-OIL), Sendai, Japan}
\address{WPI - Advanced Institute for Materials Research - Tohoku University, Sendai, Japan}

 \title[The SHE under DQ: finding heter. connect. using spatial \& spectral decompositions]{The Swift-Hohenberg Equation under directional-quenching: finding heteroclinic connections using spatial \\and spectral decompositions}

\author[R. Monteiro \&  N. Yoshinaga]{}


 \email{monteirodasilva-rafael@aist.go.jp, rafael.a.monteiro.math@gmail.com}
\email{yoshinaga@tohoku.ac.jp}

\thanks{$^*$ Corresponding author}

\begin{document}

\begin{abstract}
We study the existence of patterns (nontrivial, stationary solutions) in the one-dimensional Swift-Hohenberg Equation in a directional quenching scenario, that is, on $x\leq 0$ the energy potential associated to the equation is bistable, whereas on $x\geq 0$ it is monostable. This heterogeneity in the medium induces a symmetry break that makes the existence of heteroclinic orbits of the type point-to-periodic not only plausible but, as we prove here, true. In this search, we use an interesting result of \cite{Fefferman} in order to understand the multiscale structure of the problem, namely, how multiple scales  - fast/slow- interact with each other. In passing, we advocate for a new approach in finding connecting orbits, using what we call ``far/near decompositions'', relying both on information about the spatial behavior of the solutions and on Fourier analysis. Our method is functional analytic and PDE based, relying minimally on dynamical system techniques and making no use of comparison principles whatsoever.

\end{abstract}

\maketitle


\section{Introduction}\label{section:introduction}

In this paper we study the one-dimensional Swift-Hohenberg Equation  (SHE),
\begin{equation}\label{SH-eq}
 \partial_t u(z,t) = -(1 + \partial_z^2)^2 u(z,t)  + \delta^2\mu(z) u(z,t) - u^3(z,t), \quad z\in \R, \, t\in \R^{+},
\end{equation}
a model originally derived in the study of hydrodynamic instability due to thermal convection \cite{Swift-Hohenberg}. Here, the quantity $u(z,t)$ denotes the concentration of fluid at the point $(z,t)$ of space and time, and $\mu(\cdot)$ is a control parameter that may vary in space, representing the difference in temperature between the bottom and the top of the fluid.   Whenever $\mu(\cdot) \equiv 1$ it is known that \eqref{SH-eq} supports a 3 parameter family $(\delta,\om,\gamma)$ of $\frac{2\pi}{\om}$-periodic solutions,
\begin{align}\label{rolls}
\ur^{(\delta,\om,\gamma)}(\om z)  =\ep\cos(\om z + \gamma) + \mathcal{O}(\ep^2),
\end{align}
where $\gamma \in \R$ (or $\gamma \in \T := [0,2\pi]$), $\displaystyle{\vert 1 - \om^2\vert< \delta}$, $\displaystyle{\blue{\ep} =  \ep(\delta,\om) = \sqrt{\frac{4}{3}\left(\delta^2 -(1-\om^2)^2\right)} +\mathcal{O}\left(\delta^2 -(1-\om^2)^2\right)}$ for all  $\delta$  sufficiently small (cf. \cite[Chapter 17]{CE90} , \cite[\S 4]{Mielke}; see also Corollary \ref{rolls:reparametrization}). We contemplate the particular case of a  control parameter $\mu(\cdot)$ that varies spatially due to inhomogeneities in the media and given by 
\begin{align}\label{dq_ramp}
 \mu(x) = \left\{ \begin{array}{cl}
          1, & \text{for} \quad x\leq 0,\\
          -1, & \text{for} \quad x> 0.\\
         \end{array}\right.
\end{align}
From the mathematical point of view, besides breaking reflection and translation invariance symmetries, the parameter discontinuity in \eqref{dq_ramp} has consequences on the  energy potential associated with the equations, which jumps from bistable on $x\leq0$ to monostable on $x>0$, a scenario that describes what is known as \emph{directional quenching}; we  name  \textit{quenching-front} the boundary point  $x=0$  across which the system changes its stability. The jump in \eqref{dq_ramp} at $x=0$ closely emulates physically interesting   experiments where heterogeneities are introduced in the media aiming control of  micro-phase separation; such techniques have been applied in block copolymers \cite{Hashimoto}, dewetting and colloidal deposition \cite{berteloot2012evaporation}, patterning of surfaces \cite{langmuir};  in a similar spirit, directional quenching has also been studied in macro-phase separation models, cf. \cite{Foard-Wagner,Monteiro_Scheel}.

\begin{figure}[htb]
  \centering
 \includegraphics[width=1\textwidth]{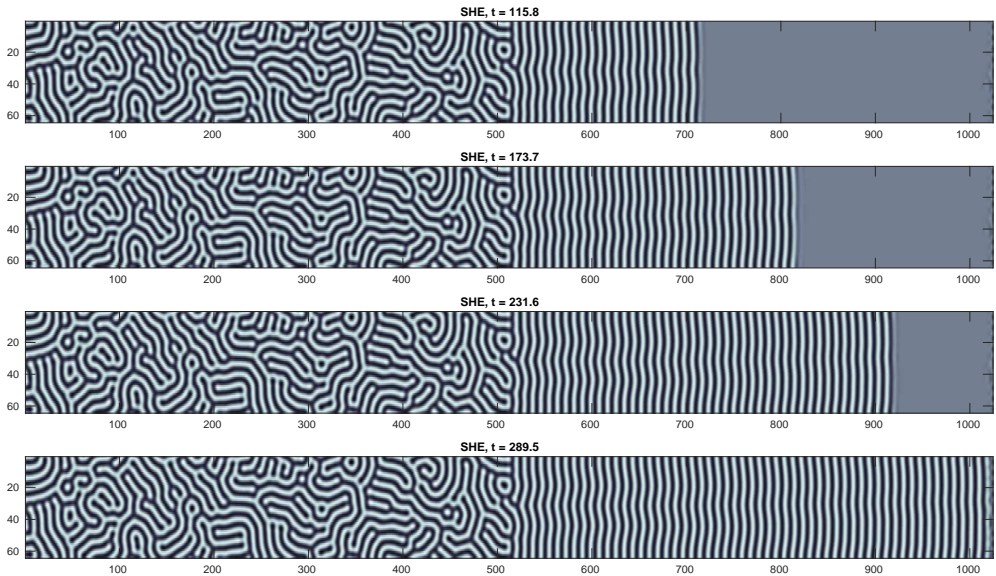}
 \caption{Snapshots of numerical simulations of SHE in 2D where, for a fixed $\delta$, a quenching-front propagates to the right. In the bulk (left), a fully phase separated state is seen, while in the wake of the quenching-front one sees rolls (right); ahead of the quenching-front we see a homogeneous state. Interestingly, the structure of which suggest the presence of zig-zag instabilities (see \S \ref{open_problems}). Similar results have been seen in numerical simulations of macro-phase separation models, cf. \cite{Foard-Wagner}. \label{fig:2D_simulations}}
\end{figure}

\begin{figure}[htb]
 \includegraphics[height=3cm]{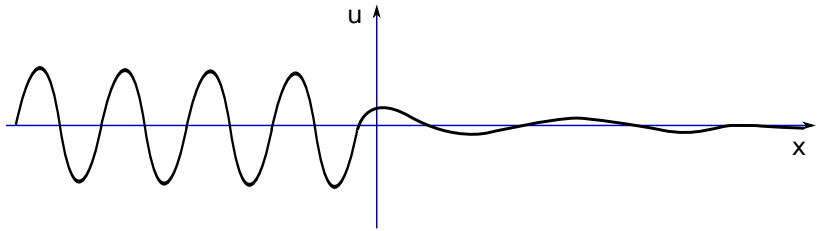}
 \caption{Sketch of heteroclinic connection between rolls and homogeneous states for 1D SHE; this sketch can be seen as a horizontal profile of the patterns seen in the wake of the quenching-front in Figure \ref{fig:2D_simulations}.\label{fig:1d_rolls}}
\end{figure}

Due to the dissimilarity in the media induced by the parameter $\mu(\cdot)$ it is reasonable to expect phase separated  states connected to homogeneous states, as one can see in numerical simulations; see Figure  \ref{fig:2D_simulations}. In the present study, we aim to understand a horizontal cross section of the pattern seen on the wake of the quenching-front of Figure \ref{fig:2D_simulations}; we shall constrain the analysis to 1D, therefore looking for point-to-periodic, time-independent solutions, as sketched in Figure \ref{fig:1d_rolls}. 

Before stepping into the discussion, we introduce a change variables $z \mapsto x/\om$ which, allied to the degree zero homogeneity property $\mu(x) = \mu\left(\om x\right)$, allow us to rewrite \eqref{SH-eq}  as 
\begin{equation}\label{SH-eq:after_change_of_variables}
 \partial_t u(x,t) =-(1 + \om^2\partial_x^2)^2 u(x,t)  + \delta^2\mu\left(x\right) u(x,t) -u^3(x,t);
\end{equation}
we shall look for a 1D, time-independent solution  $u(\cdot)$ satisfying
\begin{align}\label{asymptotics_properties}
 \lim_{x\to-\infty}\left\vert u(x) - \ur^{(\delta,\om,\gamma)}\left(x\right)\right\vert =0,\quad  \lim_{x\to+\infty}u(x)  =0, \quad 0 =-(1 + \om^2\partial_x^2)^2 u(x)  + \delta^2\mu\left(x\right) u(x) -u^3(x).
\end{align}

\noindent
Evidently, what we mean by a solution is still to be clarified, otherwise one could consider a solution like
\begin{align}\label{bad_solution}
 u(x) = \left\{\begin{array}{lll}
         \ur^{(\delta,\om,\gamma)}\left(x+x_0\right),& \text{on} &x\in (-\infty ,0);\\
          0 ,& \text{on}& x\in (0, +\infty),
        \end{array}\right.
\end{align}
which satisfies the Ordinary Differential Equation (ODE) in $\R\setminus\{0\}$, and can be adjusted for $x_0$ to be continuous on the whole line. 
A few remarks are in hand before we carefully characterize the type of solutions we are looking for.
 \begin{Remark}\label{Rmk:type_sol}
With regards to the problem \eqref{asymptotics_properties}, notice that
\begin{enumerate}[label=(\roman*), ref=\theTheorem(\roman*)]
\hangindent\leftmargin

\item\label{Rmk:type_sol-ode} Classical ODE theory shows that any solution to \eqref{asymptotics_properties} must be smooth on $\R\setminus\{0\}$. Furthermore, classical elliptic regularity theory shows that any weak solution must be locally $C^3(\R)$; the latter condition immediately rules out the ''solution'' \eqref{bad_solution}.
 
\item \label{rmk:admissibility}(Quenching-front constraint)  The parameter $\mu(\cdot)$  breaks the reversible symmetry ($x\mapsto -x$) of the ODE in \eqref{asymptotics_properties}. Nevertheless,  on the intervals $(-\infty,0)$ and $(0, +\infty)$ there exists associated Hamiltonians which are conserved quantities. Indeed, considering 
 $$-(1 + \om^2\partial_x^2)^2 \mathcal{V}(x) + \eta \mathcal{V}(x) - \mathcal{V}^3(x) = 0,$$
one has the associated Hamiltonian
\begin{equation*}
  \mathcal{H}[\mathcal{V};\eta] = -\frac{\om^4}{2}\left(\partial_x^2\mathcal{V}\right)^2 + \om^2(\partial_x\mathcal{V})^2 + \om^4\partial_x\mathcal{V} \partial_x^3\mathcal{V} + \frac{1}{4}\left(1 - \eta + \mathcal{V}^2\right)^2;
\end{equation*}
cf. \cite[\S 4.5.3, page 261]{Pismen}. Define the quantities 
$$\mathcal{H}^{(l)}[\mathcal{V}] := \mathcal{H}(\mathcal{V}; \delta^2),\quad \text{and}\quad  \mathcal{H}^{(r)}[\mathcal{V}] := \mathcal{H}(\mathcal{V}, -\delta^2).$$
Thus, for any solution $u(\cdot)\in \mathscr{C}^{(\infty)}(\R \setminus\{0\})\cup\mathscr{C}^{3}(\R)$ to the ODE \eqref{asymptotics_properties}  we must have
\begin{align}
  (-\infty,0)\ni x\mapsto \mathcal{H}^{(l)}[u(x)] \equiv \mathcal{H}^{(l)}[\ur^{(\delta,\om,\gamma)}(\cdot)], \qquad (0, +\infty)\ni x\mapsto \mathcal{H}^{(r)}[u](x) \equiv \mathcal{H}^{(r)}[0],
\end{align}
where $\mathcal{H}^{(l)}[\ur^{(\delta,\om,\gamma)}]$ and $\mathcal{H}^{(r)}[0]$ are constants that independent of $x$, but depend on $(\delta, \om,\gamma)$. Since $u(\cdot) \in \mathscr{C}^3(\R)$, the following quenching-front constraint must also be satisfied:
\begin{align}\label{phase_constraint}
  - \delta^2u^2\Big\vert_{x=0} - \delta^2 = \mathcal{H}^{(l)}[\ur^{(\delta,\om,\gamma)}(\cdot)] - \mathcal{H}^{(r)}[0].
\end{align}
In this case, we say that the parameters $(\delta,\om, \gamma)$ are admissible.

\item \label{Rmk:type_sol-changevariables} Notice that the change of variables $z\mapsto \frac{x}{\om}$ fix the period of the mapping  $x\mapsto \ur^{(\delta,\om,\gamma)}\left(x\right)$, which  is now $2\pi$-periodic, a fact that should not be overlooked due to its important consequences; see for instance Lemma \ref{Lem:to_be_not_to_be_continuous}. 
We remark that this change of variables do not eliminate the dependence of $\ur^{(\delta,\om,\gamma)}\left(x\right)$ on $\om$, because the amplitude of the rolls still depend nonlinearly on it; see \eqref{rolls}.
 \end{enumerate}
 \end{Remark}
 In this manner, we have the following
 \begin{Definition}[Heteroclinic orbits]\label{heterecolinic_def} We say that a solution  $(x,\delta,\om,\gamma) \mapsto \mathcal{U}^{(\delta,\om,\gamma)}(x)\in \mathscr{C}(\R\setminus\{0\})\cap \mathscr{C}^3(\R)$,  is a heteroclinic orbit to problem \eqref{asymptotics_properties} whenever the latter is satisfied and the quenching-front constraint  \eqref{phase_constraint} holds.  
 \end{Definition}
In order to find these heteroclinic orbits we explore a path that relies minimally on the theory of dynamical systems: we shall consider a functional-analytic based  perspective. Our first step benefits from  the  asymptotic conditions in \eqref{asymptotics_properties}, based on which we decompose the space where solutions are sought for: we shall refer to this step as a  \textit{far/near (spatial) decomposition}.  We shall consider a time-independent Ansatz of the form
\begin{align}\label{Ansatz}
 \mathcal{U}(x) = v(x) + \chi(\ep^{\beta} x)\ur^{(\delta,\om,\gamma)}\left(x\right),
\end{align}
in which a few unknowns are introduced: (i)  the parameter $\beta>0$, to be found later using matched asymptotics (see Proposition \ref{beta_is_1:sufficiency});  (ii) the function $\chi(\cdot) \geq 0$, a fixed smooth function, also to be chosen later and such that
\begin{align}\label{H1:partition_function}
 \lim_{x \to -\infty}\chi(x) = 1, \quad \lim_{x \to +\infty}\chi(x) = 0.
\end{align}
In fact, we shall see that  $\chi(\cdot)$ can be chosen to be a heteroclinic orbit satisfying  a second order ODE (see Lemma \ref{Lemma:almost_partition} and Section \ref{open_problems:invasion_fronts}). Note that the Ansatz \eqref{Ansatz} behaves as a  non-compact perturbation, i.e., as a  \textit{far field} perturbation. Moreover, with regards to \eqref{Ansatz}, whenever $(\delta, v(\cdot))=(0,0)$ we have  $\ur^{(0,\om,\gamma)}\equiv 0$ and, consequently,  $\mathcal{U}(\cdot) \equiv 0$, i.e., a trivial solution to \eqref{SH-eq}. Thus, one can consider the parameter  $\delta$ as a bifurcation parameter that, roughly speaking, \textit{turns on} a spatial-periodic perturbation in the far field as  $\delta>0$;  the function   $v(\cdot)$ plays the role of a corrector, and must  be chosen in an appropriate functional space.      

Plugging the Ansatz \eqref{Ansatz} into the ODE in \eqref{asymptotics_properties}, gives
\begin{equation}\label{SH-eq-bifurcation-1}
\begin{split}
 \mathscr{L}[v] := - \left( 1 + \om^2\partial_x^2\right)^2v & = \left\{-\left( \delta^2 \mu(x) - 3(\chi\ur^{(\delta,\om,\gamma)})^2\right)v\right\} +\{ 3(\chi \ur^{(\delta,\om,\gamma)})v^2\} +\{v^3\} \\
 & \quad + \left\{\chi( \chi^2 -1)(\ur^{(\delta,\om,\gamma)})^3 + [(1+ \om^2\partial_z^2)^2,\chi]\ur^{(\delta,\om,\gamma)} - \delta^2\chi (\mu -1)\ur^{(\delta,\om,\gamma)}\right\}.
 \end{split}
\end{equation}
Here we write  $\chi= \chi(\ep^{\beta} x)$ and  $ [\mathscr{A},\mathscr{B}] = \mathscr{A}\mathscr{B} - \mathscr{B}\mathscr{A}$ to denote the commutator of the operators $\mathscr{A}$ and $\mathscr{B}$. The functional space that contains $v(\cdot)$ consists of the domain of the operator $\mathscr{L}$, that is, $H^4(\R)$.

When plugged into the quenching-front constraint \eqref{phase_constraint}, we obtain
\begin{equation}\label{SH-eq-bifurcation-2}
   - \delta^2\left(v(0) + \chi(0)\ur^{(\delta,\om,\gamma)}(0)\right)^2 - \delta^2 = \mathcal{H}^{(l)}[\ur^{(\delta,\om,\gamma)}(\cdot)] - \mathcal{H}^{(r)}[0].
\end{equation}
We shall see in Lemma \ref{Lemma:regularity} that the constraint \eqref{SH-eq-bifurcation-2} is an essential ingredient behind a selection mechanism, where $\om$ is shown to be parametrized by by $(\delta,\gamma)$. 

\subsection{Properties of the model,  parameter choices  and main results}
\begin{SCfigure}[][h]
  \centering
 \includegraphics[width=0.4\textwidth]{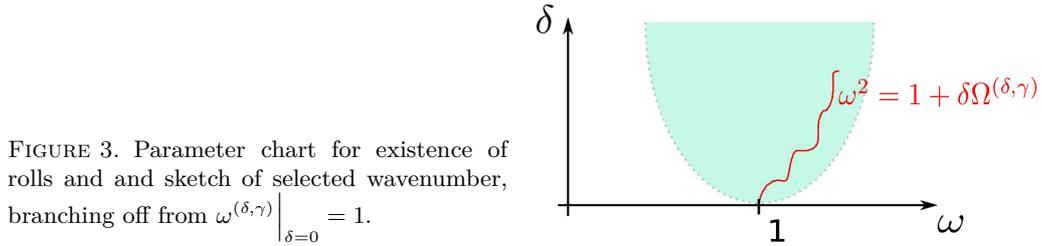}
 \caption{Parameter chart for existence of rolls and and sketch of selected  wavenumber, branching off from $\om^{(\delta,\gamma)}\Big|_{\delta =0} =1$.\label{Fig:selection_mechanism}}
\end{SCfigure}

Throughout this paper we make some a priori assumptions on  model's and Ansatz's parameters:  
\begin{enumerate}
 \item[(H1)] \label{H1}The partition function $\chi(\cdot)\in \mathscr{C}^{\infty}(\R;\R)$ satisfies \eqref{H1:partition_function}, and the convergence takes place in an exponential fashion, that is, there exists a $C_*>0$ and a $S_*>0$ for which
 \begin{equation*}
\left\vert\partial_x^{j}\left(\chi(x) -1\right)\right\vert \lesssim e^{-C_*\vert x \vert}, \quad \left\vert\partial_x^{j}\left(\chi(x) \right)\right\vert\lesssim e^{-C_*\vert x \vert}, \quad \forall j \in \{0,\ldots, 4\}, 
 \end{equation*}
 for all $x< - S_*$ and $x >S_*$, respectively.
  \item[(H2)] \label{H2} We shall further constrain $\om>0$ so  that $ \vert 1 - \om^2\vert \leq \frac{1}{3}\delta$, a range in which the existence of roll solutions is guaranteed (cf. \cite[Chapter 17]{CE90}; see also Figure \ref{Fig:selection_mechanism}). With this in mind, we introduce a parameter $\Omega$ (which, as $\ep$, is independet of $\gamma$) such that
  $$\om^2 = 1 + \delta \Omega ;\quad \text{for} \quad \Omega   \in\left(-\frac{1}{3},\frac{1}{3}\right).$$
  This readily implies that, for $\delta >0$,
  \begin{align*}
   \ep(\delta,\om) = \ep(\delta,\Omega,\gamma) =  \delta\sqrt{\frac{4}{3}\left(1 -\Omega^2\right)} +\mathcal{O}\left(\delta^2\left[1 -\Omega^2\right]\right).
  \end{align*}
Consequently, there exists\footnote{Interestingly, the role played by $\delta >0$ is equivalent to the one played by $\ep>0$. This assertion is proved later, in Corollary \ref{rolls:reparametrization}.} a $\delta_0>0$ such that $\displaystyle{\frac{1}{4}\delta \leq \ep \leq 4\delta}$, for $\delta \in [0,\delta_0)$. 
\end{enumerate}
\begin{Remark}[Parameter region blow-up]\label{rmk:parameter_blow_up}
The assumption  (\hyperref[H2]{H2}) is nothing but a parameter blow-up, for $(\delta,\Omega,\gamma)$ belongs to an open neighborhood of the zero in $\R^2$, whereas $(\delta,\om)$ does not (see Figure \ref{Fig:selection_mechanism}). This is a crucial ingredient in order to apply an Implicit Function Theorem based result; see \S \ref{section:wavenumber_selection}.
\end{Remark}

\begin{Remark}[On the role of the conserved quantities and Hamiltonian structure]
 We will see at the end of \S \ref{sec:reduced_equation+approx+solv} that,  once  $\chi(\cdot)$ is appropriately chosen, solutions to problem \eqref{SH-eq-bifurcation-1} are parametrized by  $(\delta,\om,\gamma)$. Note that under assumptions (\hyperref[H1]{H1})-(\hyperref[H2]{H2})  problem \eqref{SH-eq-bifurcation-1} is meaningful per se, regardless of the quenching-front constraint  \eqref{SH-eq-bifurcation-2}.  The  Hamiltonian structure  imposes  a severe restriction to the parameter region, yielding a \emph{selection mechanism} . In this fashion,  problem \eqref{SH-eq-bifurcation-1}, only in conjunction with the quenching-front constraint \eqref{SH-eq-bifurcation-2}, can be considered  an equivalent formulation to \eqref{asymptotics_properties}. 
\end{Remark}

\begin{figure}[htb]
  \centering
 \includegraphics[width=.48\textwidth]{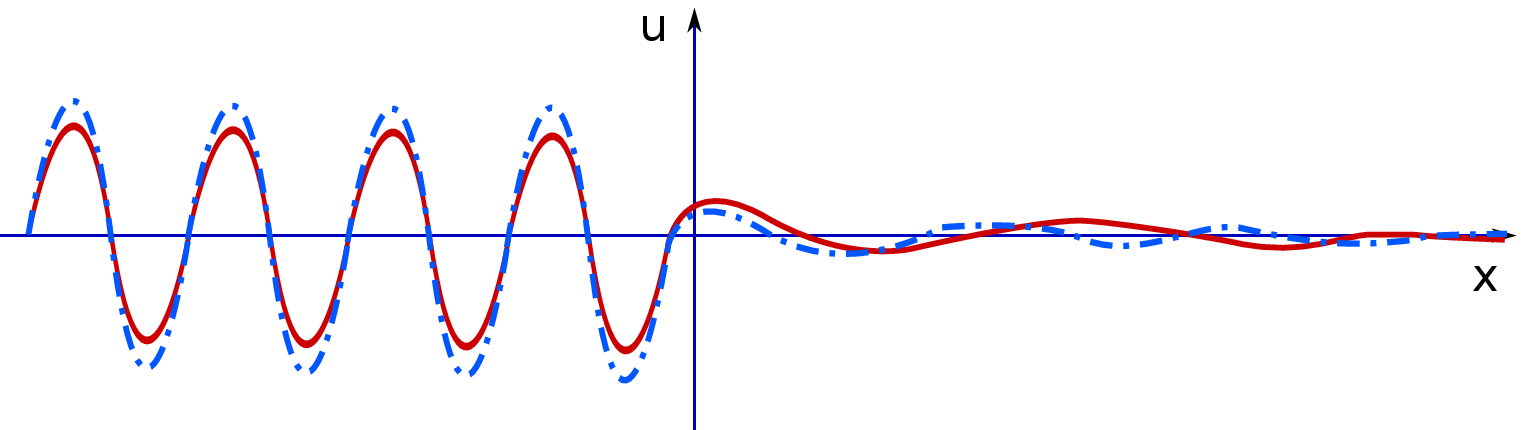}\hfill\includegraphics[width=.48\textwidth]{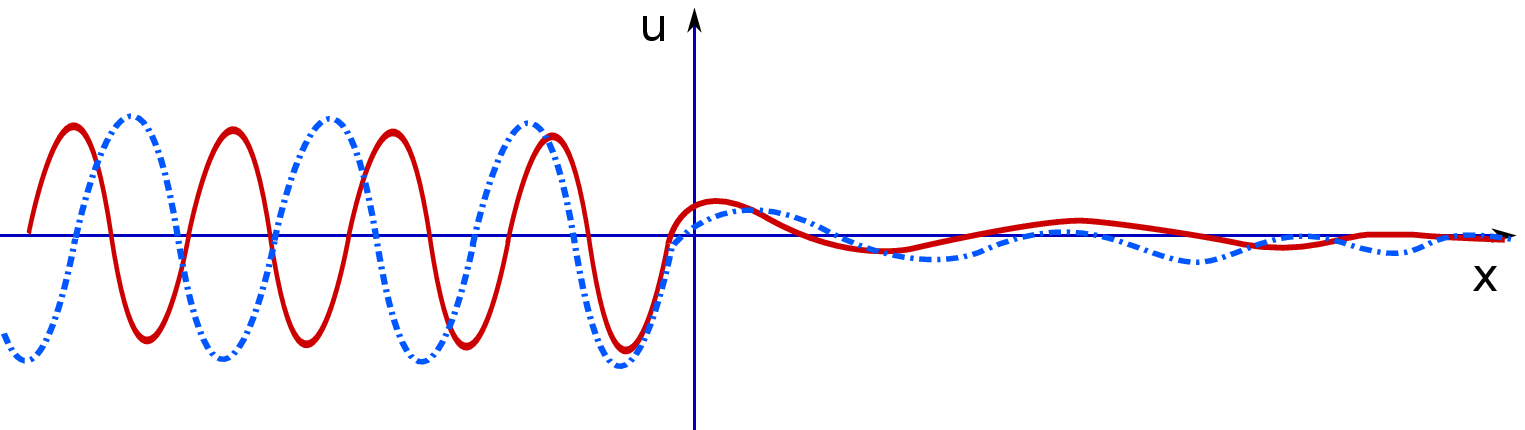}
 \caption{Sketch of bifurcating solutions illustrating Theorem \ref{Thm:main-regularity}: in the case $\om^{(\delta,\gamma)} \equiv 1$ (left) (resp., $\om^{(\delta,\gamma)} \not\equiv 1$ (right)) we have continuity (resp., discontinuity) in the sup norm.\label{fig:loss_continuity}}
\end{figure}

Our main result is the following:
\begin{Theorem}[Existence of 2-parameter family of heteroclinic connections]\label{Thm:main}
For any $0 < \tau <\frac{1}{16}$ fixed, there exists a $\delta_{**}=\delta_{**}(\tau)>0$, and a 2-parameter family of stationary solutions of the form
\begin{align}\label{Thm:main-shape}
[0,\delta_{**})\times \T \ni (\delta,\gamma) \mapsto  \mathcal{U}^{(\delta,\gamma)}(z) = v^{(\delta,\gamma)}(\om^{(\delta,\gamma)} z) + \chi(\ep \om^{(\delta,\gamma)} z)\ur^{(\delta, \om^{(\delta,\gamma)},\gamma)}(\om^{(\delta,\gamma)} z), \quad z\in \R, 
\end{align}
to \eqref{SH-eq} satisfying $\displaystyle{\mathcal{U}^{(\delta,\gamma)}(\cdot)\Big|_{\delta=0} \equiv 0}$ and with the following properties:
 \begin{enumerate}[label=(\roman*), ref=\theTheorem(\roman*)]
\hangindent\leftmargin

\item (Asymptotic properties and regularity) \label{Thm:main-regularity}The function $z\mapsto \mathcal{U}^{(\delta,\gamma)}(z)\in \mathscr{C}^{(\infty)}(\R\setminus\{0\})\cap \mathscr{C}^{(3)}(\R) $; furthemore, it satisfies
\begin{align*}
 \lim_{z\to - \infty }\left \vert \mathcal{U}^{(\delta,\gamma)}(z) - \ur^{(\delta, \om^{(\delta,\gamma)}, \gamma)}(z)\right \vert=0 ,\quad \lim_{z\to + \infty }\left \vert \mathcal{U}^{(\delta,\gamma)}(z)\right \vert =0,
\end{align*}
where $(\delta,\gamma) \mapsto \om^{(\delta,\gamma)}$ is a continuous mapping, as  defined below. Moreover, the mapping in \eqref{Thm:main-shape} is continuous in the sup norm if and only if $\om^{(\delta,\gamma)}\equiv 1$ (see Figure \ref{fig:loss_continuity}).
\item (Selection mechanism) \label{Thm:main-selection} The wavenumber $\om$ is continuously parametrized by $(\delta,\gamma)$ in  such a way that it satisfies (\hyperref[H2]{H2}), that is, 
\begin{align}
 [0,\delta_{**})\times \T\ni (\delta,\gamma) \mapsto \om^{(\delta,\gamma)}  = \sqrt{1 + \delta \Omega^{(\delta,\gamma)}};
\end{align}
the function $\displaystyle{(\delta,\gamma) \mapsto \Omega^{(\delta,\gamma)}\in \left(-\frac{1}{3},\frac{1}{3}\right)}$ is  a continuous mapping, 2$\pi$-periodic in $\gamma$, with the property that $\displaystyle{\Omega^{(\delta,\gamma)}\Big\vert_{\delta=0}=0}$.

 \item  (Envelope function) \label{Thm:main-envelope} The function $\chi(\cdot) \in \mathscr{C}^{(\infty)}(\R)$ is real valued and can be chosen in such a way that it satisfies the properties in (\hyperref[H1]{H1}). Moreover, $\chi(\cdot)$ is a fixed function  independent of $(\delta,\gamma)$  for  all $(\delta,\gamma) \in \left[0,\delta_{**}\right)\times \T.$
 \item (Fine structure of $v(\cdot)$)\label{Thm:main-fine_structure}
 We have that     $(\delta,\gamma) \mapsto v^{(\delta,\gamma)}(\cdot) \in H^4(\R)$, and this mapping can be decomposed as
\begin{align*}
 v^{(\delta,\gamma)}(z) = v\n^{(\delta,\gamma)}(z) + v\f^{(\delta,\gamma)}(z), \quad z\in \R,
\end{align*}
where both terms are 2$\pi$-periodic in $\gamma$ and, in particular,  $v\n^{(\delta,\gamma)}(\cdot)$ is a band-limited function that reads as
\begin{align*}
 v\n^{(\delta,\gamma)}(\cdot) =  \ep e^{+i  x} \gp(\ep x) + \ep e^{-i x} \gm(\ep x),
\end{align*}
where we recall from (ii) that $\ep = \ep(\delta, \om^{(\delta,\gamma)})$. Denoting the Fourier transform of $v\n^{(\delta,\gamma)}(\cdot)$ by $\widehat{v\n^{(\delta,\gamma)}}(\cdot)$,  it also holds that
\begin{align*} 
\mathrm{supp}\left(\widehat{v\n^{(\delta,\gamma)}}\right)\subset  \left\{ -1+ \ep^{\tau}\mathcal{B}\right\}\cup \left\{ 1+ \ep^{\tau}\mathcal{B}\right\},\qquad \mathcal{B} = \{x\in \R | \vert x\vert \leq 1\}.
\end{align*}
Moreover,  we can write $v\f(\cdot)$ as a function of $(v\n, \delta, \om),$ and the following bounds hold
\begin{align*}
\Vert v\n^{(\delta,\gamma)}(\cdot)\Vert_{H^{4}(\R)} = \mathcal{O}(\ep^{\frac{1}{2}}),\quad  \text{and} \quad \Vert v\f^{(\delta,\gamma)}(\cdot)\Vert_{H^{4}(\R)} = \mathcal{O}(\ep^{\frac{5}{2}- 2\tau}). 
\end{align*}
 \end{enumerate}
\end{Theorem}

\begin{Remark}[Reconciling the solution's periodicity with respect to $\gamma$ and its lack of translation invariance in space]  It is worth to highlight the interesting role played by the  parameter $\gamma$.  Notice that 
$$(x,\gamma) \mapsto \ur^{(\delta,\om,\gamma)}(x)$$
is $2\pi$-periodic in $x$, and also $2\pi$-periodic in $\gamma$. In the Ansatz \eqref{Ansatz} the far-field behavior as $x\to -\infty$  gets ``localized'' by the term $x\mapsto \chi(\ep^{\beta}x)$, where the latter function plays the role of a partition function subordinated to a neighborhood of $x =-\infty.$ Taking into account $\chi(\cdot)$'s properties (\hyperref[H1]{H1}), its non-periodicity  and independence of $\gamma$,  then the rightmost term in \eqref{Ansatz} 
\begin{align}\label{localization_ff_term}
(x,\gamma) \mapsto \chi(\ep^{\beta}x)\ur^{(\delta,\om,\gamma)}(x), 
\end{align}
is  $2\pi$-periodic in  $\gamma$, although its periodic $x$-symmetry gets broken; the latter is mostly due to the fact that the mapping only depends on $\gamma$ through $\ur^{(\delta,\om,\gamma)}$.  
A similar property can be seen in the corrector term $v^{(\delta,\gamma)}(\cdot)$ in \eqref{Ansatz}, that  is periodic with respect to $\gamma$ but doesn't have invariance with respect to translation in $x$ (due to the symmetry breaking term $\mu(\cdot)$).

\end{Remark}
\subsection{Outline of the paper and discussion} One of the main features of this article lies in the route paved to obtain our results:  we tackle and obtain Theorem \ref{Thm:main} from a functional analytic perspective, a fact that should be contrasted with the 1D dynamical systems methods used in  \cite{Weinburd}. There are advantages and drawbacks in our approach: on the positive side, one can foresee multi-dimensional extensions and applications to other contexts, like asymmetric-grain boundaries or equations with non-local terms (see discussion in \S \ref{open_problems:multiplier}); on the negative side it is clear that we obtained less information than \cite{Weinburd} with regards to the selection mechanism presented in Theorem \ref{Thm:main} (see also \S \ref{open_problems:dq_wavenumber_selec}). It should be said that the  analysis in \cite{Weinburd} is shorter than ours, although highly dependent on 1 dimensional dynamical systems theory; the functional-analytic construction of heteroclinic connections we advocate for seems more promising if one envisions  extensions to multi-D phenomena, like line defects or amplitude walls. Nevertheless, focusing on a mathematical problem that is now  rigorously understood give us a fair ground for comparison between the functional-analytic approach we advocate for and the dynamical-systems approach.

The heart of the paper lies in two decompositions:
\begin{enumerate}
 \item[(i)] the first has a spatial nature, intrinsically contained in \eqref{Ansatz}. It takes into account the far field spatial behavior of solution. We shall call it the \textit{far/near (spatial) decomposition}; 
 \item[(ii)] the second has a spectral nature, relying on the Fourier representation of the operator $\mathscr{L}[\cdot]$  in \eqref{SH-eq-bifurcation-1} as a Fourier multiplier. We shall call it the \textit{far/near (spectral) decomposition.}
\end{enumerate}
Both techniques had been seen separately in earlier works but, to the authors knowledge, have  never been used simultaneously. 
 
 The far/near (spatial) decomposition as done in \eqref{Ansatz} has been called by some authors  \textit{far field-core decomposition}:  it has been a  building block in the construction of multidimensional patterns in extended domains \cite{delPino}, in the study of perturbation effects on the far field of multidimensional patterns \cite{Monteiro_Scheel-contact_angle,Mon-horizontal}; in combination with  bifurcation techniques it is present  in the context of pattern formation as one can see in  the work of Scheel and collaborators (see for instance \cite{Monteiro_Scheel,Goh_Sch-pattern}, specially \cite{Lloyd-Sc} and \cite[\S 5]{Morrissey}).  Similar types of decompositions  have also been exploited in combination with homogenization techniques in the study and simulation of micro-structures and defects \cite{Lions_micro,Lions_micro-bigger}. Interestingly, to the authors knowledge, the combination of this type of decomposition with the introduction of a matching parameter ($\beta$ in \eqref{Ansatz})  has not been exploited before. 

The far/near (spectral) decomposition has a different origin, being deeply motivated and inspired by the work of \cite{Fefferman}. Therein, the authors consider perturbed Hamiltonian problems of the form
\begin{align*}
 \partial_x^2 \mathcal{U}(x) + V_{\text{per}}(x)\mathcal{U}(x) + \delta \kappa(\delta x)W_0(x) \mathcal{U}(x) = 0,\quad x\in \R,
\end{align*}
where $V_{\text{per}}(\cdot)$ is a known L-periodic potential, $W_0(\cdot)$ denotes an $L$-periodic perturbation and $\kappa(\cdot)$ is a smooth function with the property that $\displaystyle{ \lim_{x\to \pm \infty} \kappa(x) = \kappa_{\pm}}$. Using Bloch-Floquet theory the authors study  the spectrum of the periodic coefficients operator $\displaystyle{\partial_x^2 \mathcal{U}(x) + V_{\text{per}}(x)\mathcal{U}(x)}$ (i.e., the linear operator at $\delta =0$), whose multiplier/ eigenvalues are ``lifted'' to the case $\delta\neq 0$  using perturbation techniques. The next step, exploited in an ingenious way in \cite[\S 6]{Fefferman}, contains the seed of what we call \textit{far/near (spectral)} decomposition: it consists of decomposing the Brilloin zone for certain modes into a near and far region based on which one can derive a decomposition of the function $\mathcal{U}$ as
\begin{align*}
 \mathcal{U}(\cdot) = \mathcal{U}\n(\cdot) + \mathcal{U}\f(\cdot).
\end{align*}
The function $\mathcal{U}\n(\cdot)$ is a band limited function that parametrizes (or, say, \emph{dominates}) the far components $\mathcal{U}\f(\cdot)$; this step is proved using a Lyapunov-Schmidt reduction.  Surprisingly, as a stationary counterpart, the far/near (spectral) decomposition has similar features to  active/passive modes discussed in \cite[\S 1]{Newell-order_parameters}, where wavenumbers are divided in a group that saturates the asymptotic dynamic behavior, while passive modes are damped; when focused on stability issues, finite wavenumber instability studies dates from even earlier, cf. \cite{Newell}.

In spite of the nonlinear nature of \eqref{SH-eq-bifurcation-1}, our problem can be expressed as a far field perturbation as in \cite{Fefferman}, although a direct use of their results stumbles upon the nonlinear nature of SHE, and that  brings about issues not seen in their work. In dealing with these matters  our analysis is closer to the works of Schneider \cite{Schneider,Schneider:error,Schneider:validity}.  Schneider implements  a similar far/near (spectral) decomposition, using a frequency localization argument that he named  \textit{mode filters}. Although his aforementioned papers mostly focus on long time dynamic behavior of initial value problems, we can somehow say that the \textit{far/near (spectral) decomposition}  is essentially a stationary counterpart to the mode filters introduced therein.

We emphasize that no comparison principles are used neither in our derivation nor in \cite{Fefferman}. Even though  in the 4th order problem as we study here this comment seems apparently out of context, it should not be overlooked, in particular with regard to the second order problem studied in \cite{Fefferman},  where elliptic theory  plays an important role. Overall, another extra remark is worth of being made: neither the approach in \cite{Fefferman} nor ours rely on Center Manifold theory\footnote{For a thorough careful exposition of Center Manifold Theory, see \cite{Iooss}.  A nice discussion more tailored to the topics we are discussing, as lack of spectral gap and finite band instabilities, see   \cite[\S 1]{Newell-order_parameters}. }.

We now give a  brief outline of the paper: as we designed a far/near (spatial) decomposition in \eqref{Ansatz}, two new unknowns are introduced: a constant $\beta$ and a function $\chi(\cdot)$, both to be fully characterized  in Proposition \ref{beta_is_1:sufficiency} and Lemma \ref{Lemma:almost_partition}, respectively. In \S \ref{section:far_near} we  decompose of the Fourier wavenumbers using the aforementioned far/near (spectral) decomposition, which allows for the corrector $v(\cdot)$ to be written as
\begin{align*}
v(\cdot) =v\n(\cdot) + v\f(\cdot).
\end{align*}
Thus, one can define a coupled system of nonlinear equations for $v\n(\cdot)$ and $v\f(\cdot)$. This sets the ground for a  Lyapunov-Schmidt reduction in \S \ref{section:Lyap_Schm} where we show that under the appropriate conditions  $v\f(\cdot)$ is parametrized by  $(v\n(\cdot), \delta, \Omega)$; see Proposition \ref{Prop:near_far_domination}.
Afterwards, in \S \ref{section:blow_up}  we focus on $v\n(\cdot)$, which  satisfies
\begin{align*}
\mathrm{supp}\left(\hv\n\right)\subset  \left\{ -1+ \ep^{\tau}\mathcal{B}\right\}\cup \left\{ 1+ \ep^{\tau}\mathcal{B}\right\},\qquad \mathcal{B} = \{x\in \R | \vert x\vert \leq 1\},
\end{align*}
where $\ep, \tau>0$,  that is, the support of $\left(\hv\n\right)$ is contained in disjoint  intervals that shrink to points as $\ep\downarrow 0$, which is exactly the regime we aim to investigate. We overcome this issue by desingularizing the limit,  following the approach of \cite[\S 6.4]{Fefferman}.  In our case the consequences of this desingularization are quite interesting, giving a surprising interpretation of the near component $v\n(\cdot)$  as
$$v\n(x) = \ep^{\beta}e^{-ix}\gm(\ep^{\beta}x)+\ep^{\beta}e^{ix}\gp(\ep^{\beta}x), \qquad g_{\pm 1}(\cdot) \in H^2(\R);$$
a representation that is related  to the initial steps of the   Ginzburg-Landau formalism, commonly seem in modulation theory (cf.  \cite{Harten}, \cite[\S 3]{Schneider}; see also Proposition \ref{Prop:recentered}) and Weakly Nonlinear Theory (cf. \cite{Newell, Newell-order_parameters}).

In \S \ref{section:reduced_eq_simplifications} we use deep results from \cite[\S 6]{Fefferman} to better understand the periodic structure of the far field (``fast scale") and its interaction with the ``slow scale" structure of the correctors;  this allows for crucial simplifications in the equations. In the end we use matched asymptotic to  shown that $\beta =1$. 

In \S \ref{sec:reduced_equation+approx+solv} we finally write the reduced equation  \eqref{second_main_eqt:simplified} that contains the dominant features of the problem: a reduced nonlinear equation of the form 
\begin{align}\label{resume-sec6}
 \widetilde{\mathcal{R}^{(\delta,\Omega,\gamma)}}[\widehat{\gm}, \widehat{\gp}] - \mathbbm{1}_{\left\{ \ep^{\tau-1}\mathcal{B}\right\}}(\xi)h_{*}(\xi) = \widetilde{\mathcal{Q}^{( \delta,\Omega,\gamma;\pm)}}[\widehat{\gm}, \widehat{\gp}],
\end{align}
where the right hand side is a nonlinear function, sufficiently small in the appropriate sense; cf. \eqref{small_nonlinear_terms}. In the above equation, the mapping  $\displaystyle{(\widehat{\gm}, \widehat{\gp})\mapsto \widetilde{\mathcal{R}^{(\delta,\Omega,\gamma)}}[\widehat{\gm}, \widehat{\gp}]}$ is defined from $H^2(\R)\times H^2(\R)$ to $L^2(\R)\times L^2(\R)$, and the quantity $\displaystyle{\mathbbm{1}_{\ep^{1-\tau}\mathcal{B}}(\cdot) h_{*}(\cdot)}$ is a localized term in  $L^2(\R)\times L^2(\R)$ that can be made as small as we want upon choosing $\chi(\cdot)$ ``nicely''; see Proposition \ref{Lemma:almost_partition}.  In order to solve the nonlinear problem \eqref{resume-sec6}  we find an approximate inverse to $\widetilde{\mathcal{R}^{(\delta,\Omega,\gamma)}}[\widehat{\gm}, \widehat{\gp}]$  which, roughly speaking, corresponds to the formal limit of this same operator obtained as $\delta \downarrow 0$. This analysis culminates in another application of the Contraction Mapping Theorem (Proposition \ref{Prop:reduced}), showing that problem \eqref{resume-sec6}  has a family of solutions
\begin{align}
(0,\delta_{*})\times \left(-\delta_*,\delta_*\right)\times \T\ni (\delta,\Omega,\gamma) \mapsto \left(\gm^{(\delta,\Omega,\gamma)}(\cdot),\gp^{(\delta,\Omega,\gamma)}(\cdot)\right), 
\end{align}
where $g_{\pm 1}^{(\delta,\Omega,\gamma)}(\cdot) \in H^2(\R)$, that are band-limited, that is, $ \mathrm{supp}\left(\widehat{g_{\pm 1}^{(\delta,\Omega,\gamma)}}\right) \subset \ep^{\tau -1}\mathcal{B},$ with $\mathcal{B} = \{x\in \R | \vert x\vert \leq 1\}.$

At this point we reach  \S\ref{section:wavenumber_selection}, in which the Hamiltonian structure of the equations \eqref{SH-eq} is exploited and where Theorem \ref{Thm:main} is finally proved. An important step in its derivation comes from the study of the quenching-front constraint \eqref{SH-eq-bifurcation-2} (see  Lemma \ref{Lemma:wavenumber_selec}), where we  show that in fact, only two parameters are necessary in the parametrization of the heteroclinic orbits to \eqref{SH-eq}. Namely, for a $0<\delta_{**}<\delta_{*},$ we show that 
\begin{align*}
\left[0,\delta_{**}\right)\times \T \ni  (\delta,\gamma) \mapsto \om^{(\delta,\gamma)} = \sqrt{1 + \delta \Omega^{(\delta,\gamma)}}, 
\end{align*}
where the mapping $(\delta,\gamma) \mapsto \Omega^{(\delta,\gamma)}$ is  continuous and satisfy $\displaystyle{\Omega^{(\delta,\gamma)}\Big\vert_{\delta=0} =0}$;  this  results consists of a \emph{selection mechanism}, showing that  the wavenumber $\om$ gets parametrized by $(\delta,\gamma)$; see Figure \ref{Fig:selection_mechanism}. In conjunction with the results of   \S \ref{sec:reduced_equation+approx+solv} it implies that we can find a solution  $(\delta,\gamma) \mapsto v\n^{(\delta,\gamma)}(\cdot)$ to \eqref{SH-eq-bifurcation-1} and \eqref{SH-eq-bifurcation-2}. Undoing the change of variables $x \mapsto \om^{(\delta,\gamma)}z$ and plugging this solution in the Ansatz \eqref{Ansatz} we obtain a stationary solution 
\begin{align*}
[0,\delta_{**})\times \T\ni (\delta,\gamma) \mapsto  \mathcal{U}^{(\delta,\gamma)}(x) = v^{(\delta,\gamma)}(\om^{(\delta,\gamma)} x) + \chi(\ep \om^{(\delta,\gamma)} x)\ur^{(\delta, \om^{(\delta,\gamma)},\gamma)}(\om^{(\delta,\gamma)} x), \quad x\in \R, 
\end{align*}
as described in Theorem \ref{Thm:main}. Further properties of this mapping are studied in  Lemma \ref{Lem:to_be_not_to_be_continuous}.

In a quick summary, the construction of the pattern goes along the steps below:
\begin{align*}
 \begin{array}{c}
\text{Bifurcation equation setup: ``far/near'' spatial decomposition: Equation \eqref{SH-eq-bifurcation-1} }\\
 \big\Downarrow\\
 \text{``Far/near'' spectral decomposition: \S \ref{section:far_near}}\\
 \big\Downarrow\\
\text{Enslaving of far components by near components: \S \ref{section:Lyap_Schm}}\\
\big\Downarrow\\
\text{Blow up of the Fourier parameter, nonlinear interaction and approximation Lemmas: \S \ref{section:blow_up}}\\
\big\Downarrow\\
\text{Simplifications using a lemma of Fefferman, Thorpe and Weinstein, and matched asymptotics: \S \ref{section:reduced_eq_simplifications}}\\
\big\Downarrow\\
\text{Approximation and solvability of the reduced equation: \S \ref{sec:reduced_equation+approx+solv}}\\
\big\Downarrow\\
\text{Wavenumber selection --  proof of Theorem \ref{Thm:main}:  \S \ref{section:wavenumber_selection} }
 \end{array}
\end{align*}
Throughout the paper we give many different, but equivalent, formulations of the initial problem: in  \S\ref{section:far_near}
 problem \eqref{asymptotics_properties} is reduced to \eqref{mode_system_cutoff-a1} and \eqref{mode_system_cutoff-a2} after a far/near spectral decomposition. A Lyapunov-Schmidt reduction is then applied in \S \ref{section:Lyap_Schm}, resulting in  \eqref{second_main_eqt}.  After rewriting the latter equation with respect to the order (in $\ep$) of its nonlinearities, we reach \eqref{second_main_eqt:rewrite} in the beginning of \S 4. We devote \S 5 to the study of nonlinear interaction terms. Therein, further  simplification using ideas in \cite{Fefferman} culminate in  the reduced equation \eqref{second_main_eqt:simplified};  in passing, using matched asymptotics, we also conclude that $\beta=1$. Finally, in \S7 our Ansatz \eqref{Ansatz} is rewritten in new parameters as  \eqref{Ansatz-rewritten}; in this new formulation, it is used to  solve  \eqref{SH-eq-bifurcation-2}.
 
  The reader can track down the sequence of derivations in the next diagram:
\begin{displaymath}
    \xymatrix{
  &                               & \eqref{mode_system_cutoff-a1}  \ar[dr] & & & &&\\
\eqref{asymptotics_properties}\& \eqref{Ansatz}\ar[r]\ar[ddr]  &\eqref{SH-eq-bifurcation-1} \ar[ur]\ar[dr] &   &   \eqref{second_main_eqt}\ar[r] & \eqref{second_main_eqt:rewrite} \ar[r]& \eqref{second_main_eqt:simplified}\ar[r]&\eqref{Ansatz-rewritten}\ar[r]&\text{Theorem} \ref{Thm:main}\\
                               & & \eqref{mode_system_cutoff-a2}\ar[ur] &   &  & &&\\
  &\eqref{SH-eq-bifurcation-2}\ar[uurrrrrr]&  &   &  & &&} 
\end{displaymath}
Once the main result is proved, we discuss many open problems in  \S \ref{open_problems}, where some of the techniques we use are compared to previous methods and also put in a broader context.  An appendix  containing some important calculations  closes the paper.

\subsection{Notation and a functional-analytic settings}\label{notation}
Throughout this work we define the Fourier transform (resp., inverse Fourier transform) of  a function $f(\cdot) \in L^2(\R)$ (resp, $\widehat{f}(\cdot) \in L^2(\R)$) by 
\begin{align*}
\widehat{f}(\xi) =\mathscr{F}[f](\xi):=  \int_{\R}f(x)e^{-i  x \xi}\mathrm{d} x, \quad \left(\text{resp., } f(x) = \mathscr{F}^{-1}[\widehat{f}](x) = \frac{1}{2\pi}\int_{\R}\widehat{f}(\xi) e^{i \xi x}\mathrm{d} \xi\right). 
\end{align*}
A few properties of the  Fourier transform are used, in particular
 \begin{subequations}
 \begin{align}
   \mathscr{F}\left[f\right](\ep \xi) &= \mathscr{F}\left[\frac{1}{\ep}f\left(\frac{\cdot}{\ep}\right)\right](\xi),\label{Fourier_properties:dilation}\\
    \mathscr{F}[f](\alpha + \xi) &= \mathscr{F}[e^{- i \alpha (\cdot)}f(\cdot)]( \xi)\label{Fourier_properties:translation},\\
    \Vert \mathscr{F}[f]\Vert_{L^{\infty}(\R)} &\leq \Vert f\Vert_{L^{1}(\R)}\label{Fourier_properties:l_infty}.
\end{align}
 \end{subequations}
The pairing in Sobolev spaces $H^s(\R)$, $s\geq 0$, is defined as
\begin{align*}
 \langle f, g \rangle_{H^s(\R)} = \int_{\R} (1 + \vert \eta \vert^2)^s \widehat{f}(\eta)\overline{\widehat{g}(\eta)}\mathrm{d} \eta,
\end{align*}
where $\overline{(\cdot)}$ denoted  complex conjugation in $\mathbb{C}$; thanks to Plancherel Theorem (cf. \cite[\S 3, Proposition 3.2]{Taylor_I}), we have that $\Vert f\Vert_{L^2(\R)}= \frac{1}{\sqrt{2\pi}}\Vert f\Vert_{H^0(\R)}.$ We make repeated use of the following Sobolev embedding $H^1(\R) \hookrightarrow L^{\infty}(\R)$:
\begin{align}\label{Sobolev_embedding}
\Vert v \Vert_{L^{\infty}(\R)} \lesssim \int_{\R}\vert \hv(\xi)\vert \mathrm{d}\xi \leq \Vert v\Vert_{H^1(\R)}\sqrt{\int_{\R}\frac{1}{(1 + \xi^2)}\mathrm{d}\xi} \lesssim \Vert v\Vert_{H^1(\R)},
\end{align}
and of the Sobolev embedding $H^s(\R) \hookrightarrow \mathscr{C}^{(s-1)}(\R)$ (cf. \cite[Theorem 8.2]{Brezis}). The unit ball in $\R$ is denoted by $\mathcal{B} = \{x\in \R | \vert x\vert \leq 1\}$, while translated balls with radius $\varrho$ and centered at a point  $\alpha$ is written $\alpha + \varrho \mathcal{B} = \{x\in \R | \vert x - \alpha\vert  \leq \varrho\}.$ We also write $\T := [0,2\pi].$

The characteristic function of a Lebesgue measurable set $A$ is written $\mathbbm{1}_A(\tkp)= \mathbbm{1}_{\left\{\tkp \in A\right\}}(\tkp)$, where $\mathbbm{1}_A(x)= 1$, whenever $ x\in A$, and $\mathbbm{1}_A(x)= 0$ whenever $ x\not \in A.$ The support of a Lebesgue measurable function $f(\cdot)$ is denoted by $\mathrm{supp}(f)$.

We write $ H_{\text{near},\eta}^{s}(\R)\subset H^s(\R)$ to refer to  the space of band limited functions with Fourier transform supported in $\eta\mathcal{B}$,that is
\begin{align*}
 H_{\text{near},\eta}^{s}(\R):= \{g(\cdot)  \in H^s(\R)| \mathrm{supp}\left(\widehat{g}\right) \subset \eta \mathcal{B}\}.
\end{align*}
Given Banach spaces $X$ and $Y$, an unbounded operator $\mathscr{H}: X \to Y$ will have its domain written $\mathcal{D}\left(\mathscr{H}\right)$. Thus, given $\mathscr{H}:\mathcal{D}\left(\mathscr{H}\right) \subset X \to Y$, we write $\mathrm{Ker}\left(\mathscr{H}\right) := \{v \in X| \mathscr{H}v =0\}$. 

Throughout the paper we employ the usual ``little and big O'' convention: we write  $P = P(\ep) = \mathcal{O}(\ep^k)$ for some $k \in \N$ whenever there exists a constant $C>0$ such that $\vert P(\ep)\vert \leq C\vert\ep \vert^k$ as $\ep \to 0$. Similarly, we say that $p = p(\ep) = o(\ep^k)$ holds for some $k \in \N$ whenever, for any given  $C>0$, we have that  $\vert p(\ep)\vert \leq C\vert\ep \vert^k$ for all  $\ep$ sufficiently small; for instance, one says that $p(\ep) = o(1)$ whenever $\displaystyle{\lim_{\ep \to 0}p(\ep) = 0}$.

\begin{Remark}[Embedding for band limited functions]\label{embedding_band_limited}
As pointed out in \cite[Remark 2.2]{Fefferman}, whenever $f\in L^2(\R)$ is band limited equation then $f(\cdot)\in H^s(\R)$ for all $s \geq 1$ and 
\begin{align*}
\Vert f\Vert_{H^s(\R)} \lesssim \Vert f\Vert_{L^2(\R)}.
\end{align*}
In particular, $f(\cdot) \in \mathscr{C}^{\infty}(\R;\R)$, thanks  to the Sobolev Embedding Lemma (cf. \cite[\S 4-Corollary 1.4]{Taylor_I}). 
\end{Remark}

\subsection*{Acknowledgments}  R.M would like to thank the encouragement and feedback given by P. Sternberg (Indiana University), A. Scheel, T. Tao, and J. Weinburd (Univ. of  Minnesota), A. Nachbin, A. Maliebaev,  D. Marchesin (Fluids group at IMPA), A. Pastor and M. Martins (PDE group at UNICAMP), Y. Nishiura (AIMR/MathAM-OIL, Sendai), who kindly listened to his explanations of preliminary (and at the time, somewhat obscure) aspects of this work.  

``This is a post-peer-review, pre-copyedit version of an article published in \textit{Archive for Rational Mechanics and Analysis}. The final authenticated version is available online at: \href{http://dx.doi.org/10.1007/s00205-019-01427-z}{http://dx.doi.org/10.1007/s00205-019-01427-z}”. 

\section{The far/near (spectral) decomposition: the role of multipliers}\label{section:far_near}
To begin with, we shall represent equation \eqref{SH-eq-bifurcation-1} in a more concise form, 
\begin{equation*}
\begin{split}
 \mathscr{L}[v] := - \left( 1 + \om^2\partial_x^2\right)^2v =  \sum_{j=1}^4   \mathscr{N}^{(j)}[v, \ur^{(\delta,\om,\gamma)},\delta],
 \end{split}
\end{equation*}
where  $v(\cdot) \in \mathcal{D}\left(\mathscr{L}\right) = H^4(\R),$ the domain of the operator $\mathscr{L}$. The nonlinearities are
\begin{equation}\label{SH-eq-nonlinear_terms}
 \begin{split}  
 \mathscr{N}^{(1)}[v, \ur^{(\delta,\om,\gamma)},\delta](x) & = -\left( \delta^2 \mu(x) - 3(\chi(\ep^{\beta} x)\ur^{(\delta,\om,\gamma)}(x))^2\right)v(x), \\
 \mathscr{N}^{(2)}[v, \ur^{(\delta,\om,\gamma)},\delta](x) &=   3(\chi(\ep^{\beta} x) \ur^{(\delta,\om,\gamma)})v^2(x),\\ 
 \mathscr{N}^{(3)}[v, \ur^{(\delta,\om,\gamma)},\delta](x)&= v^3(x),\\
 \mathscr{N}^{(4)}[v, \ur^{(\delta,\om,\gamma)},\delta](x)&= \chi(\ep^{\beta} x)( \chi(\ep^{\beta} x)^2 -1)(\ur^{(\delta,\om,\gamma)}(x))^3 \\
 &\quad  + [(1+ \om^2\partial_x^2)^2,\chi(\ep^{\beta} \cdot)](x)\ur^{(\delta,\om,\gamma)}(x) - \delta^2\chi(\ep^{\beta} x) (\mu(x) -1)\ur^{(\delta,\om,\gamma)}(x).
 \end{split}
\end{equation}
In Fourier space, the operator $\mathscr{L}$ admits a multiplier representation,
\begin{equation}\label{multiplier}
\mathscr{F}\left[\mathscr{L}[v]\right](\kp)  = m(\kp; \mathscr{L})\mathscr{F}[v](\kp) = - (1 - \om^2 \kp^2)^2\hv(\kp), \quad v(\cdot)\in H^4(\R). 
\end{equation}
Taking into account the properties of  the mapping $\kp \mapsto \frac{1}{m(\kp; \mathscr{L})}$, we decompose the frequency space in two disjoint sets:
\begin{enumerate}
 \item[(i)] Near frequency region: this is the part around the zeros of $\kp \mapsto m(\kp; \mathscr{L})$ (i.e., $\kp = \pm\frac{1}{\om}$), where the mapping $\kp \mapsto \frac{1}{m(\kp; \mathscr{L})}$ has a ``bad behavior'';
 \item[(ii)] Far frequency region: the complement to the near frequency region, where we have a better behavior for $\kp \mapsto \frac{1}{m(\kp; \mathscr{L})}$.
\end{enumerate}
We do the following splitting:
\begin{equation}\label{nr_fr_frequency_regions}
 \begin{split}
\text{Near frequencies} &= \left\{- 1 + \ep^{\tau}\mathcal{B}\right\}\cup\left\{ 1+ \ep^{\tau}\mathcal{B}\right\}, \quad \text{Far frequencies} = \R \setminus\{\text{near frequencies}\},
 \end{split}
\end{equation}
where we introduce a new parameter  $\tau$ whose choice, or better saying its order, will only be defined throughout our analysis. At the present we choose its sign, arguing as follows: we shall define in  \S\ref{Splitting:sec}  spaces of functions $X_{\text{near}, \ep^{\tau}}^{s} \subset H^s(\R)$ (resp., $X_{\text{far}, \ep^{\tau}}^{s}\subset H^s(\R)$), whose elements have Fourier transform supported in the near frequency region (resp. far frequency region). We would like to work with band-limited functions, therefore,  in order to have at least one of this sets uniformly bounded in the regime  $\ep \downarrow 0$  (the one relevant to us) we must take $\tau \geq 0$. As we shall observe in Proposition \ref{Prop:irrelevant}, the endpoint $\tau =0$ must also be excluded from the analysis because in this case we do not have enough control of nonlinearities. Roughly speaking, when $\tau=0$ the corrector can be controlled in its amplitude but not its  bandwidth, and this is not enough to obtain a reduced equation approximation; as we show here, this reduction is possible when $\tau>0$, which turns out to be  the case  we  consider the most interesting.

Let's go back to the  multiplier structure \eqref{multiplier} and the decomposition  \eqref{nr_fr_frequency_regions}: as we can see, the linearized operator $v\mapsto \mathscr{L}[v]$ has continuum spectrum up to the imaginary axis in the complex plane. The decomposition \eqref{nr_fr_frequency_regions} splits the spectrum in modes that are close to the imaginary axis (near frequencies) and those that are far, or relatively far, from it (far frequencies). As we shall prove in \S \ref{section:Lyap_Schm}, the near frequencies are dominant in this problem. 
\begin{SCfigure}[][h]
  \centering
\includegraphics[height=2.5cm]{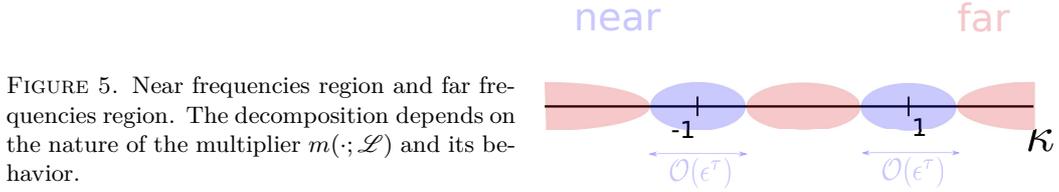}
 \caption{Near frequencies region and far frequencies region. The decomposition depends on the nature of the multiplier $m(\cdot;\mathscr{L})$ and its behavior.\label{fig:nr_fr}} 
\end{SCfigure}

We remark that the choice for a parametrization in terms of $\ep$ is somewhat a matter of convenience, and we do so because   nonlinearities scale in $\ep$. In principle, a similar analysis could be done in terms of $\delta$, thanks to the equivalence between them asserted in (\hyperref[H2]{H2}) (see also Corollary \ref{rolls:reparametrization}). Another important consequence of the next result is the fact that the lower bound constant  $C$ is independent of $\gamma$: this fact will be fundamental in deriving bounds that are uniform for all $\gamma \in \T$.

\begin{Lemma}[Multiplier behavior over frequency regions]\label{Lemma_multiplier}  The linear operator $v\mapsto \mathscr{L}[v]$ can be represented as a multiplier in Fourier space, that is, it holds the equivalence
\begin{align*}
  \mathscr{L}[v] = \mathscr{F}^{-1}\left[m(\cdot;\mathscr{L})\hv(\cdot)\right](x),  \quad \mathcal{D}\left(\mathscr{L}\right) = H^4(\R),
\end{align*}
where $m(\kp; \mathscr{L}) = -\left(1 - \om^2 \kp^2\right)^2.$ Whenever $\om$ satisfies (\hyperref[H2]{H2}), the zeros of $\kp \mapsto m(\kp;\mathscr{L})$, taken at $\kp =\pm \frac{1}{\om}$, remain in the \textit{Near frequency} region for all $\ep\geq0$ sufficiently small.  Furthermore, whenever $0 <\tau <1$ there exists a $C = C(\ep_0)>0$, that is, independent of $\gamma$,  such that 
\begin{align*}
 \vert m(\kp; \mathscr{L})\vert\Big|_{\{\kp \in \text{Far frequencies}\}} \geq C \ep^{2\tau}, \quad \forall \ep \in [0,\ep_0).
\end{align*}
\end{Lemma}
\begin{Proof}
 The first statement concerning the multiplier representation is  a consequence of standard Fourier Analysis.  Recall from (\hyperref[H2]{H2}) that $\om^{(\delta,\gamma)} \in \left(\sqrt{1-\frac{\delta}{3}}, \sqrt{1+\frac{\delta}{3}}\right)$. The proof of the case $\tau =0$ is direct, so we focus on the case $0 < \tau<1$. In the latter case, proving the inclusions $+ \frac{1}{\om} \in \text{Near frequency}$ is equivalent to verifying the inequalities
 \begin{align}\label{ineq_near_far}
  1 - \ep^{\tau} \leq 1 - \frac{\delta^{\tau}}{4^{\tau}} \leq\frac{1}{\sqrt{1+\frac{\delta}{3}}}\leq \frac{1}{\om} \leq \frac{1}{\sqrt{1-\frac{\delta}{3}}} \leq 1 + \frac{\delta^{\tau}}{4^{\tau}}\leq 1 + \ep^{\tau};
 \end{align}
 we omit the equivalent statement and proof of the case $- \frac{1}{\om} \in \text{Near frequency}$, which can be handled similarly. The first and the last inequalities are straightforward consequences of (\hyperref[H2]{H2}). For the other inequalities, we shall verify that 
\begin{align}\label{ineq_near_far-left_right}
 1 \leq \left(1 + \frac{\delta^{\tau}}{4^{\tau}}\right)^2\left(1 - \frac{\delta}{3}\right), \quad  \text{and} \quad  \left(1 - \frac{\delta^{\tau}}{4^{\tau}}\right)^2\left(1 + \frac{\delta}{3}\right)\leq 1.
\end{align}
 holds  whenever $0 < \tau <1$ and $\delta$ is sufficiently small: we first expand the inequality on the left hand side (resp. right hand side) and rearrange it, resulting in the  equivalent inequality,
\begin{align*}
 0  \leq \delta^{\tau}\left[\frac{2}{4^{\tau}} + \frac{\delta^{\tau}}{4^{2\tau}}- \frac{\delta^{1- \tau}}{3}\left( 1 + \frac{\delta^{\tau}}{4^{\tau}}\right)^2\right], \quad  \left(\text{resp.} \quad    \delta^{\tau}\left[- \frac{2}{4^{\tau}} + \frac{\delta^{\tau}}{4^{2\tau}} + \frac{\delta^{1- \tau}}{3}\left( 1 - \frac{\delta^{\tau}}{4^{\tau}}\right)^2\right]\leq 0 \right),
\end{align*}
 which holds for all $\delta \geq 0$ sufficiently small. This concludes the verification of \eqref{ineq_near_far-left_right} and, consequently, of $\eqref{ineq_near_far}$.
 
 We now prove the lower bound on the multiplier, which consists in finding
 $$\min_{\vert\xi  -1\vert\geq \ep^{\tau}}\{(1 - \om^2 \xi^2)^2\},\quad \text{for} \quad \om^{(\delta,\gamma)} \in \left(\sqrt{1-\frac{\delta}{3}}, \sqrt{1+\frac{\delta}{3}}\right).$$
 By symmetry, it suffices to consider $\xi \in \{ \xi \geq 0\} \cap \{\vert\xi  -1\vert\geq \ep^{\tau}\} $. In combination with standard critical point calculations, \eqref{ineq_near_far} shows that the minimum is attained when $\xi  \in \text{Near frequency region}$. Hence, one can reduce the minimization problem to evaluating $\xi$ in $\{1 \pm \ep^{\tau}\}$. Let's work with the case $\xi = 1 +\ep^{\tau}$; the other case is similar. Writing $\om^2 = 1 + \delta\Omega$, we are facing  the problem of minimizing 
 $$ \left(1 - (1 + \Omega \ep) \left( 1 + \ep^{\tau}\right)^2\right)^2,\quad \text{for} \quad \Omega \in \overline{\left(-\frac{1}{3},\frac{1}{3}\right)}=\left[-\frac{1}{3},\frac{1}{3}\right].$$
 It suffices to check the minimum at $\Omega =\pm \frac{1}{3}.$ At $\xi = 1+ \ep^{\tau}$ we have
 \begin{equation*}
 \begin{split}
  \left[1 - \left(1 + \frac{\ep}{3}\right) \left( 1 + \ep^{\tau}\right)^2\right]^2&=  \left[-\ep^{\tau} \left( 2 + \ep^{\tau}\right)- \frac{\ep}{3}\left( 1 + \ep^{\tau}\right)^2\right]^2  = \ep^{2\tau }\left[ \left( 2 + \ep^{\tau}\right)+ \frac{\ep^{1 -\tau}}{3}\left( 1 + \ep^{\tau}\right)^2\right]^2.
 \end{split}
 \end{equation*}
 Likewise,
 \begin{equation*}
 \begin{split}
  \left[1 - (1 - \frac{\ep}{3}) \left( 1 + \ep^{\tau}\right)^2\right]^2&=  \left[-\ep^{\tau} \left( 2 + \ep^{\tau}\right)+ \frac{\ep}{3}\left( 1 + \ep^{\tau}\right)^2\right]^2  = \ep^{2\tau }\left[ \left( 2 + \ep^{\tau}\right)- \frac{\ep^{1 -\tau}}{3}\left( 1 + \ep^{\tau}\right)^2\right]^2.
 \end{split}
 \end{equation*}
A similar analysis can be applied in the case  $\xi = 1- \ep^{\tau}$ and to the cases $\xi = -1 \pm\ep^{\tau}$, showing that  there exists a constant $C>0$ such that  $\vert m(\kp; \mathscr{L})\vert\Big|_{\{\kp \in \text{Far frequencies}\}} \geq C \ep^{2\tau}$, whenever $\ep \in [0,\ep_0)$. This finishes the proof.  \end{Proof}
\subsection{The multiplier structure and the far/near  decomposition}\label{Splitting:sec}
After  splitting the spectrum we define associated spaces that will be of fundamental importance; the idea is based on \cite[\S 6.3, Lemma 6.2]{Fefferman} (see also \cite[\S 2]{Schneider}). In the discussion below, we assume  that  $v(\cdot) \in H^4(\R)$. We define cut-offs in Fourier space, which parametrize the excised region around the zeros of the multiplier $m(\cdot; \mathscr{L})$:
%
\begin{align*}
 \hv(\cdot) \mapsto \widetilde{\mathbb{P}\n}[\hv](\kp) &=  \left(\mathbbm{1}_{\left\{ -1 +\ep^{\tau}\mathcal{B}\right\}}(\kp) + \mathbbm{1}_{\left\{ +1 +\ep^{\tau}\mathcal{B}\right\}}(\kp)\right)\hv(\kp) =: \hv\n(\kp);
\end{align*}
this operator naturally induces the following  projection in physical space,
\begin{align*}
 v(\cdot)\mapsto \mathbb{P}\n[v](x)= \mathscr{F}^{-1}\circ \widetilde{\mathbb{P}\n} \circ \mathscr{F}[v](x) =: v\n(x),
\end{align*}
or, in other words, $v\n(x) = \mathscr{F}^{-1}[\hv\n(\cdot)](x)$. Similarly, we define projections onto the  far frequencies space as
\begin{align*}
 \hv(\cdot) \mapsto \widetilde{\mathbb{P}\f}[\hv](\kp) &= \left(1 - \mathbbm{1}_{\left\{ -1 +\ep^{\tau}\mathcal{B}\right\}}(\kp) - \mathbbm{1}_{\left\{ +1 +\ep^{\tau}\mathcal{B}\right\}}(\kp)\right)\hv(\kp)=: \tv_{\text{far}}(\kp),
 \end{align*}
with associated  physical space projection,
\begin{align*}
v(\cdot)\mapsto  \mathbb{P}\f[v](x)= \mathscr{F}^{-1}\circ \widetilde{\mathbb{P}\f} \circ \mathscr{F}[v](x) =:v\f(x).
\end{align*}
It is clear that $\widetilde{\mathbb{P}\n}[\hv](\kp)$, $\widetilde{\mathbb{P}\f}[\hv](\kp)$ and associated physical space projections all depend on $(\ep,\tau)$; to keep our notation concise, we shall omit this dependence. Several properties of these mappings are readily available:
\begin{subequations}\label{nr_fr_decomposition:properties}
\begin{align}
 \hv\n(\kp) &= \widetilde{\mathbb{P}\n}\left[\hv\n\right](\kp) = \left(\mathbbm{1}_{\left\{ -1 +\ep^{\tau}\mathcal{B}\right\}}(\kp) + \mathbbm{1}_{\left\{ +1 +\ep^{\tau}\mathcal{B}\right\}}(\kp)\right)\hv\n(\kp),\\
 \hv\f(\kp) &= \widetilde{\mathbb{P}\f}\left[\hv\f\right](\kp) = \left(1 - \mathbbm{1}_{\left\{ -1 +\ep^{\tau}\mathcal{B}\right\}}(\kp) - \mathbbm{1}_{\left\{ +1 +\ep^{\tau}\mathcal{B}\right\}}(\kp)\right)\hv_{\text{far}}(\kp), \\
v(\cdot) &= v\n(\cdot) + v\f(\cdot),\\
\langle v\n,v\f \rangle_{H^s(\R)} &=0,\\
\Vert v\Vert_{H^s(\R)}^{2} &= \Vert v\n\Vert_{H^s(\R)}^{2}  +\Vert v\f\Vert_{H^s(\R)}^{2}. 
\end{align}
\end{subequations}
Thanks to these properties, we define
\begin{align}\label{spaces:X_n_X_f}
H^s(\R) = X_{\text{near}, \ep^{\tau}}^{s} \oplus X_{\text{far}, \ep^{\tau}}^{s}, \quad \text{where} \quad X_{\text{near}, \ep^{\tau}}^{s} : = \mathbb{P}\n\left(H^{s}(\R)\right)\quad  \text{and} \quad X_{\text{far}, \ep^{\tau}}^{s} : = \mathbb{P}\f\left(H^{s}(\R)\right).
\end{align}
As $X_{\text{near}, \ep^{\tau}}^{s}$,  $X_{\text{far}, \ep^{\tau}}^{s}$ are subspaces of $H^{s}(\R)$, we adopt the  induced norm they inherit.
\subsection{A Lyapunov-Schmidt reduction} As we have seen, the  near and far frequency regions are defined by taking into account the behavior of the multiplier $m(\cdot; \mathscr{L})$. In particular, the projections from $H^{s}(\R)$ onto $X_{\text{near}, \ep^{\tau}}^{s}$  and $X_{\text{far}, \ep^{\tau}}^{s}$ constructed in \S\ref{Splitting:sec}  are put in full use  to implement a Lyapunov-Schmidt reduction; our approach closely follows the ideas in \cite[\S 6]{Fefferman}. 

Let $\displaystyle{v(\cdot) = v\n(\cdot) + v\f(\cdot) \in X_{\text{near}, \ep^{\tau}}^{4}\oplus X_{\text{far}, \ep^{\tau}}^{4}}$. Thanks to \eqref{spaces:X_n_X_f}, we can rewrite \eqref{SH-eq-bifurcation-1}  in an equivalent system: the first equation results from the application of $\widetilde{\mathbb{P}\n}\circ\mathscr{F}$,
\begin{subequations}
\begin{equation}
 \begin{split}
- (1 - \om^2\kp^2)^2\hv\n(\kp)&= \widetilde{\mathbb{P}\n}\circ\mathscr{F}\left[\sum_{j=1}^4 \mathscr{N}^{(j)}[v\n + v\f, \ur^{(\delta,\om,\gamma)}, \ep]\right](\kp)\\
&= \left(\mathbbm{1}_{\left\{- 1 + \ep^{\tau}\mathcal{B}\right\}}(\kp) + \mathbbm{1}_{\left\{+1 + \ep^{\tau}\mathcal{B}\right\}}(\kp)\right) \mathscr{F}\left[\sum_{j=1}^4 \mathscr{N}^{(j)}[v\n + v\f, \ur^{(\delta,\om,\gamma)}, \ep]\right](\kp)\label{mode_system_cutoff-a1},
 \end{split}
\end{equation}
while the second equation, complimentary, is derived after an application of  $\widetilde{\mathbb{P}\f}\circ\mathscr{F}$ to \eqref{SH-eq-bifurcation-1}:
\begin{equation}
 \begin{split}
- (1 - \om^2\kp^2)^2 &\hv\f(\kappa) = \widetilde{\mathbb{P}\f}\circ\mathscr{F}\left[\sum_{j=1}^4 \mathscr{N}^{(j)}[v\n + v\f, \ur^{(\delta,\om,\gamma)}, \ep]\right](\kp)\\
&= \left(1 -\mathbbm{1}_{\left\{- 1 + \ep^{\tau}\mathcal{B}\right\}}(\kp) - \mathbbm{1}_{\left\{+1 + \ep^{\tau}\mathcal{B}\right\}}(\kp)\right)\mathscr{F}\left[\sum_{j=1}^4 \mathscr{N}^{(j)}[v\n + v\f, \ur^{(\delta,\om,\gamma)}, \ep]\right](\kp)\label{mode_system_cutoff-a2}.  
 \end{split}
\end{equation}
\end{subequations}
for all  $\kp\in \R$. Our goal is to show that the near field in fact dominates the far field components. With this in mind, we recall from \eqref{rolls} that $\ep(\cdot,\cdot)$ is a function of both $\delta$ and $\Omega$, invoking Lemma \ref{Lemma_multiplier} to then  rewrite \eqref{mode_system_cutoff-a2} as  a fixed point equation
\begin{align}\label{Fixed_point}
 v\f(x) = \mathscr{E}[v\n, v\f,  \delta,\Omega,\gamma](x),
\end{align}
where  $\mathscr{E}:X_{\text{near}, \ep^{\tau}}^{4}\times X_{\text{far}, \ep^{\tau}}^{4} \times (0,+\infty) \times \left(-\frac{1}{3},\frac{1}{3}\right)\to X_{\text{far}, \ep^{\tau}}^{4}$ reads as
\begin{equation*}
\begin{split}
\mathscr{E}[v\n, v\f,  \delta,\Omega,\gamma](x):= \mathscr{F}^{-1}\left[ -\left(\frac{1}{m(\cdot, \mathscr{L})}\right)^2\widetilde{\mathbb{P}\f}\circ\mathscr{F}\left[\sum_{j=1}^4 \mathscr{N}^{(j)}[v\n + v\f, \ur^{(\delta,\om,\gamma)}, \ep]\right](\cdot)\right](x).
\end{split}
\end{equation*}
%


\section{Near field components dominate far energy terms: a Lyapunov-Schmidt reduction.}\label{section:Lyap_Schm}
The main result of this section establishes the parametrization of the far energy terms $v\f(\cdot)$ by the near field $v\n(\cdot)$ and the other parameters $(\delta, \Omega)$. For the sake of notation, throughout this section we write  $\Vert \cdot\Vert_{H^s}$ and $\Vert \cdot\Vert_{L^p}$ to denote the  norms $\Vert \cdot\Vert_{H^s(\R)}$ and $\Vert \cdot\Vert_{L^p(\R)}$, respectively. 

\begin{Proposition}[Near field domination]\label{Prop:near_far_domination} Assume (\hyperref[H1]{H1})-(\hyperref[H2]{H2}). Let    $R$,  $\tau, \beta$ be  fixed  positive numbers, the latter two of which were introduced in \eqref{nr_fr_frequency_regions} and \eqref{Ansatz}, respectively. Assume the parameter constraints
\begin{align}\label{pmtr_range}
  0 <\tau <\frac{1}{16},\quad \text{and} \quad  1\leq \beta \leq 1+ \tau,
\end{align}
fixed ($\beta = 1$ and $\tau = \frac{1}{32}$ will do). Recall from  (\hyperref[H2]{H2}) that $\ep(\delta,\Omega,\gamma)$ is a smooth function of its arguments. Then, there exists an $\delta_*>0$ such that for all $\delta \in (0,\delta_*)$   the system \eqref{Fixed_point} has a fixed point, namely, there exists a mapping
 \begin{align}\label{Prop:nr_eng_dmn:enslaving_mappig}
 (v\n(\cdot), \delta,\Omega,\gamma) \mapsto v\f[v\n, \delta,\Omega,\gamma](\cdot) 
 \end{align}
 from $ \{ u(\cdot) \in X_{\text{near}, \ep^{\tau}}^{4} \vert\Vert  u \Vert_{H^4}\leq  R \ep^{\frac{\beta}{2}} \}\times (0,\delta_*)\times\left(-\frac{1}{3},\frac{1}{3}\right)\times \T$  to $\{ v(\cdot)\in X_{\text{far}, \ep^{\tau}}^{s} \vert \Vert  v \Vert_{H^4}\leq R \ep^{\frac{1}{16}}\}$ and satisfying 
 \begin{align*}
v\f[v\n, \delta,\Omega,\gamma](\cdot) =  \mathscr{E}[v\n,v\f[v\n, \delta,\Omega,\gamma], \ur, \delta,\Omega](\cdot).  
 \end{align*}
 Furthermore, we have that 
  \begin{enumerate}[label=(\roman*), ref=\theTheorem(\roman*)]
\hangindent\leftmargin
   \item \label{pmtr_range:bounds}the function $v\f(\cdot)$ satisfies the bound
  \begin{align*}
 \Vert v\f[v\n, \delta,\Omega,\gamma]\Vert_{H^4} \lesssim \Theta(\ep;\tau)\Vert v\n\Vert_{H^4}  +(\ep^{3- \frac{\beta}{2}} + \ep^{1+\frac{3\beta}{2}}+ \ep^{2+\frac{\beta}{2}}),
\end{align*}
 where $\displaystyle{\Theta(\ep;\tau) = \mathcal{O}\left(\sum_{j=1}^3\ep^{3-j - 2\tau+ \frac{\beta}{2}(j-1)}\right)};$ 
 \item \label{pmtr_range:extension_to_0} 
 the mapping $(v\n(\cdot), \delta,\Omega,\gamma) \mapsto v\f[v\n, \delta,\Omega,\gamma](\cdot) $ can be extended continuously at $\delta =\ep=0$, with  $\displaystyle{v\f[\cdot,\delta,\Omega,\gamma](\cdot)\Big\vert_{\delta=0}:=0;}$
 \item \label{pmtr_range:continuity} the mapping $(v\n, \delta, \Omega,\gamma) \mapsto v\f[v\n,\delta,\Omega,\gamma](\cdot)$ is continuous from $X_{\text{near}, \ep^{\tau}}^{4}  \times [0,\delta_*)\times\left(-\frac{1}{3},\frac{1}{3}\right)\times\T$ to $X_{\text{far}, \ep^{\tau}}^{4} $; furthermore, it is 2$\pi$-periodic with respect to $\gamma$.
  \end{enumerate} 
\end{Proposition}
An auxiliary Lemma concerning the scaling of  $\ur^{(\delta,\om,\gamma)}(\cdot)$ with respect to its amplitude $\ep$ is presented next.  The result  is essentially given in \cite[Section 4 and Theorem 4.1]{Mielke}. 
\begin{Lemma}[Scaling of $\ur^{(\delta,\om,\gamma)}(\cdot)$ in $\ep= \ep(\delta,\Omega,\gamma)$]\label{scaling:u_roll}
  Given the equation \eqref{SH-eq} with $\mu \equiv 1$, there exists an $\ep= \ep(\delta,\Omega,\gamma)>0$ and a mapping $(-\delta_0, \delta_0)\times \T \ni (\delta,\gamma) \mapsto \ur^{(\delta,\om,\gamma)}(\cdot)\in  H_{\text{per}}^4([0,2\pi])$, and so that $\ur^{(\delta,\om,\gamma)}(\cdot)\Big|_{\delta =0 } =0.$ Furthermore, we can define $\widetilde{\ur^{(\delta,\om,\gamma)}} = \frac{\ur^{(\delta,\om,\gamma)}}{\ep} \in L^{\infty}([0,2\pi])$ and find a $\ep_{0}>0$ such that $\ep_{0}\lesssim \delta_{0}\lesssim \ep_{0}$. 
\end{Lemma}
\begin{Corollary}[Reparametrization of rolls by their amplitudes]\label{rolls:reparametrization}
 In the region $\om \in \left(\sqrt{1-\delta},\sqrt{1+\delta}\right)$ the rolls described in \eqref{rolls} can be reparametrized as functions of $(\ep, \om,\gamma)$, where
 \begin{align*}
  \delta = \delta(\ep, \om) =\sqrt{\frac{3\ep^2}{4} + (1- \om	^2)^2} + \mathcal{O}\left(\frac{3\ep^2}{4} + (1- \om^2)^2\right);
 \end{align*}
 moreover, this mapping is a homeomorphism whenever $\ep>0$ and $\om \in \left(\sqrt{1-\delta},\sqrt{1+\delta}\right)$.
\end{Corollary}
\begin{Proof}
 This result is a consequence of Mielke's derivation of the rolls using a Lyapunov-Schmidt reduction. Inspecting the proof in \cite[Theorem 4.1]{Mielke}, one notes that the Implicit Function Theorem can be applied\footnote{Our notation of $\ep, \om$ translates to Mielke's notation in \cite[Section 4]{Mielke}  as $\alpha$ and $\ep$, respectively.} either in $\ep^2$ or in $\delta^2$ from which one obtains either a function $\ep^2=\ep^2(\delta^2,\om)$, or $\delta^2=\delta^2(\ep^2,\om).$  Last, the mapping $(\ep,\om)\mapsto \delta(\ep,\om)$ is a homeomorphism, since it is the composition of a diffeomorphism with the homeomorphism $x\mapsto \sqrt{x}$ on $(0,\infty)$.
\end{Proof}

A few  useful consequences of this result are readily available, thanks to the fact that the nonlinear terms in \eqref{SH-eq-bifurcation-1} scale in $\ep$. Indeed, we can write
\begin{align}\label{nonlinear-scaling:rolls}
  \mathscr{N}^{(j)}[v,  \ur^{(\delta,\om,\gamma)}, \ep] =\mathscr{N}^{(j)}[v, \ep\widetilde{\ur^{(\delta,\om,\gamma)}}, \ep], \quad j \in \{1,2,3, 4\}.
\end{align}
We point out that  $\mathscr{N}^{(1)}[\cdot]$ depends explicitly on $\delta^2$; however, thanks to the choice of parameters in (\hyperref[H2]{H2}) any upper bound in terms of $\delta$  can be rewritten as an upper bound in terms of $\ep$. 
\begin{Lemma}\label{controling_Nl_terms}Assume (\hyperref[H1]{H1})-(\hyperref[H2]{H2}).  Let $u(\cdot)\in X_{\text{near}, \ep^{\tau}}^{4}$ be fixed and parameters  $\beta$ and $\tau$ satisfying the constraints in \eqref{pmtr_range}.   Given the nonlinear mappings $X_{\text{far}, \ep^{\tau}}^{4} \ni v(\cdot)\mapsto \mathscr{N}^{(j)}[u + v, \ur^{(\delta,\om,\gamma)}, \ep]$ for $j \in \{1, 2, 3, 4\}$, defined in \eqref{SH-eq-bifurcation-1},the following properties hold:
\begin{enumerate}[label=(\roman*), ref=\theTheorem(\roman*)]
\hangindent\leftmargin
\item\label{controling_Nl_terms:a}  Let $\Vert u\Vert_{H^4}<1$, and $\Vert v^{(1)}\Vert_{H^4} < 1$. For every $j \in \{1, 2, 3, 4\}$ there exists quantities $M_1>0$ and $E_j[ v^{(1)},u, \ep] \geq 0$ such that
  $$\displaystyle{\Vert \mathscr{N}^{(j)}[u +  v^{(1)}, \ur^{(\delta,\om,\gamma)}, \ep] \Vert_{L^2}\leq E_j[u,v^{(1)}, \ep],}$$
where
  \begin{align*}
E_j[u,v^{(1)},  \ep]  \leq \left\{   \begin{array}{cll}
    M_1 \ep^{3-j} \left(\Vert u^j\Vert_{L^2}+\Vert v^{(1)}\Vert_{H^4}^{j}\right) \leq M_1 \ep^{3-j} \left(\Vert u\Vert_{H^4}^{j}+\Vert v^{(1)}\Vert_{H^4}^{j}\right), & \text{for}& j=1,2,3;\\
    M_1(\ep^{3- \frac{\beta}{2}} + \ep^{1+\frac{3\beta}{2}}+ \ep^{2+\frac{\beta}{2}}), & \text{for}& j=4.
   \end{array}\right.
  \end{align*}
\item \label{controling_Nl_terms:b} For every $j \in \{1, 2, 3, 4\}$ there exists  quantities $M_2 >0$ and $D_j[v^{(2)}, v^{(3)}, \ep] \geq 0$ such that 
$$\displaystyle{\Vert \mathscr{N}^{(j)}[u +  v^{(2)}, \ur^{(\delta,\om,\gamma)}, \ep] - \mathscr{N}^{(j)}[u +  v^{(3)}, \ur^{(\delta,\om,\gamma)}, \ep] \Vert_{L^2}\leq D_j[v^{(2)},v^{(3)},u, \ep]\Vert v^{(2)} - v^{(3)} \Vert_{H^4},}$$
  where
  \begin{align*}
D_j[v^{(2)},v^{(3)},u, \ep]  \leq \left\{   \begin{array}{cll}
    M_2\ep^{3-j}\left(\Vert v^{(2)}\Vert_{H^4}^{j-1} + \Vert v^{(3)}\Vert_{H^4}^{j-1} + \Vert u\Vert_{H^4}^{j-1}\right), & \text{for}& j=1,2,3;\\
    0, & \text{for}& j=4.
   \end{array}\right.
  \end{align*}

 \end{enumerate}
\end{Lemma}
\begin{Proof}
For the sake of notation, we write $\ur(\cdot)$ to denote $\ur^{(\delta,\om,\gamma)}(\cdot).$ 

Whenever $j \in \{ 1, 2 , 3\}$  the proof exploits the scaling \eqref{nonlinear-scaling:rolls} in $\ep$,  the polynomial nature of the nonlinearity, and the Sobolev embedding $H^4(\R) \hookrightarrow L^{\infty}(\R)$ given in \eqref{Sobolev_embedding}.  Indeed, for all $p \geq 1$, $w \in H^4(\R)$, we have
   \begin{align*}
    \Vert w^p \Vert_{L^2} \leq \Vert w \Vert_{L^2}\Vert w \Vert_{L^{\infty}}^{p-1} \lesssim \Vert w \Vert_{L^2}\Vert w \Vert_{H^4}^{p-1};
   \end{align*}
in the particular case $j=2$ we resort to property (\hyperref[H2]{H2}), using the similarity $\delta \approx \ep$.  

We now prove the case $j =4$: the inequalities in  (ii) are trivial because there is no dependence  of  $\mathscr{N}^{(4)}[\ur^{(\delta,\om,\gamma)}, \ep](x)$  on $w(\cdot)$. Still in the case $j=4$, before proving (i) we recall that
\begin{align*}
 \mathscr{N}^{(4)}[\ur^{(\delta,\om,\gamma)}, \ep](x) &=    \chi(\ep^{\beta} x)( \chi^2(\ep^{\beta} x) -1)\left(\ur^{(\delta,\om,\gamma)}(x)\right)^3 \\
 &\quad + [(1+ \partial_x^2)^2,\chi(\ep^{\beta} x)]\ur^{(\delta,\om,\gamma)}(x) - \delta^2\chi(\ep^{\beta} x)(\mu -1) \ur^{(\delta,\om,\gamma)}(x).
\end{align*}
Lemma \ref{scaling:u_roll} implies that the first term has the  form
$$\chi(\ep^{\beta} x)( \chi^2(\ep^{\beta} x) -1)\left(\ur^{(\delta,\om,\gamma)}\right)^3(x) = \ep^3 \chi(\ep^{\beta} x)( \chi^2(\ep^{\beta} x) -1)(\widetilde{\ur^{(\delta,\om,\gamma)}})^3(x),$$  
hence a change of variables gives
\begin{align*}
\left\Vert\chi(\ep^{\beta} x)( \chi^2(\ep^{\beta} x) -1)\left(\ur^{(\delta,\om,\gamma)}\right)^3\right\Vert_{L^2(\R)} \leq \ep^{3 - \frac{\beta}{2}}\left\Vert \widetilde{\ur^{(\delta,\om,\gamma)}}\right\Vert_{L^{\infty}}\sqrt{\int_{\R} \vert \chi(y)( \chi^2(y) -1) \vert^2 \mathrm{d}y},
\end{align*}
where the integrand in the latter integral is localized, thanks to the properties of $\chi(\cdot)$ in (\hyperref[H1]{H1}). Similar reasoning shows that the  $L^2(\R)$ norm of the last integral is $\mathcal{O}(\ep^{2+\frac{\beta}{2}})$. With regards to the second term, initially notice that
$  [(1+ \partial_x^2)^2,\chi(\ep^{\beta} x)]\ur^{(\delta,\om,\gamma)} = \ep [(1+ \partial_x^2)^2,\chi(\ep^{\beta} x)]\widetilde{\ur^{(\delta,\om,\gamma)}}$
Then notice that the lowest order terms in $\ep$ come from terms in the commutant that contain one derivative of $\chi(\cdot)$, that is,
$$ \ep [(1+ \partial_x^2)^2,\chi(\ep^{\beta} x)]\widetilde{\ur^{(\delta,\om,\gamma)}}(x) = 4\ep^{1+\beta}\chi^{'}(\ep^{\beta} x)\partial_x(1 + \partial_x^2)\widetilde{\ur^{(\delta,\om,\gamma)}}(x)  +\mathcal{O}(\ep^{2\beta+1}).$$
   However, $\widetilde{\ur^{(\delta,\om,\gamma)}}(\cdot) = \cos(\cdot) + \mathcal{O}(\ep)$, hence  $\displaystyle{ 4\ep^{1+ \beta}\chi^{'}(\ep^{\beta} x)\partial_x(1 + \partial_x^2)\widetilde{\ur^{(\delta,\om,\gamma)}} = \mathcal{O}(\ep^{\beta+2})}.$ A change of variables then gives
\begin{align*}
\int_{\R}\vert \ep^{2+\beta} \chi'\left(\ep^{\beta}x\right)\vert^2\mathrm{d}x = \ep^{4 + \beta}\int_{\R}\vert \chi'\left(z\right)\vert^2\mathrm{d}z = \mathcal{O}\left(\ep^{4+\beta}\right),
\end{align*}
and this finishes the the proof of (i) when $j=4$.
\end{Proof}
We are now ready to prove the main result of this section. 
\begin{Proof}[of Proposition \ref{Prop:near_far_domination}] Throughout the proof we write $\mathscr{N}^{(j)}[v\n + v\f]$ to denote $\mathscr{N}^{(j)}[v\n + v\f,  \ur^{(\delta,\om,\gamma)}, \ep]$. We shall achieve the implicit parametrization of $v\f(\cdot)$ by $v\n(\cdot)$, $\delta$ and $\Omega$ by applying  the Contraction Mapping Theorem \cite[\S 2, Theorem 2.2]{Hale}. Without loss of generality we assume that $\ep<1$.   Recall from \eqref{nr_fr_decomposition:properties} that $ v(\cdot) = v\f(\cdot) + v\n(\cdot)$. Writing  $\displaystyle{\tilde{\mathbbm{1}}(\kp) = \left(1 - \mathbbm{1}_{\left\{+1 + \ep^{\tau}\mathcal{B}\right\}}(\kp) - \mathbbm{1}_{\left\{-1 + \ep^{\tau}\mathcal{B}\right\}}(\kp)\right)},$
a direct application of  Lemma \ref{Lemma_multiplier} implies that
 \begin{align}\label{multiplier_bound}
 \left \vert\frac{\tilde{\mathbbm{1}}(\kp)}{m(\kp; \mathscr{L})} \right\vert^{2} = \left \vert\frac{\tilde{\mathbbm{1}}(\kp)}{(1-\om^{2}\kp^{2})^{2}}\right\vert^{2}  \lesssim \frac{1}{\ep^{4\tau}}.
 \end{align}
 We use this inequality to estimate   \eqref{Fixed_point}:
\begin{align*}
 \Vert \mathscr{E}[v\n, v\f,  \delta,\Omega,\gamma]\Vert_{H^4}^2 &= \int_{\R}(1+ \vert \kp \vert^2)^4\vert  \mathscr{E}[v\n, v\f,  \delta,\Omega,\gamma](\kp)\vert^2\mathrm{d}\kp\\
&\leq \int_{\R}\frac{(1+ \kp^2)^4}{(1 - \om^2\kp^2)^4} \left(\tilde{\mathbbm{1}}(\kp)\right)^2 \left\vert\mathscr{F}\left[\sum_{j=1}^4 \mathscr{N}^{(j)}[v\n + v\f]\right]\right\vert^2(\kp)\mathrm{d} \kp\\
&\leq  \sum_{j=1}^4 \int_{\R}\frac{C^2}{ \ep^{4\tau} }\left\vert\mathscr{F}\left[ \mathscr{N}^{(j)}[v\n + v\f]\right]\right\vert^2(\kp)\mathrm{d} \kp\\
&\leq  \sum_{j=1}^4 \frac{C^2}{ \ep^{4\tau} }\left\Vert  \mathscr{N}^{(j)}[v\n + v\f]\right\Vert_{L^2}^2.
\end{align*}
An application of Lemma \ref{controling_Nl_terms} gives
\begin{align*}
 \Vert \mathscr{E}[v\n, v\f,  \delta,\Omega,\gamma]\Vert_{H^4} &\leq  \sum_{j=1}^4 \frac{C}{\ep^{2\tau}}\Vert\mathscr{N}^{(j)}[v\n + v\f] \Vert_{L^2}\\
 & \leq \sum_{j=1}^4 \frac{C}{\ep^{2\tau }} E_j[v\n, v\f,  \ur, \delta,\Omega,\gamma]  \nonumber\\
 &\leq  CM_1\left\{\sum_{j=1}^3 \ep^{3-j - 2\tau}\left( \Vert v\n\Vert_{H^4}^{j}+\Vert v\f\Vert_{H^4}^{j}\right)+\ep^{3- \frac{\beta}{2}-2\tau} + \ep^{1+\frac{3\beta}{2}-2\tau}+ \ep^{2+\frac{\beta}{2}-2\tau}\right\}.
\end{align*}
Now choosing $\Vert v\n\Vert_{H^4} \leq R\ep^{\frac{\beta}{2}}$, and  $\Vert v\f \Vert_{H^4}\leq R \ep^{\alpha}$, where $\alpha = \frac{1}{16}$, we get 
\begin{align*}
 \Vert \mathscr{E}[v\n, v\f,  \delta,\Omega,\gamma]\Vert_{H^4} & \leq C M_1\left\{\sum_{j=1}^3  \ep^{3-j - 2\tau}R^j\left( \ep^{j\frac{\beta}{2}}+ \ep^{j \alpha}\right)+ \ep^{3- \frac{\beta}{2}-2\tau} + \ep^{1+\frac{3\beta}{2}-2\tau}+ \ep^{2+\frac{\beta}{2}-2\tau}\right\}.
\end{align*}
Hence,
\begin{align*}
 \Vert \mathscr{E}[v\n, v\f,  \delta,\Omega,\gamma]\Vert_{H^4} & \leq CM_1\left\{\sum_{j=1}^3 \ep^{3-j - 2\tau}R^j\left(\ep^{j\frac{\beta}{2}}+\ep^{j \alpha}\right)+\ep^{3- \frac{\beta}{2}-2\tau} + \ep^{1+\frac{3\beta}{2}-2\tau}+ \ep^{2+\frac{\beta}{2}-2\tau}\right\}\nonumber \\
 & \leq 3C M_1R\left\{ \ep^{\frac{3}{2}- 2\tau} + \left(\sum_{j=1}^3\ep^{3+j(\alpha -1) - 2\tau}\right)+\ep^{3- \frac{\beta}{2}-2\tau} + \ep^{1+\frac{3\beta}{2}-2\tau}+ \ep^{2+\frac{\beta}{2}-2\tau}\right\},
\end{align*}
the last inequality being a consequence of the monotonic decay of the mapping $x \mapsto \ep^x$ in $x>0$, for $0<\ep<1$ and $\beta\geq 1$; we shall use this property repeatedly in what follows. After further simplification we get
\begin{equation*}
 \begin{split}
 &\Vert \mathscr{E}[v\n, v\f,  \delta,\Omega,\gamma]\Vert_{H^4}\\
 &\quad \leq  3C M_1R\left\{ \ep^{\frac{3}{2}- 2\tau} + \left(\sum_{j=1}^3\ep^{3+j(\alpha -1) - 2\tau}\right)+\ep^{3- \frac{\beta}{2}-2\tau} + \ep^{1+\frac{3\beta}{2}-2\tau}+ \ep^{2+\frac{\beta}{2}-2\tau}\right\}\\
 &\quad   \leq  3C M_1\left\{ \ep^{\frac{3}{2}- 2\tau-\alpha} + \left(\sum_{j=1}^3\ep^{3- j- 2\tau + \alpha(j-1)}\right)+\ep^{3- \frac{\beta}{2}-2\tau-\alpha} + \ep^{1+\frac{3\beta}{2}-2\tau-\alpha}+ \ep^{2+\frac{\beta}{2}-2\tau-\alpha}\right\}R \ep^{\alpha}\\
 &\quad =: \Theta_1( \alpha, \tau, \ep, \beta)R \ep^{\alpha}.
 \end{split}
\end{equation*}
In our second step, we rely again on \eqref{multiplier_bound} and Lemma \ref{controling_Nl_terms:b} to derive a few more estimates:
\begin{align*}
 &\Vert \mathscr{E}[v\n, v\f^{(1)},  , \delta, \Omega,\gamma]-\mathscr{E}[v\n, v\f^{(2)},  , \delta, \Omega,\gamma]\Vert_{H^4}\leq \sum_{j=1}^4 \frac{ C}{\ep^{2\tau}}\Vert\mathscr{N}^{(j)}[(v\n + v\f^{(1)})]- \mathscr{N}^{(j)}[(v\n + v\f^{(2)})]\Vert_{L^2}\nonumber \\
 &\qquad\leq \sum_{j=1}^4 \frac{ C}{\ep^{2\tau}} D_j[ v\f^{(1)}), v\f^{(2)},v\n, \ep] \Vert v\f^{(1)} - v\f^{(2)} \Vert_{H^4}  \nonumber\\
  &\qquad\stackrel{\text{Lemma}\, \ref{controling_Nl_terms:b}}{\leq} \sum_{j=1}^3 \frac{ C M_2}{\ep^{2\tau}} \ep^{3-j}\left(\Vert v\f^{(1)} \Vert_{H^4}^{j-1}+\Vert v\f^{(2)}   \Vert_{H^4}^{j-1}+\Vert v\n \Vert_{H^4}^{j-1} \right)\Vert v\f^{(1)} - v\f^{(2)} \Vert_{H^4}  \nonumber\\
  &\qquad\leq \sum_{j=1}^3 2 C M_2 R^{j-1}\ep^{3-j- 2\tau}\left(\ep^{\alpha(j-1)} +  \ep^{(j-1)\frac{\beta}{2}} \right)\Vert v\f^{(1)} - v\f^{(2)} \Vert_{H^4}  \nonumber\\
  &\qquad\leq 2 C M_2 \sum_{j=1}^3 \ep^{3-j- 2\tau}\left(\ep^{\alpha(j-1)} +  \ep^{(j-1)\frac{1}{2}} \right)\Vert v\f^{(1)} - v\f^{(2)} \Vert_{H^4}  \nonumber\\
  &\qquad\leq 2 C M_2\left(\sum_{j=1}^3  \left(\ep^{3-j- 2\tau +  \alpha(j-1)} +  \ep^{\frac{5}{2}-\frac{j}{2}- 2\tau} \right)\right)\Vert v\f^{(1)} - v\f^{(2)} \Vert_{H^4}.
\end{align*}
Thus, 
\begin{equation}\label{Prop:nr_eng_dmn-contraction}
\begin{split}
\Vert \mathscr{E}[v\n, v\f^{(1)},  \delta, \Omega,\gamma]-\mathscr{E}[v\n, v\f^{(2)},  \delta, \Omega,\gamma]\Vert_{H^4}& \leq 6 C M_2 \left(\sum_{j=1}^3  \left(\ep^{3-j- 2\tau +  \alpha(j-1)} \right) + \ep^{1- 2\tau}\right)\Vert v\f^{(1)} - v\f^{(2)} \Vert_{H^4}\\
 & =: \Theta_{2}( \alpha, \tau, \ep, \beta)\Vert v\f^{(1)} - v\f^{(2)} \Vert_{H^4}.
\end{split}
\end{equation}
Hence, uniformly for all $R, \tau$ and $\beta$ satisfying conditions  \eqref{pmtr_range} (recall that $\alpha = \frac{1}{16}$), one can choose $\ep_1>0$ and $\ep_2>0$ such that 
\begin{align*}
 \sup_{\ep \in [0,\ep_1)} \Theta_1( \alpha, \tau, \ep, \beta) \leq \frac{1}{2}, \quad \text{and} \quad \sup_{\ep \in [0,\ep_2)} \Theta_2(\alpha, \tau, \ep, \beta) \leq \frac{1}{2}.
\end{align*}
Now, choosing $\ep_0>0$ as in Lemma \ref{scaling:u_roll} and $\ep_*:= \min\{\ep_0, \ep_1, \ep_2\}$ we get that for all $\ep \in [0, \ep_*)$ the mapping $(v\n, v\f,  \ur, \delta,\Omega,\gamma) \mapsto \mathscr{E}[v\n, v\f, \delta, \Omega,\gamma], $ as we fix parameters $v\n$, $\delta$, $\Omega$,  and $\gamma$,  maps the set
$$(v\n, v\f,\delta, \Omega,\gamma)\in\{v\in X_{\text{near}, \ep^{\tau}}^{4}|\Vert v\Vert_{H^4}\leq R \ep^{\frac{\beta}{2}}\}\times\{v\in X_{\text{far}, \ep^{\tau}}^{4}|\Vert v\Vert_{H^4}\leq R\ep^{\frac{1}{16}}\}\times(0,\delta_*)\times \left(-\frac{1}{3},\frac{1}{3}\right)\times \T,$$
into $\{v\in H^4|\Vert v\Vert_{H^4}\leq R\ep^{\frac{1}{16}}\}$. Therefore, we can apply the Contraction Mapping Theorem to obtain the existence of a fixed point $v\f = v\f[v\n,\delta, \Omega,\gamma]$ for all $\ep\in [0, \ep_*)$, where $v\f[v\n,\delta, \Omega,\gamma]\Big\vert_{\ep =0} =0.$ This finishes the proof of the implicit parametrization of $v\f(\cdot)$ by  $(v\n(\cdot),\delta, \Omega, \gamma)$.

Now that  $v\f[v\n,\delta, \Omega,\gamma](\cdot)$ is well defined, we can finally study its properties. The bounds in  (i) are a consequence of 
\begin{align*}
& \Vert v\f[v\n,\delta, \Omega,\gamma]\Vert_{H^4} = \Vert \mathscr{E}[v\n, v\f,  \delta,\Omega,\gamma]\Vert_{H^4} \nonumber \\
 &\qquad \leq  CM_1\left\{\sum_{j=1}^3  \ep^{3-j - 2\tau} \Vert v\n\Vert_{H^4}^{j}+\sum_{j=1}^3 \ep^{3-j - 2\tau}\Vert v\f\Vert_{H^4}^{j}+(\ep^{3- \frac{\beta}{2}- 2\tau } + \ep^{1+\frac{3\beta}{2}- 2\tau }+ \ep^{2+\frac{\beta}{2}- 2\tau })\right\} \nonumber \\
 &\qquad \leq   CM_1 \sum_{j=1}^3  \left(\ep^{3-j - 2\tau+ \frac{\beta}{2}(j-1)}\right)\Vert v\n\Vert_{H^4}\nonumber   +  CM_1 \sum_{j=1}^3  \left(\ep^{3-j - 2\tau+\frac{1}{16}(j-1)}\right)\Vert v\f\Vert_{H^4}\nonumber\\ 
 &\qquad \quad + CM_1(\ep^{3- \frac{\beta}{2}- 2\tau } + \ep^{1+\frac{3\beta}{2}- 2\tau }+ \ep^{2+\frac{\beta}{2}- 2\tau }) 
 \end{align*}
Thanks to the parameter conditions in \eqref{pmtr_range} and the definition of $\ep_*$,  the term dependent of  $\Vert v\f\Vert_{H^4}$ on the right hand side can be absorbed, hence
\begin{align*}
 &\Vert v\f[v\n,\delta, \Omega,\gamma]\Vert_{H^4} \lesssim  \left(\sum_{j=1}^3\ep^{3-j - 2\tau+ \frac{\beta}{2}(j-1)}\right) \Vert v\n\Vert_{H^4}  + (\ep^{3- \frac{\beta}{2}- 2\tau } + \ep^{1+\frac{3\beta}{2}- 2\tau }+ \ep^{2+\frac{\beta}{2}- 2\tau }),
\end{align*}
whenever $\ep \in (0, \ep_*)$, and this finishes the proof of (i).

  Item (ii) is an immediate consequence of (i) taking  $\ep\downarrow 0$. 
  
  With regards to (iii), we must show that the mapping $(v\n, \delta, \Omega,\gamma) \mapsto v\f[v\n,\delta, \Omega,\gamma](\cdot)$ is continuous, this is also a consequence of the Contraction Mapping Principle, namely, it suffices to show that 
\begin{align*}
(v\n, \delta, \Omega,\gamma) \mapsto \mathscr{E}[v\n, v\f,  \delta,\Omega,\gamma]
\end{align*}
is a continuous mapping from $H^4(\R) \times \left[0,\delta_*\right)\times \left(-\frac{1}{3},\frac{1}{3}\right)\times \T$ to $H^4(\R)$; continuity with respect to $v\n(\cdot)$ is easily obtained by exchanging the roles of $(v\n,v\f)$ for $(v\f,v\n)$ in the estimates \eqref{Prop:nr_eng_dmn-contraction}. In this fashion, one obtains the similar bound
\begin{equation*}
\begin{split}
\Vert \mathscr{E}[v\n^{(1)}, v\f,  \delta, \Omega,\gamma]-\mathscr{E}[v\n^{(2)}, v\f,  \delta, \Omega,\gamma]\Vert_{H^4}  \leq \Theta_{2}\left( \frac{1}{16}, \tau, \ep, \beta\right)\Vert v\n^{(1)} - v\n^{(2)} \Vert_{H^4}.
\end{split}
\end{equation*}
Exploiting the pointwise continuity of the mappings $(\delta,\Omega,\gamma)\mapsto \mathscr{E}[v\n, v\f,  \delta,\Omega,\gamma](x)$ for any fixed $x\in \R$ and the fact that the $H^4(\R)$ norm of $\mathscr{E}[v\n, v\f,  \delta,\Omega,\gamma]$ is uniform in $\ep$, an application of the Lebesgue Dominated Convergence Theorem shows that this mapping is continuous  with respect to $(\delta, \Omega,\gamma)$. 
In order to prove $2\pi$-periodicity with respect to $\gamma$  note that $\gamma \mapsto \mathscr{E}[v\n^{(1)}, v\f,  \delta, \Omega,\gamma]$ only depends on $\gamma$ through $\ur^{(\delta,\om,\gamma)}$, which is $2\pi$-periodic in this parameter. By uniqueness of the fixed point, the result then follows, and this finishes the proof.  
\end{Proof}

\begin{Remark}[The fine balance between the blow-up rate in $\ep$ and the scaling of the nonlinearities in $\ep$] In \cite{Fefferman} the near component domination relies on two main ingredients:
\begin{enumerate}[label=(\roman*), ref=\theTheorem(\roman*)]
\item  counterbalancing the blow-up of the multiplier as $\ep \downarrow 0$, which in our case is encoded by Lemma \ref{Lemma_multiplier} and the inequality \eqref{multiplier_bound};
\item  an appropriate rescaling of the solutions, a feature that  strongly depends on the linear nature of the problem.
\end{enumerate}
 In the nonlinear case we are dealing with, resorting to rescaling is also crucial, otherwise singularities like the upper bound $\ep^{-4\tau}$ in \eqref{multiplier_bound} would be harmful upon understanding the regime $\ep\downarrow 0$. Hence, we rely on the fact that nonlinearity enhances the order of $\ep$ dependent parameters to balance out the singularity of $\displaystyle{\frac{1}{m(\kp;\mathscr{L})}}$ as $\ep\downarrow0$, a fact that should be contrasted with the case in \cite{Fefferman}, where this balance has to be found between the rate of spectral gap separation and the linear scaling.
\end{Remark}
%



\section{Desingularization, nonlinear interaction estimates and approximation results }\label{section:blow_up}
Proposition \ref{Prop:near_far_domination} establishes the parametrization by the near frequency components, whose behavior is still to be understood and  where we now concentrate our concerns at. We briefly recall that we have introduced unknown constants and  functions: the constant  $\beta$  first seen in \eqref{Ansatz}, the constant  $\tau$ introduced in the far/near decomposition \eqref{multiplier}, and the function $\chi(\cdot)$ that plays the role of an envelope function and was introduced in \eqref{Ansatz}. We begin by observing  that
\begin{align*}
\mathrm{supp}\left(\hv\n\right)\subset  \left\{ -1+ \ep^{\tau}\mathcal{B}\right\}\cup \left\{ 1+ \ep^{\tau}\mathcal{B}\right\}.
\end{align*}
Clearly, the set on the right hand side gets reduced to two points as $\ep \downarrow 0$, that is, $$X_{\text{near}, \ep^{\tau}}^{s}\Big \vert_{\ep=0} = \{0\} \subset H^s(\R), \quad \forall \tau >0,$$ which is exactly the regime we are interested at. In order to circumvent this issue, we extract the relevant properties of $v\n(\cdot)$  by desingularizing  this limit; that is,   using  blow-up variables in the frequency space, as done in  \cite[\S 6.4]{Fefferman}. As we shall see, this approach readily gives   another representation of $v\n(\cdot)$. But first a slight modification of the operator $\widetilde{\mathbb{P}\n}[\cdot]$ introduced in section \ref{Splitting:sec} is necessary: since the Near frequency set has two components, we can define the operators
\begin{equation}\label{n_n_p-secondary:Fourier}
\begin{split}
 \hv\mapsto  \widetilde{\mathbb{P}\n^{(\pm)}}[\hv](\tkp) &=  \mathbbm{1}_{\left\{ \tkp \in \ep^{\tau}\mathcal{B}\right\}}(\tkp)\hv\left(\pm1 + \tkp\right) =: \hv\n^{(\pm)}(\tkp),
\end{split}
\end{equation}
and associated physical space action
\begin{equation}\label{n_n_p-secondary:physical}
\begin{split}
 v\mapsto  \mathbb{P}\n^{(\pm)}[v](x) &=  \mathscr{F}^{-1}\left[ \widetilde{\mathbb{P}\n^{(\pm)}}\circ \mathscr{F}[v](\cdot)\right](x) = v\n^{(\pm)}(x),\quad x\in \R,
\end{split}
\end{equation}
whose properties are summarized below.
\begin{Proposition}[Recentered projections]\label{Prop:recentered}  Assume $ \tau >0$ fixed. Let $f(\cdot) \in L^2(\R)$, with decomposition $f(\cdot) = f\n(\cdot) + f\f(\cdot)\in X_{\text{near}, \ep^{\tau}}^{(0)}\oplus X_{\text{far}, \ep^{\tau}}^{(0)}$. Consider the operators $\widetilde{\mathbb{P}\n^{(\pm)}}[\cdot]$ as defined in \eqref{n_n_p-secondary:Fourier}. Then, the following properties hold:
\begin{enumerate}[label=(\roman*), ref=\theTheorem(\roman*)]\label{Reparametrized}
\hangindent\leftmargin
\item\label{Reparametrized:n_n_p-secondary:relation_near_prjct} Whenever $0 < \ep_{0}<1$,  it holds that
\begin{equation*}
\begin{split}
   \widetilde{\mathbb{P}\n^{(\pm)}}[\widehat{f}]\left(\tkp \mp 1\right) &= \mathbbm{1}_{\left\{ \tkp \in \pm1 +\ep^{\tau}\mathcal{B}\right\}}(\tkp)\widetilde{\mathbb{P}\n}[\widehat{f}]\left(\tkp\right), \quad \forall 0 <\ep <\ep_{0}.\end{split}
\end{equation*}
\item \label{Reparametrized:w_nr:before_blow_up} Writing  $f\n^{(\pm)}(x) = \mathscr{F}^{-1}\left[\widehat{f}\n^{(\pm)}\right](x)$, we have
\begin{align*}
f\n(x)=e^{+ix}f\n^{(+)}(x) + e^{-ix}f\n^{(-)}(x). 
\end{align*}
\item \label{Reparametrized:supports}For any given $\alpha \in \R$ we have 
$$ \mathrm{supp}\left(\mathscr{F}\left[e^{i\alpha (\cdot) }f\n^{(-)}(\cdot)\right]\right) \subset \alpha+ \ep^{\tau}B, \quad \text{and} \quad \mathrm{supp}\left(\mathscr{F}\left[e^{i\alpha (\cdot) }f\n^{(+)}(\cdot)\right]\right) \subset \alpha+ \ep^{\tau}B;$$
We can say then that $\displaystyle{f\mapsto \mathbb{P}\n^{(\pm)}[f](x): X_{\text{near}, \ep^{\tau}}^{s}\to  H_{\text{near}, \ep^{\tau}}^{s}(\R)}$;
 \item\label{Reparametrized:symmetry} $ \widetilde{\mathbb{P}\n^{(+)}}[\widehat{f}](\tkp) = \overline{\widetilde{\mathbb{P}\n^{(-)}}[\widehat{f}](-\tkp) }$;
 \item  \label{Reparametrized:conjugation} If  $f(\cdot)$ is real-valued, then  $f\n^{(+)}(x)  =  \overline{f\n^{(-)}(x)}$;
\item \label{Reparametrization:mplr_rpmt} Let $\mathcal{T}: L^2(\R) \mapsto L^2(R)$ be a mapping with a multiplier $m(\cdot; \mathcal{T})$, that is, $\widehat{\left(\mathcal{T}f\right)}(\xi) = m(\xi; \mathcal{T})\widehat{f}(\xi).$ Then, $\mathbb{P}\n^{(\pm)}[\widehat{\mathcal{T}f}](\tkp) = m\left(\pm1 + \tkp; \mathcal{T}\right)\mathbb{P}\n^{(\pm)}[\widehat{f}](\tkp).$
\end{enumerate}
\end{Proposition}
\begin{Proof}  Property (i) is obtained after a simple consequence of the definition of $\widetilde{\mathbb{P}\n}[\cdot]$ given in  section \ref{Splitting:sec}. 
Property (ii) also has a simple proof: writing  $f\n^{(\pm)}(x) = \mathscr{F}^{-1}\left[\widehat{f}\n^{(\pm)}\right](x)=\mathscr{F}^{-1}\circ\widetilde{\mathbb{P}\n}\circ\mathscr{F}[f](x)$, we have
\begin{equation*}
\begin{split}
 f\n(x) & = \frac{1}{2\pi}\int_{\{\tkp \in 1 +\ep^{\tau}B \}}\widehat{f}\left(\tkp\right)e^{i x\tkp}\mathrm{d}\tkp + \frac{1}{2\pi}\int_{\{\tkp \in -1 +\ep^{\tau}B \}}\widehat{f}\left(\tkp\right)e^{i x\tkp}\mathrm{d}\tkp.\\
 & = \frac{1}{2\pi}\int_{\{\tkp \in \ep^{\tau}B \}}\widehat{f}\left(1 + \tkp\right)e^{ix}e^{i x \tkp}\mathrm{d}\tkp + \frac{1}{2\pi}\int_{\{\tkp \in \ep^{\tau}B \}}\widehat{f}\left(-1 + \tkp\right)e^{-ix}e^{i x \tkp}\mathrm{d}\tkp\\
 &=e^{+ix}f\n^{(+)}(x) + e^{-ix}f\n^{(-)}(x)
\end{split}
\end{equation*}
Property (iii) is a direct consequence of the definition in \eqref{n_n_p-secondary:Fourier} allied to the Fourier transform property \eqref{Fourier_properties:translation};  the fact that $\displaystyle{\mathbb{P}\n^{(\pm)}[f](x)\in  H_{\text{near}, \ep^{\tau}}^{s}(\R)}$ whenever $f(\cdot)\in  X_{\text{near}, \ep^{\tau}}^{s}$ is a simple consequence of the definition of the spaces $  H_{\text{near}, \ep^{\tau}}^{s}(\R)$ given in \S\ref{notation}. In order to prove (iv) we study the action of the operators $\widetilde{\mathbb{P}\n^{(\pm)}}\circ \mathscr{F}[\cdot]$ on real-valued functions $f(\cdot)\in L^2(\R;\R)\cap L^1(\R;\R)$; a direct computations establishes the result:
\begin{equation*}
\begin{split}
  \widetilde{\mathbb{P}\n^{(+)}}[\widehat{f}](\tkp) &=   \mathbbm{1}_{\left\{ \tkp \in \ep^{\tau}\mathcal{B}\right\}}(\tkp)\mathscr{F}[f]\left(1 + \tkp\right)\\
  &=   \mathbbm{1}_{\left\{ \tkp \in \ep^{\tau}\mathcal{B}\right\}}(\tkp)\overline{\mathscr{F}[f]\left(-1 - \tkp\right)}\\
  &=   \mathbbm{1}_{\left\{ \tkp \in \ep^{\tau}\mathcal{B}\right\}}(-\tkp)\overline{\mathscr{F}[f]\left(-1 - \tkp\right)}\\
  & = \overline{\widetilde{\mathbb{P}\n^{(-)}}[\widehat{f}](-\tkp) }.
\end{split}
\end{equation*}
 We now turn to (v): thanks to definition of $\widehat{f}\n(\cdot)$ given in \S \ref{Splitting:sec}, we can use item (i) to write
\begin{align*}
\widehat{f}\n^{(\pm)}\left(\tkp\right) =  \mathbbm{1}_{\{\tkp \in \ep^{\tau}B\}}\widehat{f}\n\left(\pm 1 + \tkp\right),
\end{align*}
Using property (iii) we can write 
 $\displaystyle{\widehat{f}\n^{(+)}\left(\tkp\right) = \overline{\widehat{f}\n^{(-)}\left(-\tkp\right)} }$
and consequently, in physical space, we have
\begin{equation*}
 \begin{split}
  f\n^{(+)}(x) &= \mathscr{F}^{-1}\left[\widehat{f}\n^{(+)}\left(\cdot\right)\right](x)= \mathscr{F}^{-1}\left[\overline{\widehat{f}\n^{(-)}\left( - \cdot\right)}\right](x)= \frac{1}{2\pi}\int_{\R}\overline{\widehat{f}\n^{(-)}\left( - \tkp\right)}e^{i\tkp x}\mathrm{d}\tkp\\
  &= \frac{1}{2\pi}\int_{\R}\overline{\widehat{f}\n^{(-)}\left( - \tkp\right)e^{-i\tkp x}}\mathrm{d}\tkp\\
  & =  \overline{f\n^{(-)}(x)},
\end{split}
\end{equation*}
which proves (v). 

Last, the proof of (vi) is obtained after a direct computation: 
\begin{align*}
 \mathbb{P}\n^{(\pm)}[\widehat{\mathcal{T}f}](\tkp) = \mathbbm{1}_{\left\{ \ep^{\tau}\mathcal{B}\right\}}(\tkp)m\left(\pm1 + \tkp; \mathcal{T}\right)\widehat{f}\left(\pm1 + \tkp\right)= m\left(\pm1 + \tkp; \mathcal{T}\right)\mathbb{P}\n^{(\pm)}[\widehat{f}](\tkp),
\end{align*}
and we are done.
\end{Proof}
\subsection{Desingularization in Fourier space and the representation of \texorpdfstring{$v\n(\cdot)$}{vnear(.)} as a Ginzburg-Landau type approximation} Now that we are able to center our parametrizations, we apply the results derived in the previous section to construct the functions 
\begin{align*}
v\n^{(+)}(\cdot)= \mathbb{P}^{(+)}\n[v\n](\cdot) \in  H_{\text{near}, \ep^{\tau}}^{4}(\R),\quad  v\n^{(-)}(\cdot)=\mathbb{P}^{(-)}\n[v\n](\cdot) \in  H_{\text{near}, \ep^{\tau}}^{4}(\R), 
\end{align*}
and their corresponding Fourier transforms, given respectively by  
\begin{align*}
 \hv\n^{(+)}(\cdot)= \widetilde{\mathbb{P}\n^{(+)}}[\hv](\cdot), \quad \text{ and} \quad  \hv\n^{(-)}(\cdot)=\widetilde{\mathbb{P}\n^{(-)}}[\hv](\cdot).
\end{align*}
Our construction is motivated by that in  \cite[Equation 6.49]{Fefferman}): we first define functions $g\p(\cdot)$ and $g\m(\cdot)$ in the following manner,
\begin{align*}
 \widehat{g}_{\pm1}\left(\xi\right)= \widehat{g}_{\pm1}\left(\frac{\tkp}{\ep^{\beta}}\right):= \hv\n^{(\pm)}\left(\ep^{\beta}\xi\right),\quad\text{where} \quad  \xi = \frac{\tkp}{\ep^{\beta}}\quad \text{and}\quad g_{\pm1}(x) = \mathscr{F}^{-1}\left[\widehat{g}_{\pm1}\right](x).
\end{align*}
By construction,  $\displaystyle{g{\pm1}(\cdot) \in  H_{\text{near}, \ep^{\tau-\beta}}^{4}(\R)}$.  Using the identity $\displaystyle{\mathbbm{1}_{\{A\}}\left(\frac{x}{\zeta}\right) = \mathbbm{1}_{\{\zeta A\}}(x) }$ (whenever $\zeta >0$) and the properties \eqref{nr_fr_decomposition:properties} we get
\begin{align*}
\widehat{g}_{\pm1}(\xi)&= \mathbbm{1}_{\left\{\ep^{\tau}\mathcal{B}\right\}}(\ep^{\beta} \xi)\widehat{g}_{\pm1}(\xi)=\mathbbm{1}_{\left\{\ep^{\tau-\beta}\mathcal{B}\right\}}(\xi)\widehat{g}_{\pm1}(\xi).
\end{align*}
The blow-up in Fourier space induces a rescaling in the physical space; indeed, thanks to the identities \eqref{Fourier_properties:dilation}, we get that
\begin{align*}
 v\n^{(\pm)}(x) = \mathscr{F}^{-1}\left[\hv\n^{(\pm)}(\cdot)\right](x) = \mathscr{F}^{-1}\left[\widehat{g}_{\pm1}\left(\frac{\cdot}{\ep^{\beta}}\right)\right](x) = \ep^{\beta}\mathscr{F}^{-1}\left[\widehat{g}_{\pm1}\right](\ep^{\beta}x) = \ep^{\beta}g_{\pm1}(\ep^{\beta}x). 
\end{align*}
When combined with Proposition \ref{Reparametrized:w_nr:before_blow_up} it implies that
\begin{equation}\label{blow_up:near}
\begin{split}
 v_{near}(x)& = \ep^{\beta} e^{+i  x} \gp(\ep^{\beta} x) + \ep^{\beta} e^{-i x} \gm(\ep^{\beta} x),
 \end{split}
\end{equation}
with $ \gp(\ep^{\beta}x)  = \overline{g\m(\ep^{\beta}x)}$ due to  Proposition \ref{Reparametrized:conjugation}. In fact,  this representation  goes back to modulation theory and the Ginzburg-Landau formalism (see for instance  \cite[\S 3]{Schneider}; see also the Ansatz in \cite[Equation (1.2)]{Mielke_Schneider-cubic}). 

\subsection{The functions $\gm(\cdot)$, $\gp(\cdot)$ and the topology of the space they are in}
According to Remark \ref{embedding_band_limited}, we have
\begin{align}\label{Prop:nr_eng_dmn:trivial}
\Vert v\n\Vert_{H^4(\R)}\approx\Vert v\n\Vert_{L^2(\R)}\approx  \ep^{\frac{\beta}{2}}\Vert g\m\Vert_{L^2(\R)} +\ep^{\frac{\beta}{2}}\Vert g_1\Vert_{L^2(\R)},
\end{align}
therefore we can write  $v\f = v\f[v\n,\delta, \Omega,\gamma]$ whenever  $\ep$ sufficiently small because  $v\n(\cdot)$ is indeed an element in the sets used in the Contraction Mapping argument as applied in Proposition \ref{Prop:near_far_domination}. 

In fact, \eqref{Prop:nr_eng_dmn:trivial} immediately implies  the following result:
\begin{Lemma}[The Ginzburg-Landau representation as a mapping]\label{GL:representation}  For any fixed $\ep \geq0$, consider  \eqref{blow_up:near},
$$(\gp(\cdot), \gm(\cdot)) \mapsto v_{near}(x) = v_{near}[\gp(\cdot), \gm(\cdot)](x) = \ep^{\beta} e^{+i  x} \gp(\ep^{\beta} x) + \ep^{\beta} e^{-i x} \gm(\ep^{\beta} x).$$
This mapping is continuous from $H_{\text{near}, \ep^{\tau-\beta}}^{0}(\R)\times H_{\text{near}, \ep^{\tau-\beta}}^{0}(\R)$ to $X_{\text{near}, \ep^{\tau}}^{4}(\R)$. Furthemore, using continuity of the mapping $H^2(\R) \to L^2(\R)$, we conclude that this mapping is also continuous from  $H_{\text{near}, \ep^{\tau-\beta}}^{2}(\R)\times H_{\text{near}, \ep^{\tau-\beta}}^{2}(\R)$ to $X_{\text{near}, \ep^{\tau}}^{4}(\R).$
\end{Lemma}

We must highlight a few things concerning this result. First, the inequality \eqref{Prop:nr_eng_dmn:trivial} is misleading, for it gives the impression that control of $L^2(\R)$ norms of $\gm(\cdot)$ and $\gp(\cdot)$ is enough to control nonlinearities and then reduce the problem to a simpler bifurcation equation. Second, the discussion of the mapping $(\gp(\ep^{\beta} x), \gm(\ep^{\beta} x)) \mapsto v_{near}(x)$ in the topology $H_{\text{near}, \ep^{\tau-\beta}}^{2}(\R)\times H_{\text{near}, \ep^{\tau-\beta}}^{2}(\R)$ to $H_{\text{near}, \ep^{\tau}}^{4}(\R),$ seems a bit far-fetched: why $H_{\text{near}, \ep^{\tau-\beta}}^{2}(\R)$? As we shall soon see, this is a matter of convenience: this norm will provide an easier  way to control nonlinear and ``higher order" terms in $\ep$. In \S \ref{sec:reduced_equation+approx+solv} we  show that we can approximate and reduce our problem even further. Therefore, henceforth we shall look for solutions to problem \eqref{SH-eq-bifurcation-1} under the condition that
$$(\gm(\cdot), \gp(\cdot)) \in H_{\text{near}, \ep^{\tau-\beta}}^{2}(\R)\times H_{\text{near}, \ep^{\tau-\beta}}^{2}(\R)\subset H^{2}(\R)\times H^{2}(\R).$$  
Last, we point out that the proportionality $\Vert v\n\Vert_{H^4(\R)}\approx\Vert v\n\Vert_{L^2(\R)}$ asserted in \eqref{Prop:nr_eng_dmn:trivial} is uniform for all $\ep>0$ sufficiently small because the support of $\displaystyle{\widehat{v\n}(\cdot)}\in X_{\text{near}, \ep^{\tau}}^{4}$ grows with order $\mathcal{O}(\ep^{\tau})$. This property is in enormous contrast to the cases $\displaystyle{g_{\pm1}(\cdot) \in H_{\text{near}, \ep^{\tau-\beta}}}$, for which the equivalence in norm still holds for all $\ep>0$, but is not uniform for all  $\ep>0$ sufficiently small.
\subsection{Irrelevant nonlinearities and interaction Lemmas} Plugging $v\f(\cdot) = v\f[v\n,\delta, \Omega,\gamma](\cdot)$ into \eqref{mode_system_cutoff-a1}, yields
\begin{align*}
m(\kp; \mathscr{L})\hv\n(\kappa)  &=\sum_{j=1}^4 \widetilde{\mathbb{P}\n}\circ\mathscr{F}\left[ \mathscr{N}^{(j)}[v\n + v\f[v\n,\delta, \Omega,\gamma], \ur^{(\delta, \om,\gamma)},\delta \right](\kp).
\end{align*}
 Proposition \ref{Reparametrized:n_n_p-secondary:relation_near_prjct} allow us to  rewrite this equation in terms of $\hv\n^{(\pm)}(\cdot)$; we obtain two equations
which, in the blow-up variable $\kp = \ep^{\beta}\xi$, correspond to
\begin{align}\label{red_eqt:step_0}
m\left(\pm1+ \ep^{\beta}\xi; \mathscr{L}\right)\hv\n^{(\pm)}(\ep^{\beta}\xi) &=\sum_{j=1}^4 \widetilde{\mathbb{P}\n^{(\pm)}}\circ\mathscr{F}\left[ \mathscr{N}^{(j)}[v\n + v\f[v\n,\delta, \Omega,\gamma], \ur^{(\delta, \om,\gamma)},\delta] \right](\ep^{\beta}\xi).
\end{align}
We shall rewrite these equations in such a way that: (i) the linear terms in $v\n(\cdot)$ are written explicitly, and  (ii) the non-homogeneous term (that stems from $\mathscr{N}^{(4)}[v\n + v\f[v\n,\delta, \Omega,\gamma]](\cdot)$) is highlighted. 

The starting point consists of simplifying the left hand side of \eqref{red_eqt:step_0}: plugging $\om^2 = 1 + \delta\Omega$ as given by (\hyperref[H2]{H2}), we expand to get
\begin{equation*}\label{red_eqt:step_1}
\begin{split}
 m\left(\pm1+ \tkp; \mathscr{L}\right)\hv\n^{(\pm)}(\tkp) &= -\om^2 \ep^{2\beta}\xi^2 \left[2 \mp  \om \ep^{\beta}\xi\right]^2\widehat{g}_{\mp1}(\xi)= - 4 \ep^{2\beta}\xi^2 \widehat{g}_{\pm1}(\xi) + \mathscr{A}_1^{(\pm)}[v\n, \ep], 
\end{split}
\end{equation*}
where 
$$\mathscr{A}_1^{(\pm)}[v\n, \ep] = \mathcal{O}\left((\delta\Omega\ep^{2\beta}\xi^2 + \ep^{3\beta}\xi^{3} + \ep^{4\beta}\xi^{4})\widehat{g}_{\pm1}(\xi)\right)=\mathcal{O}\left((\ep^{2\beta+1}\xi^2 + \ep^{3\beta}\xi^{3} + \ep^{4\beta}\xi^{4})\widehat{g}_{\pm1}(\xi)\right).$$
The last equality is a consequence of (\hyperref[H2]{H2}), namely, $\delta \approx \ep$. Since $v\mapsto \mathscr{N}^{(1)}[v , \ur^{(\delta, \om,\gamma)},\delta]$ is linear and $\mathscr{N}^{(4)}[v , \ur^{(\delta, \om,\gamma)},\delta]=  \mathscr{N}^{(4)}[ \ur^{(\delta, \om,\gamma)},\delta]$ is independent  of  $v(\cdot)$, we rewrite   \eqref{mode_system_cutoff-a1} as 
\begin{equation}\label{second_main_eqt}
 \begin{split}
 &- 4\om^2 \ep^{2\beta}\xi^2 \widehat{g}_{\pm1}(\xi) - \widetilde{\mathbb{P}\n^{(\pm)}}\circ\mathscr{F}\left[\mathscr{N}^{(1)}[v\n , \ur^{(\delta, \om,\gamma)},\delta] +  \mathscr{N}^{(4)}[\ur^{(\delta, \om,\gamma)},\delta]\right](\ep^{\beta}\xi)  \\
  &\qquad =: -  \mathscr{A}_1^{(\pm)}[v\f[v\n,\delta, \Omega,\gamma], \ep] + \mathscr{A}_2^{(\pm)}[v\f[v\n,\delta, \Omega,\gamma], \ep],
 \end{split}
\end{equation}
with new remaining term
\begin{align*}
&\mathscr{A}_2^{(\pm)}[v\f[v\n,\delta, \Omega,\gamma], \ep]=\widetilde{\mathbb{P}\n^{(\pm)}}\circ\mathscr{F}\left[ \mathscr{N}^{(1)}[v\f[v\n,\delta, \Omega,\gamma]]+\sum_{j=2}^3 \mathscr{N}^{(j)}[v\n + v\f[v\n,\delta, \Omega,\gamma], \ur^{(\delta, \om,\gamma)},\delta] \right].
\end{align*}
The main goal in this section is showing that the right hand side is small in $\ep$, $\Vert \gm \Vert_{H^2(\R)}$, and $\Vert \gp \Vert_{H^2(\R)}$.
\begin{Proposition}[Irrelevant nonlinearities]\label{Prop:irrelevant} Given $v\n(\cdot)$ and $v\f(\cdot)$ obtained in Proposition \ref{Prop:near_far_domination}, where $v\n(\cdot)$ is of the form \eqref{blow_up:near}. Assume the parameters $\beta$ and $\tau$ satisfying the constraints \eqref{pmtr_range}. Then, assuming $g_{\pm 1}(\cdot) \in H_{\text{near}, \ep^{\tau-\beta}}^{2}(\R)$, we have 
 
 \begin{enumerate}[label=(\roman*), ref=\theTheorem(\roman*)]
\hangindent\leftmargin
  \item \label{Prop:irrelevant:part1}$\Vert\mathscr{A}_1[v\n,\ep]\Vert_{L^2(\R_{\xi})} \approx \Vert(\ep^{2\beta+1}\xi^2 +\ep^{3\beta}\xi^{3} + \ep^{4\beta}\xi^{4})\widehat{g_{\pm 1}}(\xi)\Vert_{L^2(\R)} \lesssim \ep^{\tau +  2\beta}\Vert g_{\pm 1} \Vert_{H^2(\R)} = o(\ep^{2\beta})$; 

  \item \label{Prop:irrelevant:part2}$\Vert \mathscr{A}_2^{(\pm)}[v\f[v\n,\delta, \Omega,\gamma], \ep]\Vert_{L^2(\R_{\xi})} = o(\ep^{2\beta}) + \ep^{2\beta}\mathcal{O}\left(\Vert \gm\Vert_{L^2(\R)}^2 +\Vert \gp\Vert_{L^2(\R)}^2\right).$ 
 \end{enumerate}
\end{Proposition}
We prove only part (i) for now; the result is essentially a consequence of the interactions studied in this section. Part (ii) is derived as a consequence of more refined analysis, presented in the next section. 

\begin{Observation}[The importance of $\tau>0$]
Lemma \ref{Prop:irrelevant} do not hold in the case $\tau=0$.
\end{Observation}
\begin{Proof}[of Proposition \ref{Prop:irrelevant:part1}]
 Using triagle inequality, it suffices to bound each term separately. It is straightforward to find an upper bound to the first term:
 $$\left\Vert \ep^{2\beta+1}\xi^2\widehat{g_{\pm 1}}\right\Vert \lesssim \ep^{2\beta+1}\Vert\widehat{g_{\pm 1}} \Vert_{H^2(\R)}\lesssim \ep^{\tau + 2\beta}\Vert\widehat{g_{\pm 1}} \Vert_{H^2(\R)}= o(\ep^{2\beta}).$$
  For the other terms, we rely on the identities \eqref{nr_fr_decomposition:properties}, specifically,  $\widehat{g}_{\pm1}(\xi)=\mathbbm{1}_{\left\{\ep^{\tau-\beta}\mathcal{B}\right\}}(\xi)\widehat{g}_{\pm1}(\xi)$. Indeed, 
 \begin{align*}
  &\left\Vert\ep^{3\beta}\xi^3\widehat{g_{\pm 1}}(\xi)\right\Vert_{L^2(\R_{\xi})}^2 \\
  &\lesssim \vert\ep\vert^{6\beta}\int_{\vert\xi\vert \leq \ep^{\tau -\beta}} \vert \xi\vert^6\left\vert \widehat{g_{\pm 1}}(\xi)\right\vert^2\mathrm{d}\xi \lesssim \vert\ep\vert^{2(\tau +   2 \beta)}\int_{\vert\xi\vert \leq \ep^{\tau -\beta}} \vert \xi\vert^4\left\vert \widehat{g_{\pm 1}}(\xi)\right\vert^2\mathrm{d}\xi \lesssim \vert\ep\vert^{2(\tau + 2 \beta)}\left\Vert g_{\pm 1}\right\Vert_{H^2(\R)}^2 .
 \end{align*}
 Similarly, one obtains
 \begin{align*}
  &\left\Vert \ep^{4\beta}\xi^4\widehat{g_{\pm 1}}(\xi)\right\Vert_{L^2(\R_{\xi})}^2 \lesssim \vert\ep\vert^{8\beta}\int_{\vert\xi\vert \leq \ep^{\tau -\beta}} \vert \xi\vert^8\left\vert \widehat{g_{\pm 1}}(\xi)\right\vert^2\mathrm{d}\xi  
  \lesssim \vert\ep\vert^{4(\tau + \beta)}\left\Vert g_{\pm 1}\right\Vert_{H^2(\R)}^2,
 \end{align*}
and this finishes the proof of part (i).
 \end{Proof}
 
For the second part of this Proposition, we first need to understand how to control $g_{\pm1}(\cdot)$ controls $v\n(\cdot)$; later on we study the nonlinear interaction between $v\n(\cdot)$ and $v\f(\cdot)$.

\begin{Lemma}[Nonlinear interaction Lemma I]\label{Nl_itr_Lm_1} Let $\ep_*>0$ be given in Proposition \ref{Prop:near_far_domination} and  $\ep \in (0, \ep_*)$. Assume $v\n(\cdot)$ is of the form \eqref{blow_up:near}. Then,
 \begin{align}\label{Nl_itr_inequ:1}
  \Vert v\n^p\Vert_{H^4(\R)} \lesssim  \ep^{\beta(p -\frac{1}{2})}\left(\Vert \gm^p\Vert_{L^2(\R)} + \Vert \gp^p \Vert_{L^2(\R)}  \right), \quad p\in \mathbb{N}\setminus \{0\}.
 \end{align}
 Furthermore, we can use the above inequality to improve the estimate  given in Proposition \ref{pmtr_range:bounds}, that is,
 \begin{align}\label{Prop:nr_eng_dmn:far_bounded_by_near_improved}
   \Vert v\f\Vert_{H^4(\R)}\lesssim \Lambda(\ep, \tau) \ep^{\frac{\beta}{2}}\Vert g_{\pm 1}\Vert_{H^1(\R)}  + \ep^{-2\tau}(\ep^{3- \frac{\beta}{2}} + \ep^{1+\frac{3\beta}{2}}+ \ep^{2+\frac{\beta}{2}}),
 \end{align}
where $\Lambda(\ep, \tau) := \mathcal{O}\left(\ep^{2- 2\tau} + \ep^{1+ \frac{\beta}{2}- 2\tau} +  \ep^{3\frac{\beta}{2}- 2\tau}\right)$. Furthermore,  $v\f(\cdot)$ scales in $\ep$ as
 \begin{align}\label{Prop:nr_eng_dmn:far_bounded_by_near_improved_II}
   \Vert v\f\Vert_{H^4(\R)}\lesssim  \left(\ep^{2 - 2\tau +  \frac{\beta}{2}} + \ep^{1 - 2\tau +  \frac{3\beta}{2}}+ \ep^{\frac{5\beta}{2}- 2\tau}+ \ep^{3- \frac{\beta}{2}- 2\tau }\right). 
 \end{align}

\end{Lemma}
\begin{Remark}
Inequality  \eqref{Prop:nr_eng_dmn:far_bounded_by_near_improved} is not a direct consequence of the inequality \eqref{Prop:nr_eng_dmn:trivial}. Indeed, the latter implies only that
$$ \Vert v\n^p\Vert_{H^4(\R)}\leq \Vert v\n\Vert_{H^4(\R)}^p \lesssim \mathcal{O}(\ep^{\frac{p\beta}{2}}), \qquad p\in \mathbb{N}\setminus \{0\}.$$ 
 Unfortunately, this upper bound is not good enough: as we will see later on, in order to obtain a reduced equation we need to derive better estimates.
\end{Remark}

\begin{Proof}
This is analogous to the proof of \cite[Lemma 6.9]{Fefferman}. Using \eqref{blow_up:near},
\begin{align*}
 v\n^p(x) = \ep^{p\beta}\left(g\m (\ep^{\beta} x)e^{-ix}+ g_1 (\ep^{\beta} x)e^{ix}\right)^p.
\end{align*}
In Fourier space, for any $p \in \mathbb{N}\setminus\{0\}$ the function $v\n^p(\cdot)$ corresponds to a convolution of band-limited functions, therefore it is also a band limited function (cf. \cite[Proposition 4.18]{Brezis}). We conclude that $\Vert v\n^p\Vert_{H^4}\lesssim \Vert v\n^p\Vert_{L^2}$ holds. Finally, the proof of \eqref{Nl_itr_inequ:1} follows upon integration,  using a change of variables:
\begin{align*}
 \Vert v\n^p \Vert_{H^4} \leq\Vert v\n^p \Vert_{L^2} & \leq \ep^{p\beta} \Vert \left(g\m(\ep^{\beta} \cdot)\right)^p \Vert_{L^2} + \ep^{p\beta} \Vert\left( g\p(\ep^{\beta} \cdot)\right)^p \Vert_{L^2} \\
 & \lesssim \ep^{\beta(p - \frac{1}{2})} \Vert g\m^p(\cdot) \Vert_{L^2} + \ep^{\beta(p - \frac{1}{2})} \Vert g\p^p( \cdot) \Vert_{L^2}.
\end{align*}
Now we show \eqref{Prop:nr_eng_dmn:far_bounded_by_near_improved}: using the previous inequality, Lemma \ref{controling_Nl_terms:a}, and the Sobolev Embedding $H^4(\R) \hookrightarrow L^{\infty}(\R)$, we obtain
\begin{align*}
 &\Vert v\f[v\n,\delta, \Omega,\gamma]\Vert_{H^4}= \Vert \mathscr{E}[v\n, v\f,  \delta,\Omega,\gamma]\Vert_{H^4} \nonumber \\
 & \leq  CM_1\left\{\sum_{j=1}^3  \ep^{3-j - 2\tau} \Vert v\n^{j}\Vert_{L^2}+\sum_{j=1}^3 \ep^{3-j - 2\tau}\Vert v\f\Vert_{H^4}^{j}+(\ep^{3- \frac{\beta}{2}- 2\tau } + \ep^{1+\frac{3\beta}{2}- 2\tau }+ \ep^{2+\frac{\beta}{2}- 2\tau })\right\} \nonumber \\
 & \leq  CM_1\left\{\sum_{j=1}^3  \ep^{3-j - 2\tau} \Vert v\n^{j-1}\Vert_{L^2}\Vert v\n\Vert_{H^4}+\sum_{j=1}^3 \ep^{3-j - 2\tau}\Vert v\f\Vert_{H^4}^{j}+(\ep^{3- \frac{\beta}{2}- 2\tau } + \ep^{1+\frac{3\beta}{2}- 2\tau }+ \ep^{2+\frac{\beta}{2}- 2\tau })\right\} \nonumber \\
 & \leq 3 CM_1(1 + R^2)\left(\ep^{2- 2\tau} + \ep^{1+ \frac{\beta}{2}- 2\tau}+ \ep^{3\frac{\beta}{2}- 2\tau}\right)\Vert v\n\Vert_{H^4}\nonumber \\
 & \quad  +\sum_{j=1}^3  \left( CM_1\ep^{3-j - 2\tau+\alpha(j-1)}\right)\Vert v\f\Vert_{H^4}  + CM_1(\ep^{3- \frac{\beta}{2}- 2\tau } + \ep^{1+\frac{3\beta}{2}- 2\tau }+ \ep^{2+\frac{\beta}{2}- 2\tau }). 
 \end{align*}
Thanks to the parameter conditions in \eqref{pmtr_range},  the term depending in $v\f(\cdot)$ on the right hand side can be absorbed to the left hand side, yielding
\begin{align*}
 \Vert v\f[v\n,\delta, \Omega,\gamma]\Vert_{H^4} &\lesssim \left(\ep^{2- 2\tau} + \ep^{1+ \frac{\beta}{2}- 2\tau}+ \ep^{3\frac{\beta}{2}- 2\tau}\right) \Vert v\n\Vert_{H^4}  + (\ep^{3- \frac{\beta}{2}- 2\tau } + \ep^{1+\frac{3\beta}{2}- 2\tau }+ \ep^{2+\frac{\beta}{2}- 2\tau }),
\end{align*}
and this finishes the proof of  \eqref{Prop:nr_eng_dmn:far_bounded_by_near_improved}.  

We conclude with a proof of \eqref{Prop:nr_eng_dmn:far_bounded_by_near_improved_II}:
\begin{align*}
 &\Vert v\f[v\n,\delta, \Omega,\gamma]\Vert_{H^4} =  \Vert \mathscr{E}[v\n, v\f,  \delta,\Omega,\gamma]\Vert_{H^4} \nonumber \\
 &\qquad  \leq  CM_1\left\{\sum_{j=1}^3  \ep^{3-j - 2\tau} \Vert v\n^{j}\Vert_{L^2}+\sum_{j=1}^3 \ep^{3-j - 2\tau}\Vert v\f\Vert_{H^4}^{j}+ (\ep^{3- \frac{\beta}{2}- 2\tau } + \ep^{1+\frac{3\beta}{2}- 2\tau }+ \ep^{2+\frac{\beta}{2}- 2\tau })\right\} \nonumber \\
 &\qquad  \leq  CM_1\left\{\sum_{j=1}^3  \ep^{3-j - 2\tau} \ep^{\beta(j - \frac{1}{2})}+\sum_{j=1}^3 \ep^{3-j - 2\tau}\Vert v\f\Vert_{H^4}^{j}+ (\ep^{3- \frac{\beta}{2}- 2\tau } + \ep^{1+\frac{3\beta}{2}- 2\tau }+ \ep^{2+\frac{\beta}{2}- 2\tau })\right\}.
 \end{align*}
As in the previous case, we finalize the proof absorbing the right hand side term depending on $v\f(\cdot)$; using \eqref{Prop:nr_eng_dmn:far_bounded_by_near_improved} we obtain
\begin{align*}
 \Vert v\f[v\n,\delta, \Omega,\gamma]\Vert_{H^4}  & \lesssim \ep^{2 - 2\tau +  \frac{\beta}{2}} + \ep^{1 - 2\tau +  \frac{3\beta}{2}}+ \ep^{\frac{5\beta}{2}- 2\tau}+ \ep^{3- \frac{\beta}{2}- 2\tau }.
 \end{align*}
\end{Proof} 
Before we prove the next result, we recall a Lemma from \cite{Fefferman}.

\begin{Lemma}{\cite[Lemma 6.12]{Fefferman}}\label{fef_Ineqality} For all $f(\cdot) \in L^2(\R)$ and any fixed $\xi_0 \in \R $ we have
\begin{align*}
  \left\Vert \mathbbm{1}_{\left\{ \xi \,\in \,\ep^{\tau-\beta}\mathcal{B}\right\}}(\xi)\mathscr{F}[ f](\xi_0 + \ep^{\beta}\xi)\right\Vert_{L^2(\R_{\xi})} \lesssim \frac{1}{\ep^{\frac{\beta}{2}}} \Vert f\Vert_{L^2(\R)}.
\end{align*} 
\end{Lemma}
\begin{Proof}
 We first use  Plancherel Theorem to  get the bound
 \begin{align*}
 \frac{1}{(2\pi)^2}\int_{\left\{ \eta\,\in \,\xi_0 + \ep^{\tau}\mathcal{B}\right\}}\vert \widehat{f}(\eta)\vert^2 d\eta\leq \Vert f\Vert_{L^2(\R)}^2,  
 \end{align*}
followed by a change of variables $\eta = \xi_0 + \ep^{\beta}\xi $.
\end{Proof}

\begin{Lemma}[Nonlinear interaction Lemma II]\label{Lem:Nl_itr-2} Recall the choice of parameters \eqref{pmtr_range} of Proposition \ref{Prop:near_far_domination}, that is, $ 0 < \tau <\frac{1}{16}$,and $\beta \geq 1$.
Let $f, g \in H^4$ be given, with $\displaystyle{\Vert f \Vert_{H^4}= \mathcal{O}(\ep^{\frac{\beta}{2}})},$ $\displaystyle{\Vert f^2 \Vert_{H^4}= \mathcal{O}(\ep^{\frac{3\beta}{2}}) },$ and $\displaystyle{\Vert g \Vert_{H^4}= \mathcal{O}\left(\ep^{2 - 2\tau +  \frac{\beta}{2}} + \ep^{1 - 2\tau +  \frac{3\beta}{2}}+ \ep^{\frac{5\beta}{2}- 2\tau}+ \ep^{3- \frac{\beta}{2}- 2\tau } \right)}.$
Then
\begin{enumerate}[label=(\roman*), ref=\theTheorem(\roman*)]
\hangindent\leftmargin
 \item \label{Nl_ap_1}$  \max\{ \Vert f^2 g\Vert_{L^2}, \Vert f g^2 \Vert_{L^2}, \Vert g^3 \Vert_{L^2}\} = o(\ep^{2\beta}\ep^{\frac{\beta}{2}}).$
\item \label{Nl_ap_2}
  $\max\{ \Vert f g\Vert_{L^2}, \Vert g^2 \Vert_{L^2}\} = o(\ep^{2\beta-1}\ep^{\frac{\beta}{2}}).$
\end{enumerate}
  In particular, the result holds whenever  $f(\cdot) = v\n(\cdot)$ and $g(\cdot) = v\f(\cdot)$.
\end{Lemma}
\begin{Proof}
 We prove case by case, making repeated use of the  Sobolev Embedding \eqref{Sobolev_embedding}. First, we estimate each term in (i): 
 \begin{align*}
  \Vert f^2 g\Vert_{L^2} \lesssim \Vert f^2\Vert_{H^4}\Vert g\Vert_{H^4}&\lesssim \ep^{\frac{3\beta}{2}}\left(\ep^{2 - 2\tau +  \frac{\beta}{2}} + \ep^{1 - 2\tau +  \frac{3\beta}{2}}+ \ep^{\frac{5\beta}{2}- 2\tau}+ \ep^{3- \frac{\beta}{2}- 2\tau }  \right)\\
  &\lesssim \ep^{\frac{5\beta}{2}}\left(\ep^{2 - 2\tau - \frac{\beta}{2}} + \ep^{1 - 2\tau +  \frac{\beta}{2}}+ \ep^{\frac{3\beta}{2}- 2\tau}+ \ep^{3- \frac{3\beta}{2}- 2\tau }  \right),
 \end{align*}
which is $o(\ep^{2\beta})$, due to \eqref{pmtr_range}. Similarly, we have
\begin{align*}
  \Vert f g^2\Vert_{L^2} \lesssim \Vert f\Vert_{H^4}\Vert g\Vert_{H^4}^2  &\lesssim \ep^{\frac{\beta}{2}}\left(\ep^{2 - 2\tau +  \frac{\beta}{2}} + \ep^{1 - 2\tau +  \frac{3\beta}{2}}+ \ep^{\frac{5\beta}{2}- 2\tau}+ \ep^{3- \frac{\beta}{2}- 2\tau }  \right)^2\\
  &\lesssim \ep^{\frac{5\beta}{2}}\left(\ep^{4 - 4\tau -\beta} + \ep^{2 - 4\tau + \beta}+ \ep^{3\beta- 4\tau}+ \ep^{6- 3\beta- 4\tau }  \right),
 \end{align*}
which is also $o(\ep^{2\beta})$.  The last term is bounded as
\begin{align*}
 \Vert g^3 \Vert_{L^2}& \lesssim\left(\ep^{2 - 2\tau +  \frac{\beta}{2}} + \ep^{1 - 2\tau +  \frac{3\beta}{2}}+ \ep^{\frac{5\beta}{2}- 2\tau}+ \ep^{3- \frac{\beta}{2}- 2\tau }  \right)^3\\
 & \lesssim \ep^{\frac{5\beta}{2}}\left(\ep^{6 - 6\tau -\beta} + \ep^{3 - 6\tau +  2\beta}+ \ep^{5\beta -6\tau}+ \ep^{9 - 4b- 6\tau }  \right) = o(\ep^{2\beta}).
\end{align*}
Thus, (i) holds. To prove (ii), note that
\begin{align*}
 \Vert fg  \Vert_{L^2} & \lesssim 
 \ep^{\frac{5\beta}{2} -1}\left(\ep^{3 - 2\tau - \frac{3\beta}{2}} + \ep^{2 - 2\tau - \frac{\beta}{2}}+ \ep^{1+ \frac{\beta}{2}- 2\tau}+ \ep^{4- \frac{5\beta}{2}- 2\tau }  \right).
\end{align*}
Thanks to \eqref{pmtr_range}, we have $3 - \frac{3\beta}{2} - 2\tau = \left(3 - \frac{1}{4}- \frac{3\beta}{2}\right) + \left(\frac{1}{4} - 2\tau\right) >0$ and  $4 - \frac{5\beta}{2} - 2\tau = \left(4- \frac{1}{4}- \frac{5\beta}{2} - 2\tau\right) + \left( \frac{1}{4} - 2\tau\right) >0$. Therefore,  $\Vert fg  \Vert_{L^2}  = o(\ep^{2\beta -1}\ep^{\frac{\beta}{2}})$.

 Using \eqref{pmtr_range} once more, it is straightforward to verify that  $5 - 4\tau -3\frac{\beta}{2} = \left(2 - 2\tau - \frac{\beta}{2}\right) + \left(3 - 2\tau - \beta\right) >0 $ and  $7- 7\frac{\beta}{2}- 4\tau  = 2\left(\frac{1}{4} - 2\tau \right)  +\left(\frac{13}{2} - 7\frac{\beta}{2}\right)>0 $ holds. Hence,
\begin{align*}
 \Vert g^2  \Vert_{L^2} & \lesssim \left(\ep^{2 - 2\tau +  \frac{\beta}{2}} + \ep^{1 - 2\tau +  \frac{3\beta}{2}}+ \ep^{\frac{5\beta}{2}- 2\tau}+ \ep^{3- \frac{\beta}{2}- 2\tau }  \right)^2\\
 & \lesssim \ep^{\frac{5\beta}{2} -1}\left(\ep^{5 - 4\tau -3\frac{\beta}{2}} + \ep^{3 - 4\tau + \frac{\beta}{2}}+ \ep^{1 + \frac{5\beta}{2}- 4\tau}+ \ep^{7- 7\frac{\beta}{2}- 4\tau }  \right) =  o(\ep^{2\beta -1}\ep^{\frac{\beta}{2}}),
\end{align*}
and we are done.
 \end{Proof}
These results imply that we can use $v\n(\cdot)$ in the form \eqref{blow_up:near} to approximate the equation \eqref{mode_system_cutoff-a1}. The main result of this section is the following.

\begin{Proof}[of Proposition \ref{Prop:irrelevant:part2}]
 For simplicity, we shall write, $\mathscr{N}^{(j)}[v\n + v\f[v\n,\delta, \Omega,\gamma]]$ to denote $\mathscr{N}^{(j)}[v\n + v\f[v\n,\delta, \Omega,\gamma],  \ur, \ep]$.  To begin with, we apply Lemma \ref{fef_Ineqality}, getting
 \begin{align*}
  \left\Vert\mathbbm{1}_{\left\{ \xi \,\in \,\ep^{\tau-\beta}\mathcal{B}\right\}}(\xi) \mathscr{F}\left[\mathscr{N}^{(j)}[v\n + v\f]\right]\left(\pm 1  + \ep^{\beta} \xi\right)\right\Vert_{L^2(\R_{\xi})}\lesssim\ep^{-\frac{\beta}{2}} \Vert \mathscr{N}^{(j)}[v\n + v\f]\Vert_{L^2(\R)} =: Q_j.
 \end{align*}
 In the case $j=1$ we make use of  \eqref{Prop:nr_eng_dmn:far_bounded_by_near_improved_II}, of Lemma \ref{Lem:Nl_itr-2}, and of the linearity of $v\f(\cdot) \mapsto \mathscr{N}^{(j)}[v\f]$  to get
 \begin{align*}
Q_1&\lesssim \ep^{2-\frac{\beta}{2}} \Vert v\f\Vert_{L^2(\R)}   \lesssim \ep^{2-\frac{\beta}{2}}\left(\ep^{2 - 2\tau +  \frac{\beta}{2}} + \ep^{1 - 2\tau +  \frac{3\beta}{2}}+ \ep^{\frac{5\beta}{2}- 2\tau}+ \ep^{3- \frac{\beta}{2}- 2\tau } \right)\\
  &\lesssim \ep^{2\beta}\left(\ep^{4 - 2\tau -2\beta } + \ep^{3 - 2\tau-  \beta}+ \ep^{2- 2\tau}+ \ep^{5- 3\beta- 2\tau } \right),
 \end{align*}
 which is $o(\ep^{2\beta})$, due to the constraints in \eqref{pmtr_range}. The estimate for $j=2$ is a consequence of Lemma \ref{Lem:Nl_itr-2}, 
 \begin{align*}
   Q_2 & \lesssim\ep^{-\frac{\beta}{2}} \Vert \mathscr{N}^{(2)}[v\n + v\f]\Vert_{L^2(\R)}\\
   &\lesssim \ep^{1-\frac{\beta}{2}}\max\{ \Vert v\n v\f\Vert_{L^2}, \Vert v\f^2 \Vert_{L^2}, \Vert v\n^2 \Vert_{L^2}\}\\
   &= \ep^{1-\frac{\beta}{2}}\max\{o(\ep^{2\beta}\ep^{\frac{\beta}{2}-1}),  \ep^{\beta(\frac{3}{2})}\left(\Vert \gm^2\Vert_{L^2(\R)} + \Vert \gp^2 \Vert_{L^2(\R)}\right)\}\\ 
   &=\max\{o(\ep^{2\beta}), \ep^{1+ \beta}\left(\Vert \gm^2\Vert_{L^2(\R)} + \Vert \gp^2 \Vert_{L^2(\R)}\right)\}.
 \end{align*}
Finally, we use Lemma \ref{Lem:Nl_itr-2} to obtain the estimate when $j=3$:
\begin{align*}
  Q_3 & \lesssim\ep^{-\frac{\beta}{2}} \Vert \mathscr{N}^{(3)}[v\n + v\f]\Vert_{L^2(\R)}\\
  & \lesssim \ep^{-\frac{\beta}{2}}\max\{ \Vert v\n^2 v\f\Vert_{L^2}, \Vert v\n v\f^2 \Vert_{L^2}, \Vert v\f^3 \Vert_{L^2}, \Vert v\n^3 \Vert_{L^2}\} \\
  &= \ep^{-\frac{\beta}{2}}\max\{o(\ep^{2\beta}\ep^{\frac{\beta}{2}}), \ep^{\frac{5}{2}\beta}\left(\Vert \gm^3\Vert_{L^2(\R)} + \Vert \gp^3 \Vert_{L^2(\R)}\right)\}\\
  &= \max\{o(\ep^{2\beta}), \ep^{2\beta}\left(\Vert \gm^3\Vert_{L^2(\R)} + \Vert \gp^3 \Vert_{L^2(\R)}\right)\},
 \end{align*}
 and this finishes the proof. \end{Proof}
 
To conclude this section, we study the decay of $ v\n(x)$ and $v\f(x)$ as $\vert x\vert \to +\infty.$
 \begin{Proposition}[Decay of \texorpdfstring{$v\n(x)$}{v-near(x)} and \texorpdfstring{$v\f(x)$}{v-far(x)} as \texorpdfstring{$\vert x\vert \to +\infty$}{|x| ->infty}]\label{Prop:decay}
  For any fixed $\delta >0$ and $j \in {0,\ldots,3}$,  we have that 
  $$\lim_{\vert x\vert \to 0} \partial_x^{j}v\n(x) =  \lim_{\vert x\vert \to 0} \partial_x^{j}v\f(x) = 0.$$
 \end{Proposition}
 \begin{Proof}
  The result follows from classical Fourier analysis once we show that $\hv\n(\cdot)$, $\hv\f(\cdot)$ are $L^1(\R)$ functions (cf. \cite[Theorem 1.2]{Stein_Weiss}). In the case of $\hv\n(\cdot)$ this is straightforward: since the support of $\hv\n(\cdot)$ is bounded we can use the fact that $L_{loc}^1(\R) \subset L_{loc}^2(\R)$ to derive the result. 
  In the case  $v\f(\cdot)$, the conclude from the embedding $v\f(\cdot) \in H^1(\R)\subset H^4(\R)$ that $\lim_{|x|\to \infty} v(x) = 0$ (cf. \cite[Corollary 8.9]{Brezis}). Successive applications of this reasoning to $\partial_x^{j}v\f(\cdot)\in H^1(\R),$ for $j \in {1,2, 3}$ establishes the result.
  \end{Proof}


\section{Simplifications using a lemma of Fefferman, Thorpe and Weinstein, and matched asymptotics}\label{section:reduced_eq_simplifications}

Written as a system,  \eqref{second_main_eqt} reads as
\begin{equation}\label{second_main_eqt:rewrite}
 \left\{\begin{array}{l}
 - 4 \ep^{2\beta}\xi^2 \widehat{\gm}(\xi) - \widetilde{\mathbb{P}\n^{(-)}}\circ\mathscr{F}\left[\mathscr{N}^{(1)}[v\n , \ur^{(\delta, \om,\gamma)},\delta] +  \mathscr{N}^{(4)}[\ur^{(\delta, \om,\gamma)},\delta]\right](\ep^{\beta}\xi)  = -  \mathscr{A}_1^{(-)} + \mathscr{A}_2^{(-)}\\
  - 4 \ep^{2\beta}\xi^2 \widehat{\gp}(\xi) - \widetilde{\mathbb{P}\n^{(+)}}\circ\mathscr{F}\left[\mathscr{N}^{(1)}[v\n , \ur^{(\delta, \om,\gamma)},\delta] +  \mathscr{N}^{(4)}[\ur^{(\delta, \om,\gamma)},\delta]\right](\ep^{\beta}\xi)  = -  \mathscr{A}_1^{(+)} + \mathscr{A}_2^{(+)},  
 \end{array}\right. 
\end{equation}
where, for simplicity, we write  $\mathscr{A}_1^{(\pm)}$ and $\mathscr{A}_2^{(\pm)}$ to denote   $\mathscr{A}_1^{(\pm)}[v\f[v\n, \delta], \ep]$  and  $\mathscr{A}_2^{(\pm)}[v\f[v\n, \delta], \ep]$, respectively.  In virtue of \eqref{rolls}, we consider $\ep = \ep(\delta,\om)$.   Recall that the blow-up in (Fourier) parameter  \eqref{blow_up:near} gives the following structure to $v\n(\cdot)$:
$$v\n(x) = \ep^{\beta} e^{+i  x} \gp(\ep^{\beta} x) + \ep^{\beta} e^{-i x} \gm(\ep^{\beta} x).$$
Our goal from here to the end of the paper is solving  \eqref{mode_system_cutoff-a1}. An insatisfactory aspect of this equation lies on its dependence on \textit{artificially introduced parameters}: $\beta$, $\chi(\cdot)$ and $\tau$ . We mitigate this issue in this section,  finding a  fixed value for $\beta$. We shall adopt the following criterion:
\begin{Definition}[Parameter selection using matched asymptotics]\label{Def:matching}
  Within the parameter range \eqref{pmtr_range}, $\beta$ is a valid matching parameter whenever  the quantity
  \begin{align}\label{Def:matching_limit}
   \lim_{\ep \downarrow 0}\left(\frac{- 4 \ep^{2\beta}\xi^2 \widehat{g}_{\pm1}(\xi) - \widetilde{\mathbb{P}\n^{(\pm)}}\circ\mathscr{F}\left[\mathscr{N}^{(1)}[v\n , \ur^{(\delta, \om,\gamma)},\delta] +  \mathscr{N}^{(4)}[\ur^{(\delta, \om,\gamma)},\delta]\right](\ep^{\beta}\xi)}{\ep^{2\beta}}\right) 
  \end{align}
exists and is finite in the $L^2(\R)$ topology for   $v\n(\cdot)=v\n[\gm(\cdot),\gp(\cdot)]$ represented as in \eqref{blow_up:near} and any fixed $g_{\pm1}(\cdot)$ in $H^2(\R)$. 
\end{Definition}
In this section we show that $\beta =1$ is the only valid parameter we can use. Furthermore, we show  that the limit in Def. \ref{Def:matching} not only exists, but it is also computable quantity, a fact to which most of this section is devoted to. 
\subsection{Some auxiliary lemmas on two-scale interactions}
Initially we state some results that will help us to understand the multiple scales nature of the problem. We aim to simplify \eqref{second_main_eqt:rewrite}, in which the most relevant terms consist of the non-homogeneous and the linear ones. We shall write
\begin{align}\label{red_eqt:to_be_simplified}
\mathcal{I}^{(j; \pm)}(\xi)  = \widetilde{\mathbb{P}\n^{(\pm)}}\circ\mathscr{F}\left[ \mathscr{N}^{(j)}[v\n,\ur,\ep]\right](\ep^{\beta}\xi),\quad  \text{for} \quad j \in \{1, 4 \}. 
\end{align}
We decompose and split these terms in several pieces, whose analysis are fundamental to us. The following multi-scale result, which we see as   a  far/near (``slow/fast'') scales interaction estimate,  is one of our main tools.
\begin{Lemma}{\cite[Lemma 6.5 and Lemma 6.11]{Fefferman}}\label{Lemma:fef_pois}
 Let $f(x,\xi)$ and $g(x)$  denote smooth functions of $(x,\xi)\in \R\times \R$ that are 1-periodic in x. Let $\Gamma(x, X)$ be defined for $(x, X)$ and such that the two next conditions hold:
 $$\Gamma(x +1, X) = \Gamma(x, X), \quad \sum_{j=0}^2 \int_0^1\Vert \partial_x^j \Gamma(x,X) \Vert_{L^2(\R_X)}^2\mathrm{d}x <+\infty.$$
 Denote by $\widehat{\Gamma}(x,\omega)$ its Fourier transform with respect to the $X$ variable. Then, 
 \begin{align}\label{poisson_summation}
  \frac{\theta}{2\pi}\int_{\R}g(x) \Gamma(x,\theta x)\overline{f(x,\theta \xi)}e^{-i\theta \xi x}\mathrm{d}x = \sum_{n\in \mathbb{Z}}\int_0^1 e^{2\pi i n x} \widehat{\Gamma}\left(x,\frac{2\pi n}{\theta} + \xi\right)\overline{f(x,\theta \xi)}g(x)\mathrm{d}x
 \end{align}
Assume further that 
\begin{align}\label{assump:fef_pois}
C_f := \sup_{0\leq x\leq 1, \vert\omega\vert \leq \theta^{\tau}}\vert f(x,\omega)\vert< \infty,\quad  D_g := \Vert g \Vert_{L^{\infty}[0,1]} <\infty, \quad \text{and}\quad \Big\Vert \sup_{0\leq x\leq 1} \widehat{\Gamma}(x,\zeta) \Big\Vert_{H^1(\R_{\zeta})} < \infty. 
\end{align}
Define 
 $\displaystyle{\mathscr{I}_n(\xi, \theta) :=\int_0^1 e^{2\pi i n x} \widehat{\Gamma}\left(x,\frac{2\pi n}{\theta} + \xi\right)\overline{f(x,\theta \xi)}g(x)\mathrm{d}x}.$
Then, the following bounds hold
\begin{subequations}
 \begin{align*}
\left\Vert \mathbbm{1}( \vert\xi\vert \leq \theta^{\tau -1})\sum_{\vert n\vert \geq 1}\mathscr{I}_n(\xi;\theta)\right\Vert_{L^2(\R_{\xi})} &\lesssim C_f D_g \theta\Big\Vert \sup_{0\leq x\leq 1} \widehat{\Gamma}(x,\zeta) \Big\Vert_{H^1(\R_{\zeta})}, \\
\left\Vert \mathbbm{1}( \vert\xi\vert \leq \theta^{\tau -1})\sum_{ n\in \mathbb{Z}} \mathscr{I}_n(\xi;\theta)\right\Vert_{L^2(\R_{\xi})} &\lesssim C_f D_g\Big\Vert \sup_{0\leq x\leq 1} \widehat{\Gamma}(x,\zeta) \Big\Vert_{H^1(\R_{\zeta})}. 
 \end{align*}
\end{subequations}
\end{Lemma}
We adapt the previous  result  to our context  in the following manner:
\begin{Corollary}\label{Cor:fef_pois}  Let  constants, $f(x,\xi)$, $g(x)$,   $\Gamma(x, X)$ be as in the previous Lemma.  
Assuming  (\hyperref[H1]{H1})-(\hyperref[H2]{H2}) and $\tau$ and $\beta$ as in  and satisfying the constraints \eqref{pmtr_range}. Then the following bounds hold
\begin{subequations}
 \begin{align*}
\left\Vert \mathbbm{1}( \vert\xi\vert \leq \theta^{\tau -\beta})\sum_{\vert n\vert \geq 1}\mathscr{I}_n(\xi;\theta)\right\Vert_{L^2(\R_{\xi})} &\lesssim C_f D_g \theta\Big\Vert \sup_{0\leq x\leq 1} \widehat{\Gamma}(x,\zeta) \Big\Vert_{H^1(\R_{\zeta})}, \\
\left\Vert \mathbbm{1}( \vert\xi\vert \leq \theta^{\tau -\beta})\sum_{ n\in \mathbb{Z}} \mathscr{I}_n(\xi;\theta)\right\Vert_{L^2(\R_{\xi})} &\lesssim C_f D_g\Big\Vert \sup_{0\leq x\leq 1} \widehat{\Gamma}(x,\zeta) \Big\Vert_{H^1(\R_{\zeta})}. 
 \end{align*}
\end{subequations}
\end{Corollary}
\begin{Proof}
The proof is a slight modification of the proof in \cite{Fefferman}, using the fact that, whenever $\vert\xi \vert \leq \theta^{\tau - \beta}$ and $\vert\theta\vert <1$ we have
\begin{align*}
\left\vert\frac{2\pi n}{\theta} + \xi \right \vert \gtrsim \frac{\vert n\vert}{\theta}, \quad \forall \vert n\vert \geq 1,\quad \text{and}\quad 1 + \left\vert\frac{2\pi n}{\theta} + \xi \right \vert^2 \gtrsim 1 + \vert n\vert^2,\quad \forall n\in \mathbb{Z}.
\end{align*}
\end{Proof}
Before we embark into  more calculations, we derive another useful result, whose proof is straightforward.
\begin{Lemma}\label{Lem:simplf_cos}
  Recall from Lemma \ref{scaling:u_roll} that  
  \begin{align*}
  \ur^{(\delta,\om,\gamma)}(\cdot) = \ep \widetilde{\ur^{(\delta,\om,\gamma)}}(x),\quad \text{for} \quad \widetilde{\ur^{(\delta,\om,\gamma)}}(\cdot) = \cos(x+\gamma) + \ep^2 h(x); 
  \end{align*}
   where $x\mapsto h(x)$ is a $2\pi$-periodic $L^2$ mapping. Consider $g(x)$ and $\Gamma(x, X)$  as in the previous Lemma. Then, for any $k, p\in \mathbb{N}$ we have
 \begin{align*}
  \int_{0}^1e^{ip2\pi x}\left(\widetilde{\ur^{(\delta,\om,\gamma)}}(2\pi x)\right)^k\mathrm{d}x  =  \int_{0}^1e^{ip2\pi x}\left(\cos(2\pi x +\gamma)\right)^k\mathrm{d}x  + \mathcal{O}(\ep^2).
 \end{align*}  
\end{Lemma}

\begin{Observation}
In the next section  the following identities will be used several times:
 \begin{align*}
  \int_0^1\cos^4(2\pi z) \mathrm{d} z = \frac{3}{8}, \quad  \int_0^1\cos^2(2\pi z) \mathrm{d} z = \frac{1}{2}, \quad  \int_0^1\cos^2(2\pi z)\cos(4\pi z) \mathrm{d} z = \frac{1}{4}.
 \end{align*} 
\end{Observation}

\subsection{Simplifying \texorpdfstring{$\mathcal{I}^{(4; \, \pm)}(\xi)$}{I [+-1](4)(xi)}, or ``when we finally choose \texorpdfstring{$\beta=1$}{b=1}''}\label{sec:5_2}
In what follows, we analyze $\mathcal{I}^{(4; \, +)}(\xi)$; we omit the analysis of the case $\mathcal{I}^{(4; \, -)}(\xi)$, which is similar. 
The dependence of this term in $\ep$ is twofold: first due to the scalings $\mathbbm{1}_{\left\{ \xi \,\in \,\ep^{\tau-\beta}\mathcal{B}\right\}}(\ep^{\beta} \xi)$; and second, due to the term $\ur(\cdot) = \ep\widetilde{\ur}(\cdot)$. We initially break the interaction as
\begin{align*}
\mathcal{I}^{(4; \, +)}(\xi)  &= \widetilde{\mathbb{P}\n^{(+)}}\circ\mathscr{F}\left(   \mathscr{N}^{(4)}[v\n,\ur^{(\delta,\om,\gamma)},\ep]\right)(\ep^{\beta}\xi)\nonumber\\
 & = \widetilde{\mathbb{P}\n^{(+)}}\circ\mathscr{F}\left(  \chi( \chi^2 -1)\left(\ur^{(\delta,\om,\gamma)}\right)^3 \right)(\ep^{\beta}\xi)  + \widetilde{\mathbb{P}\n^{(+)}}\circ\mathscr{F}\left( [(1+ \partial_x^2)^2,\chi]\ur^{(\delta,\om,\gamma)} \right)(\ep^{\beta}\xi)\\
 & \quad + \widetilde{\mathbb{P}\n^{(+)}}\circ\mathscr{F}\left(\delta^2(\mu -1)\chi \ur^{(\delta,\om,\gamma)} \right)(\ep^{\beta}\xi)\nonumber\\
 & = \mathcal{S}_{\text{I}}^{(4; \, +)}(\xi) + \mathcal{S}_{\text{II}}^{(4; \, +)}(\xi) + \mathcal{S}_{\text{III}}^{(4; \, +)}(\xi),
\end{align*} 
which we analyze separately. We rewrite the fist term as
\begin{align*}
 \mathcal{S}_{\text{I}}^{(4; \, +)}(\xi) & = \mathbbm{1}_{\left\{ \xi \,\in \,\ep^{\tau-\beta}\mathcal{B}\right\}}(\ep^{\beta}\xi)\int_{\R}\chi(\ep^{\beta} x)( \chi^2(\ep^{\beta} x) -1)\left(\ur^{(\delta,\om,\gamma)}\right)^3(x) e^{-i(1 + \ep^{\beta} \xi)x }\mathrm{d}x \nonumber \\
 & \stackrel{x= 2\pi z, }{=}2\pi\ep^3\mathbbm{1}_{\left\{ \xi \,\in \,\ep^{\tau-\beta}\mathcal{B}\right\}}(\ep^{\beta}\xi)\int_{\R}\chi(2\pi\ep^{\beta} z)( \chi^2(2\pi \ep^{\beta} z) -1)\widetilde{\ur^{(\delta,\om,\gamma)}}^3(2\pi z) e^{-i(1 + \ep^{\beta} \xi)2\pi z }\mathrm{d}z \nonumber \\
 & \stackrel{\theta = 2\pi \ep^{\beta} }{=} \frac{4\pi^2\ep^3}{\theta}\mathbbm{1}_{\left\{ \xi \,\in \,\ep^{\tau-\beta}\mathcal{B}\right\}}(\ep^{\beta}\xi)\left(\frac{\theta}{2\pi}\int_{\R}\chi(\theta z)( \chi^2(\theta z) -1)(\widetilde{\ur^{(\delta,\om,\gamma)}})^3(2\pi z) e^{-i2\pi z} e^{-i \theta\xi z }\mathrm{d}z \right).
\end{align*}
Setting $Z = \theta z$ and writing $\Gamma_1(z, Z) = \Gamma_1( Z) =\chi(Z)( \chi^2(Z) -1)$, $f(z, \xi) = (\widetilde{\ur^{(\delta,\om,\gamma)}})^3(2\pi z)$, $g(z) = e^{-i2\pi z}$, we apply Lemma \ref{Lemma:fef_pois} and its Corollary \ref{Cor:fef_pois}, getting
\begin{align*}
 \mathcal{S}_{\text{I}}^{(4; \, +)}(\xi) & = \frac{4\pi^2\ep^3}{\theta}\mathbbm{1}_{\left\{ \xi \,\in \,\ep^{\tau-\beta}\mathcal{B}\right\}}(\ep^{\beta}\xi)\left[\widehat{\Gamma_1}\left( \xi\right)\int_{0}^1g(z)\overline{f(z)}\mathrm{d} z  + \sum_{\vert n\vert\geq 1}\int_{0}^1 e^{2\pi i n z}\widehat{\Gamma_1}\left(\frac{2\pi n}{\theta} + \xi\right)g(z)\overline{f(z)}\mathrm{d} z\right] \\
 & = \frac{4\pi^2\ep^3}{\theta}\mathbbm{1}_{\left\{ \xi \,\in \,\ep^{\tau-\beta}\mathcal{B}\right\}}(\ep^{\beta}\xi)\widehat{\Gamma_1}\left( \xi\right)\int_{0}^1g(z)\overline{f(z)}\mathrm{d} z  + \widetilde{F}_{\text{I}},
\end{align*}
with $\Vert \widetilde{F}_{\text{I}}\Vert_{L^2(\R)} = \mathcal{O}\left( \ep^3\right).$ Now we make a change of variables  and use the computations in Lemma \ref{Lem:simplf_cos},
\begin{align*}
 \mathcal{S}_{\text{I}}^{(4; \, +)}(\xi) & = \frac{4\pi^2\ep^3}{\theta}\mathbbm{1}_{\left\{ \xi \,\in \,\ep^{\tau-\beta}\mathcal{B}\right\}}(\ep^{\beta}\xi)\widehat{\Gamma_1}\left( \xi\right)\int_{0}^1e^{-2\pi i z}\cos^3(2\pi z + \gamma)\mathrm{d} z  + F_{\text{I}} \\
 & = \frac{4\pi^2\ep^3}{\theta}\mathbbm{1}_{\left\{ \xi \,\in \,\ep^{\tau-\beta}\mathcal{B}\right\}}(\xi)\widehat{\Gamma_1}\left( \xi\right)e^{i\gamma}\int_{0}^1\cos^4(2\pi z)\mathrm{d} z + F_{\text{I}} \\
  & = \frac{3\pi^2 \ep^3}{2 \theta}\mathbbm{1}_{\left\{ \xi \,\in \,\ep^{\tau-\beta}\mathcal{B}\right\}}(\xi)e^{i\gamma}\widehat{\Gamma_1}\left( \xi\right) + F_{\text{I}},
\end{align*}
where $\Vert F_{\text{I}}\Vert_{L^2(\R)}= \mathcal{O}\left(\frac{\ep^4}{\theta} + \ep^3\right).$ Substituting back $\theta = 2\pi \ep^{\beta}$ we get 
\begin{align}\label{s_4_plus_reduction}
 \mathcal{S}_{\text{I}}^{(4; \, +)}(\xi)  & = \frac{3\pi^2 \ep^3}{4\pi \ep^{\beta}}\mathbbm{1}_{\left\{ \xi \,\in \,\ep^{\tau-\beta}\mathcal{B}\right\}}(\xi)e^{i\gamma}\widehat{\Gamma_1}\left( \xi\right) +F_{\text{I}},
\end{align}
where the remainder has order $\Vert F_{\text{I}}\Vert_{L^2(\R)} = \mathcal{O}\left(\ep^{4-\beta} + \ep^3\right)$.  Note that $\widehat{\Gamma_1}(\cdot) \in L^2(\R)$, due to (\hyperref[H1]{H1}).

The terms $\mathcal{S}_{\text{II}}^{(4; \, +)}(\xi)$ and $\mathcal{S}_{\text{III}}^{(4; \, +)}(\xi)$ have similar derivations, and can be found in the appendix \ref{appendix:comput-4}. We put all together to derive a final form to $\mathcal{I}^{(4; \, \pm)}(\xi)$, which reads as
\begin{align}\label{bounds_I_4}
 \mathcal{I}^{(4; \, \pm)}(\xi)   = \mathbbm{1}_{\left\{ \xi \,\in \,\ep^{\tau-\beta}\mathcal{B}\right\}}(\xi)\left\{\frac{3 \pi \ep^{3-\beta}}{4}e^{\pm i \gamma}\widehat{\Gamma_1}\left( \xi\right) - 4\pi \ep^{1+\beta}e^{ \pm i \gamma}\widehat{\Gamma_2}( \xi) - 2\delta^2 \ep^{1-\beta}\widehat{\Gamma_3^{(\pm)}}(\xi) \right\}  +\mathscr{M}^{(4; \, \pm)}(\xi),
\end{align}
where $\Gamma_{2}(\cdot), \Gamma_{3}^{(\pm)}(\cdot)\in L^2(\R)$ have explicit expressions given in \ref{appendix:comput-4} and further 
\begin{align*}
\left \Vert\mathscr{M}^{(4; \, \pm)}(\xi) \right\Vert_{L^2(\R)} = \mathcal{O}(\ep^{6+\beta}+\ep^{6+3\beta}+\ep^{5\beta}+ \ep^{2(1 + 2\beta)} + \ep^{2(5-\beta)}), 
\end{align*}
with $\mathscr{M}^{(4; \, \pm)}(\xi) = \mathbbm{1}_{\left\{ \xi \,\in \,\ep^{\tau-1}\mathcal{B}\right\}}(\xi)\mathscr{M}^{(4; \, \pm)}(\xi)$.%

At this point we have enough information to choose a value for $\beta$; the reasoning lies in the following result:

\begin{Lemma}[Matched asymptotics I]\label{Lemma:beta_is_necess_1} Assume $\beta$ satisfying \eqref{pmtr_range}. In order to the limit \eqref{Def:matching_limit} exist it is necessary to have $\beta =1.$
\end{Lemma}
\begin{Proof} Since the limit in Definition \eqref{Def:matching} must exist for all $ g_{\pm1}(\cdot)\in  H^2(\R)$ we can set $g_{\pm1}(\cdot)\equiv 0$. We must show that, as $\ep\downarrow0$, the limit
\begin{align}\label{beta_necesseary_eqt}
 \frac{ \widetilde{\mathbb{P}\n^{(+)}}\circ\mathscr{F}\left[ \mathscr{N}^{(4)}[\ur^{(\delta, \om,\gamma)},\delta]\right](\ep^{\beta}\xi)}{\ep^{2\beta}} = \frac{\mathcal{S}_{\text{I}}^{(4; \, +)}(\xi) + \mathcal{S}_{\text{II}}^{(4; \, +)}(\xi)+\mathcal{S}_{\text{III}}^{(4; \, +)}(\xi)}{\ep^{2\beta}},
\end{align}
exists, being the treatment the same in the case $\frac{ \widetilde{\mathbb{P}\n^{(-)}}\circ\mathscr{F}\left[ \mathscr{N}^{(4)}[\ur^{(\delta, \om,\gamma)},\delta]\right](\ep^{\beta}\xi)}{\ep^{2\beta}}$. In what follows we show that
\begin{align}\label{beta_dichotomy}
\lim_{\ep\downarrow 0}\left \Vert  \frac{\mathcal{S}_{\text{I}}^{(4; \, +)}(\xi) + \mathcal{S}_{\text{II}}^{(4; \, +)}(\xi)+\mathcal{S}_{\text{III}}^{(4; \, +)}(\xi)}{\ep^{2\beta}}
\right\Vert_{L^2(\R)} =\left\{\begin{array}{ll}                                                                                                                +\infty, &  \quad \text{when} \quad \beta >1,\\
\text{exists}, & \quad  \text{when} \quad \beta =1,                                                                                                               \end{array}\right.
 \end{align}
which then asserts the necessity of $\beta =1$. Considering \eqref{bounds_I_4}, and defining 
$\psi(\xi):= \frac{3 \pi \ep^{3-\beta}}{4}e^{\pm i \gamma}\widehat{\Gamma_1}\left( \xi\right) - 4\pi \ep^{1+\beta}e^{ \pm i \gamma}\widehat{\Gamma_2}( \xi) - 2\delta^2 \ep^{1-\beta}\widehat{\Gamma_3^{(\pm)}}(\xi)$ we rewrite \eqref{beta_necesseary_eqt} as
\begin{align*}
 \frac{ \widetilde{\mathbb{P}\n^{(+)}}\circ\mathscr{F}\left[ \mathscr{N}^{(4)}[\ur^{(\delta, \om,\gamma)},\delta]\right](\ep^{\beta}\xi)}{\ep^{2\beta}} = \frac{\mathbbm{1}_{\left\{ \xi \,\in \,\ep^{\tau-\beta}\mathcal{B}\right\}}(\xi)\psi(\xi)}{\ep^{2\beta}}.
\end{align*}
As (\hyperref[H1]{H1}) implies that $\psi(\cdot) \in L^2(\R)\cap \mathscr{C}^{\infty}(\R)$ and $\Gamma_3(\cdot)$ discontinuous, we can conclude that $\psi(\cdot) \not \equiv 0$. Therefore, for a sufficiently large $c>0$, we have $\displaystyle{\Vert \mathbbm{1}_{c\mathcal{B}}(\cdot) \psi(\cdot) \Vert > 0 }.$ 
 Thanks to the choice of parameters in \eqref{pmtr_range}, for a sufficiently small $\ep>0$ we have that $c\mathcal{B} \subset \ep^{\tau - \beta}\mathcal{B},$ for the latter set exhaust the real line as $\ep\downarrow 0$. We now use this fact to do the following minoration:
 \begin{align*}
  \left \Vert  \frac{\mathcal{S}_{\text{I}}^{(4; \, +)}(\xi) + \mathcal{S}_{\text{II}}^{(4; \, +)}(\xi)+\mathcal{S}_{\text{III}}^{(4; \, +)}(\xi)}{\ep^{2\beta}}
\right\Vert_{L^2(\R)} \gtrsim \left \Vert \frac{\mathbbm{1}_{\left\{ \xi \,\in \,\ep^{\tau-\beta}\mathcal{B}\right\}}(\xi)\psi(\xi)}{\ep^{2\beta}}\right\Vert_{L^2(\R)}\gtrsim   \left \Vert \frac{\mathbbm{1}_{c\mathcal{B}}(\xi)\psi(\xi)}{\ep^{2\beta}} \right \Vert.
 \end{align*}
Inspecting the coefficient of each term on the last term we have 
\begin{align*}
 \frac{\mathbbm{1}_{c\mathcal{B}}(\xi)\psi(\xi)}{\ep^{2\beta}} &= \mathbbm{1}_{c\mathcal{B}}(\xi)\left(\frac{3 \pi \ep^{3-3\beta}}{4}e^{\pm i \gamma}\widehat{\Gamma_1}\left( \xi\right) - 4\pi \ep^{1-\beta}e^{ \pm i \gamma}\widehat{\Gamma_2}( \xi) - 2\delta^2 \ep^{1-3\beta}\widehat{\Gamma_3^{(\pm)}}(\xi) \right)\\
 & = \ep^{3-3\beta}\mathbbm{1}_{c\mathcal{B}}(\xi)\left(\frac{3 \pi }{4}e^{\pm i \gamma}\widehat{\Gamma_1}\left( \xi\right) - 4\pi \ep^{2\beta-2}e^{ \pm i \gamma}\widehat{\Gamma_2}( \xi) - 2\frac{\delta^2}{\ep^2}\widehat{\Gamma_3^{(\pm)}}(\xi) \right)\\
 & = \ep^{3-3\beta}\mathbbm{1}_{c\mathcal{B}}(\xi)\widetilde{\psi}(\xi),
\end{align*}
where $\widetilde{\psi}(\cdot) = \frac{\psi(\cdot)}{\ep^{3-\beta}} \in L^2(\R)$ is a non-zero term that is uniformly bounded in $L^2(\R)$ for all $\ep\geq 0$ sufficiently small. Hence, given the constraint \eqref{pmtr_range}, then the left hand side of \eqref{beta_dichotomy} blows-up as $\ep\downarrow 0$ whenever $\beta >1$. 

It remains to prove the  case $\beta =1$:  recall that $\mathbbm{1}_{\left\{ \xi \,\in \,\ep^{\tau-1}\mathcal{B}\right\}}(\xi) \uparrow1$; thus, invoking  $(\hyperref[H2]{H2})$, namely,  $\displaystyle{\lim_{\ep\downarrow0} \frac{\delta^2}{\ep^2} =1}$, and the Lebesgue Dominated Convergence Theorem, we conclude that $\frac{\mathcal{S}_{\text{I}}^{(4; \, +)}(\cdot) + \mathcal{S}_{\text{II}}^{(4; \, +)}(\cdot)+\mathcal{S}_{\text{III}}^{(4; \, +)}(\cdot)}{\ep^{2}}$ converges to $\frac{3 \pi}{4}e^{\pm i \gamma}\widehat{\Gamma_1}\left( \cdot\right) - 4\pi e^{ \pm i \gamma}\widehat{\Gamma_2}( \cdot) - 2\widehat{\Gamma_3^{(\pm)}}(\cdot)$ in $L^2(\R)$ as $\ep \downarrow 0$, and this finishes the proof.
\end{Proof}

Henceforth we shall take $\beta=1$, devoting the rest of this section to prove that the limit in Definition \ref{Def:matching} holds true in this case. To conclude the analysis of $\mathcal{I}^{(4; \, \pm)}(\xi)$, we use computations in appendix \ref{appendix:comput-4}.

\subsection{The term \texorpdfstring{$\mathcal{I}^{(1; \, \pm)}(\xi)$}{I [+-1](1)(xi) }}\label{sec:5_3}
Calculations in the appendix \ref{appendix:comput-1} show that%
\begin{align*}
\mathcal{I}^{(1; \, \pm)}(\xi)& =-\mathbbm{1}_{\left\{ \ep^{\tau-1}\mathcal{B}\right\}}(\xi) \mathscr{F}\left[\left(\delta^2\mu(x) - 3 (\widetilde{\ur^{(\delta,\om,\gamma)}})^2\right) v\n(\cdot)\right](1 + \ep\xi)
\end{align*}
gets simplified to 
\begin{align}\label{bounds_I_1}
 \mathcal{I}^{(1; \, \pm)}(\xi) & = \mathbbm{1}_{\left\{ \ep^{\tau-1}\mathcal{B}\right\}}(\xi)\left\{ -\delta^2  \int_{\R} \mu(X)\left[ g_{\pm 1}(X) \right] e^{-i  \xi X}\mathrm{d}X  +  3\pi \ep^2\widehat{g_{\pm 1}}(\xi) + 3\frac{\pi\ep^2}{2}\widehat{g_{\mp 1}}(\xi)\right\} +\mathscr{M}^{(1; \, \pm)}(\xi),
\end{align}
where $\mathscr{M}^{(1; \, \pm)}(\xi) = \mathbbm{1}_{\left\{ \xi \,\in \,\ep^{\tau-1}\mathcal{B}\right\}}(\xi)\mathscr{M}^{(1; \, \pm)}(\xi)$, with bounds $\left\Vert \mathscr{M}^{(1; \, \pm)}(\xi) \right \Vert_{L^2(\R)} =  \mathcal{O}\left(\ep^3\right)$.

With this result in hands, we finally conclude the matched asymptotics argument:
\begin{Proposition}[Matched asymptotics II]\label{beta_is_1:sufficiency}
 In order to the limit in Definition \ref{Def:matching} to exist, it is necessary and  sufficient to take $\beta =1$.
\end{Proposition}
\begin{Proof} Necessity has been proven in Lemma \ref{Lemma:beta_is_necess_1}. The other terms in \eqref{Def:matching_limit} are given by \eqref{bounds_I_1}, which is clearly of order $\mathcal{O}(\ep^2)$ and the term $\xi^2\widehat{g}_{\pm 1}$, which is bounded in $L^2(\R)$ due to $g_{\pm 1}(\cdot) \in H^2(\R)$. This establishes the result. 
\end{Proof}
In summary, with $\beta=1$ we have obtained
\begin{equation}\label{almost_all_together}
\begin{split}
&\widetilde{\mathbb{P}\n^{(\pm)}}\circ\mathscr{F}\left[\mathscr{N}^{(1)}[v\n , \ur^{(\delta, \om,\gamma)},\delta] +  \mathscr{N}^{(4)}[\ur^{(\delta, \om,\gamma)},\delta]\right](\ep\xi) \\ 
&\stackrel{\eqref{red_eqt:to_be_simplified}}{=}   \mathcal{I}^{(1; \pm)}(\xi) +  \mathcal{I}^{(4; \pm)}(\xi)\\
&\stackrel{\eqref{bounds_I_1}\&\eqref{bounds_I_4}}{=}
\mathbbm{1}_{\left\{ \ep^{\tau-1}\mathcal{B}\right\}}(\xi)\left\{ -\delta^2  \int_{\R} \mu(X)\left[ g_{\pm 1}(X) \right] e^{-i  \xi X}\mathrm{d}X  +  3\pi \ep^2\widehat{g_{\pm 1}}(\xi) + 3\frac{\pi\ep^2}{2}\widehat{g_{\mp 1}}(\xi)\right\} +
\\&\quad\quad\quad \quad+ \mathbbm{1}_{\left\{ \xi \,\in \,\ep^{\tau-1}\mathcal{B}\right\}}(\xi)\left\{\frac{3 \pi \ep^{2}}{4}e^{\pm i \gamma}\widehat{\Gamma_1}\left( \xi\right) - 4\pi \ep^{2}e^{ \pm i \gamma}\widehat{\Gamma_2}( \xi) - 2\delta^2 \widehat{\Gamma_3^{(\pm)}}(\xi) \right\} \\
&\quad \quad \quad\quad +\mathscr{M}^{(1; \, \pm)}(\xi) +\mathscr{M}^{(4; \, \pm)}(\xi).
\end{split} 
\end{equation}

\subsection{Putting things together}\label{sec:5_4}
Plugging \eqref{almost_all_together} into \eqref{second_main_eqt:rewrite}, we  divide it by $\ep^2 >0$, obtaining a  system of \emph{reduced equations},
\begin{align}\label{second_main_eqt:simplified}
 \widetilde{\mathcal{R}^{(\delta,\Omega,\gamma)}}[\widehat{\gm}, \widehat{\gp}] - \mathbbm{1}_{\left\{ \ep^{\tau-1}\mathcal{B}\right\}}(\xi)h_{*}(\xi) = \widetilde{\mathcal{Q}^{( \delta,\Omega,\gamma;\pm)}}[\widehat{\gm}, \widehat{\gp}];
\end{align}
where $\ep = \ep(\delta,\om)\stackrel{(\hyperref[H2]{H2})}{=}\ep\left(\delta,\sqrt{1+ \delta\Omega}\right)$. We carefully explain each term: the left hand side depends linearly on $\gm(\cdot)$ and $\gp(\cdot)$, namely, 

\begin{equation*}
\begin{split}
 \widetilde{\mathcal{R}^{(\delta,\Omega,\gamma)}}[\widehat{\gm}, \widehat{\gp}]&:= \left( \begin{array}{c}
                                                                                   \widetilde{\mathcal{R}^{(\delta,\Omega,\gamma;-)}}[\widehat{\gm}, \widehat{\gp}]\\
                                                                                   \widetilde{\mathcal{R}^{(\delta,\Omega,\gamma;+)}}[\widehat{\gm}, \widehat{\gp}]
                                                                                  \end{array}\right) 
\end{split},
\end{equation*}
with
\begin{align*}
\widetilde{\mathcal{R}^{(\delta,\Omega;\pm)}}[\widehat{\gm}, \widehat{\gp}]=  - 4 \xi^2 \widehat{g_{\pm 1}}(\xi)  - \mathbbm{1}_{\left\{ \ep^{\tau-1}\mathcal{B}\right\}}(\xi)\left[  3\pi \widehat{g_{\pm 1}}(\xi)  + \frac{3 \pi}{2}\widehat{g_{\mp 1}}(\xi) - \mathscr{F}\left[\mu(\cdot)g_{\pm 1}(\cdot)\right](\xi)\right].
\end{align*}
Since the operator $\widetilde{\mathcal{R}^{(\delta,\Omega,\gamma)}}[\widehat{\gm}, \widehat{\gp}] $ is a multiplier, we can define its associated operator $\mathcal{R}^{(\delta,\Omega,\gamma)}[g\m, \gp]$ in  physical space  as
\begin{align}\label{Obs:operator_in_physical_space}
 \mathcal{R}^{(\delta,\Omega,\gamma)}[g\m, \gp](x) := \mathscr{F}^{-1}\left[\widetilde{\mathcal{R}^{(\delta,\Omega,\gamma)}}[\widehat{\gm}, \widehat{\gp}]\right](x),
\end{align} 
an operator with domain $\mathcal{D} \left(\mathcal{R}^{(\delta,\Omega,\gamma)}[\cdot, \cdot] \right) = H^2(\R) \times H^2(\R)$. On  the right hand side are non-homogeneous and nonlinear terms, 
\begin{align*}
 \widetilde{\mathcal{Q}^{(\delta,\Omega,\gamma)}}[\widehat{\gm}, \widehat{\gp}] +h_*(\xi) = \left(\begin{array}{c}
                                                                                   \widetilde{\mathcal{Q}^{(\delta,\Omega,\gamma;-)}}[\widehat{\gm}, \widehat{\gp}]\\
                                                                                   \widetilde{\mathcal{Q}^{(\delta,\Omega,\gamma;+)}}[\widehat{\gm}, \widehat{\gp}]
                                                                                  \end{array}\right),
\end{align*}
with
\begin{align*}
\widetilde{\mathcal{Q}^{( \delta,\Omega,\gamma;\pm)}}[\widehat{\gm}, \widehat{\gp}] &= -  \mathscr{A}_1^{(\pm)} + \mathscr{A}_2^{(\pm)} + \frac{\mathscr{M}^{(1; \, \pm)}(\xi) + \mathscr{M}^{(4; \, \pm)}(\xi)}{\ep^2},
\end{align*}
and
\begin{align}\label{h_star}
h_*(\xi) = \left(\begin{array}{c}
h_{*}^{(-)}(\xi)\\
h_{*}^{(+)}(\xi)
\end{array}\right) =
\mathbbm{1}_{\left\{ \ep^{\tau-1}\mathcal{B}\right\}}(\xi)\left(\begin{array}{c}
\frac{3 \pi }{4}e^{-i \gamma}\widehat{\Gamma_1}\left( \xi\right) - 4\pi e^{-i \gamma}\widehat{\Gamma_2}( \xi) - 2\frac{\delta^2}{\ep^2} \widehat{\Gamma_3^{(-)}}(\xi) \\
\frac{3 \pi }{4}e^{+i \gamma}\widehat{\Gamma_1}\left( \xi\right) - 4\pi e^{+i \gamma}\widehat{\Gamma_2}( \xi) - 2\frac{\delta^2}{\ep^2} \widehat{\Gamma_3^{(+)}}(\xi) 
\end{array}\right).
\end{align}
We remark that  $h_{*}(\xi) = \mathbbm{1}_{\left\{ \ep^{\tau-1}\mathcal{B}\right\}}(\xi)h_*(\xi)$. Thanks to Proposition \ref{Prop:irrelevant}, estimates \eqref{bounds_I_4} and \eqref{bounds_I_1}, we have 
\begin{align}\label{small_nonlinear_terms}
\Vert\widetilde{\mathcal{Q}^{( \delta,\Omega,\gamma;\pm)}}[\widehat{\gm}, \widehat{\gp}] \Vert_{L^2(\R)}  = o\left(1 \right)+\mathcal{O}\left(\Vert \gm\Vert_{L^2(\R)}^2 +\Vert \gp\Vert_{L^2(\R)}^2\right).
\end{align}
%

\section{Approximation and solvability of the reduced equation: final steps} \label{sec:reduced_equation+approx+solv} 
All the previous simplifications led us so far to \eqref{second_main_eqt:simplified}, what we have called a system of  \emph{reduced equations},
\begin{align*}
 \widetilde{\mathcal{R}^{(\delta,\Omega,\gamma)}}[\widehat{\gm}, \widehat{\gp}] - \mathbbm{1}_{\left\{ \ep^{\tau-1}\mathcal{B}\right\}}(\xi)h_{*}(\xi) = \widetilde{\mathcal{Q}^{( \delta,\Omega,\gamma;\pm)}}[\widehat{\gm}, \widehat{\gp}].
\end{align*}
The main result of this section shows that this system admits a nontrivial solution  $\left(\gm(\cdot), \gp(\cdot)\right)$,  parametrized  by $(\delta,\Omega, \gamma)$.  After this preamble, we are able to state the one of the main result of this section, which consists of another application of the Lyapunov-Schmidt reduction method; Theorem \ref{Thm:main} will be a obtained as a direct consequence of it. 

\begin{Proposition}\label{Prop:reduced} Assume (\hyperref[H1]{H1})-(\hyperref[H2]{H2}) and  a fixed parameter $\tau$ satisfying the constraints  \eqref{pmtr_range}   ($\beta=1$ also fixed, due to Proposition \ref{beta_is_1:sufficiency}). Then, there exists a small $\delta_{*}>0$ and a continuous mapping
$$ (0,\delta_{*})\times \left(-\delta_*,\delta_*\right)\times \T \ni (\delta,\Omega,\gamma) \mapsto \left(\gm^{(\delta,\Omega,\gamma)}(\cdot),\gp^{(\delta,\Omega,\gamma)}(\cdot)\right) \in  H_{\text{near}, \ep^{\tau-1}}^{2}(\R)\times H_{\text{near}, \ep^{\tau-1}}^{2}(\R),$$
which is $2\pi$-periodic in $\gamma$ and satisfies \eqref{second_main_eqt:simplified}, hence \eqref{second_main_eqt:rewrite}. These functions are band-limited, namely,
\begin{align*}
g_{\pm 1}^{(\delta,\Omega,\gamma)}(\cdot) \in H^2(\R), \quad \text{and further,} \quad  \mathrm{supp}\left(\widehat{g_{\pm 1}^{(\delta,\Omega,\gamma)}}\right) \subset \ep^{\tau -1}\mathcal{B},
\end{align*}
where $\ep=\ep(\delta,\Omega)$.
\end{Proposition}
In passing, combining the previous result with Lemma \ref{GL:representation} and the extension to $\delta =0$ given by Proposition \ref{pmtr_range:extension_to_0}, we immediately obtain the next Corollary.
\begin{Corollary}\label{Prop_6-1-corollary}
For a sufficiently small $\delta_{*}= \delta_{*}(\tau)>0$, the mapping
\begin{align}\label{parametrizations-near_by_delta}
 [0,\delta_{*})\times \left(-\delta_*,\delta_*\right) \times \T\ni (\delta,\Omega,\gamma) \mapsto v\n^{(\delta,\Omega,\gamma)}(\cdot) \in X_{\text{near}, \ep^{\tau}}^{4}\subset H^4(\R).
\end{align}
is continuous. Furthermore, allying this result to that of Proposition \ref{Prop:near_far_domination} implies continuity of
\begin{align}\label{parametrizations-far_by_delta}
 [0,\delta_{*})\times \left(-\delta_*,\delta_*\right)\times \T \ni (\delta,\Omega,\gamma) \mapsto v\f^{(\delta,\Omega,\gamma)}(\cdot) =  v\f[ \delta,\Omega,v\n^{(\delta,\Omega,\gamma)}](\cdot) \in X_{\text{far}, \ep^{\tau}}^{4} \subset H^4(\R),
 \end{align}
where we write
$\om= \sqrt{1+ \delta \Omega}$
thanks to (\hyperref[H2]{H2}).
\end{Corollary}
\begin{Remark}\label{Rmk:digression} We make a brief digression before we tackle this problem; our discussion is similar to that in \cite[Remark 6.3]{Fefferman}. First, recall  that we seek for band limited functions $\gm(\cdot)$, $\gp(\cdot)$ as a solution to this problem. Thus, it is natural to investigate the invertibility of the mapping $(\widehat{\gm}(\cdot), \widehat{\gp}(\cdot))\mapsto  \widetilde{\mathcal{R}^{(\delta,\Omega,\gamma)}}[\widehat{\gm}, \widehat{\gp}]: H_{\text{near}, \ep^{\tau-1}}^{2}(\R)\times H_{\text{near}, \ep^{\tau-1}}^{2}(\R)\to L^{2}(\R).$ An approach to this problem could go along the following line of reasoning: the parameter choice in \eqref{pmtr_range} implies that  $\displaystyle{\lim_{\ep \downarrow 0}\mathbbm{1}_{\left\{ \xi \,\in\, \ep^{\tau- 1}\mathcal{B}\right\}}(\xi)\uparrow1}$. In this fashion, as we take $\ep\downarrow 0$ and use $ \Vert\widetilde{\mathcal{Q}^{( \delta,\Omega,\gamma;\pm)}}[\widehat{\gm}, \widehat{\gp}] \Vert_{L^2(\R)}  = o\left(\ep^2 \right)$, we \textit{formally} obtain on the left hand side
\begin{align}
\lim_{\ep \downarrow 0}\widetilde{\mathcal{R}^{(\delta,\Omega,\gamma)}}[\widehat{\gm}, \widehat{\gp}] \approx \widetilde{\mathcal{R}^{(0,\Omega,\gamma)}}[\widehat{\gm}, \widehat{\gp}]:= \left( \begin{array}{c}
                                                                                   \widetilde{\mathcal{R}^{(0,\Omega;-)}}[\widehat{\gm}, \widehat{\gp}]\\
                                                                                   \widetilde{\mathcal{R}^{(0,\Omega;+)}}[\widehat{\gm}, \widehat{\gp}]
                                                                                  \end{array}\right),
\end{align}
where $\approx$ should read as a \textit{formal} limit.\footnote{It must be highlighted that the operator $\widetilde{\mathcal{R}^{(0,\Omega,\gamma)}}[\widehat{\gm}, \widehat{\gp}]$ in fact does not depend on $\Omega$. In spite of the risk of being misleading,  we kept the notation in this form for consistency.} The operator $\widetilde{\mathcal{R}^{(0,\Omega;\pm)}}$ can be written explicitly as
\begin{align}\label{Operator_formal_limit}
\widetilde{\mathcal{R}^{(0,\Omega;\pm)}}[\widehat{\gm}, \widehat{\gp}]=  - 4 \xi^2 \widehat{g_{\pm 1}}(\xi)  -  3\pi \widehat{g_{\pm 1}}(\xi)  - \frac{3 \pi}{2}\widehat{g_{\mp 1}}(\xi) + \mathscr{F}\left[\mu(\cdot)g_{\pm 1}(\cdot)\right](\xi). 
\end{align}
It is illuminating and worth to reflect upon this erroneous approach to solve \eqref{second_main_eqt:simplified}:  plugging $\widetilde{\mathcal{R}^{(0,\Omega;\pm)}}[\widehat{\gm}, \widehat{\gp}]$ in place of $\widetilde{\mathcal{R}^{(\delta,\Omega;\pm)}}[\widehat{\gm}, \widehat{\gp}]$, the system to be solved becomes
\begin{align*}
\widetilde{\mathcal{R}^{(0,\Omega,\gamma)}}[\widehat{\gm}, \widehat{\gp}] - \mathbbm{1}_{\left\{ \ep^{\tau-1}\mathcal{B}\right\}}(\xi)h_{*}(\xi) = \widetilde{\mathcal{Q}^{( \delta,\Omega,\gamma;\pm)}}[\widehat{\gm}, \widehat{\gp}].
\end{align*}
Thus, we would ideally solve this system in $(\gm, \gp) \in H^2(\R)\times H^2(\R)$ to, posteriorly, truncate the solution in Fourier space (for we want $(\gm, \gp) \in H_{\text{near}, \ep^{\tau-1}}^{2}\times H_{\text{near}, \ep^{\tau-1}}^{2}$). 
The issue we stumble upon with this argument is that, for some $f(\cdot) \in L^2(\R)\times L^2(\R)$, we have
\begin{align*}
\mathbbm{1}_{\left\{ \xi \,\in\, \ep^{\tau- 1}\mathcal{B}\right\}}(\xi)\left(\left(\widetilde{\mathcal{R}^{(0,\Omega,\gamma)}}\right)^{-1}[f]\right)(\xi) \neq \left(\widetilde{\mathcal{R}^{(0,\Omega,\gamma)}}\right)^{-1}\left[\mathbbm{1}_{\left\{ \xi \,\in\, \ep^{\tau- 1}\mathcal{B}\right\}}(\cdot)f(\cdot)\right](\xi).
\end{align*}
A more effective strategy, the one that we shall pursue next, still uses  $ \widetilde{\mathcal{R}^{(0,\Omega,\gamma)}}[\widehat{\gm}, \widehat{\gp}] $ to approximate the operator $ \widetilde{\mathcal{R}^{(\delta,\Omega,\gamma)}}[\widehat{\gm}, \widehat{\gp}] $.
\end{Remark}
We appeal to the formal construction of the operator $(\widehat{\gm}, \widehat{\gp})\mapsto \widetilde{\mathcal{R}^{(0,\Omega,\gamma)}}[\widehat{\gm}, \widehat{\gp}]$ to rewrite the system \eqref{second_main_eqt:simplified} in an equivalent form as
\begin{align}\label{Approximated_inverse:problem}
\widetilde{\mathcal{R}^{(0,\Omega,\gamma)}}[\widehat{\gm}, \widehat{\gp}] + \left(\widetilde{\mathcal{R}^{(\delta,\Omega,\gamma)}} - \widetilde{\mathcal{R}^{(0,\Omega,\gamma)}}\right)[\widehat{\gm}, \widehat{\gp}] - \mathbbm{1}_{\left\{ \xi \,\in\, \ep^{\tau- 1}\mathcal{B}\right\}}(\xi)h_{*}(\xi)= \widetilde{\mathcal{Q}^{( \delta,\Omega,\gamma;\pm)}}[\widehat{\gm}, \widehat{\gp}].
\end{align}
Three ingredients are involved in finding nontrivial solutions to the above system of equations: 
\begin{enumerate}[label=(\roman*), ref=\theTheorem(\roman*)]
\item  proving that the operator  $\widetilde{\mathcal{R}^{(\delta,\Omega,\gamma)}}[\widehat{\gm}, \widehat{\gp}]$ is well approximated by  $\widetilde{\mathcal{R}^{(0,\Omega,\gamma)}}[\widehat{\gm}, \widehat{\gp}]$, and that the latter is invertible; 
\item  showing that solutions to problem \eqref{Approximated_inverse:problem} are indeed band-limited; 
\item  showing that $h_{*}(\cdot)$ can be taken with a ``small'' $\Vert \cdot \Vert_{L^2(\R)}$ norm.
\end{enumerate}
We follow this agenda in the next Lemmas.
\begin{Lemma}[$\mathcal{R}^{(0,\Omega,\gamma)}$ is a good, invertible, approximation]\label{Lemma:reduced_equation:coercivity}  Let $\mathcal{R}^{(\delta,\Omega,\gamma)}$ be defined as in \eqref{Obs:operator_in_physical_space}, while we assume that 
$\mathcal{R}^{(0,\Omega,\gamma)}$ is defined as the limiting operator \eqref{Operator_formal_limit}, with  $\mathcal{D}\left(\mathcal{R}^{(0,\Omega,\gamma)}\right) = H^2(\R)\times H^2(\R)\subset L^2(\R)\times L^2(\R) \to L^2(\R)\times L^2(\R).$ Then, the following properties  hold:
 \begin{enumerate}
 \item[(i)] \label{reduced_equation:coercivity-1} There exists a constant $C>0$ (independent of $\ep,\delta,\Omega$) such that 
 \begin{align}
  \left(\mathcal{R}^{(0,\Omega,\gamma)}\right)^{-1}: L^2(\R)\times L^2(\R) \to H^2(\R)\times H^2(\R), \quad \left\Vert  \left(\mathcal{R}^{(0,\Omega,\gamma)}\right)^{-1} f\right\Vert_{ H^2(\R)\times H^2(\R)} \leq C \Vert f\Vert_{ L^2(\R)\times L^2(\R)}.
 \end{align}

  \item[(ii)]\label{reduced_equation:coercivity-2} Assume that $g_{\pm 1}(\cdot) \in H^2(\R)$. Then
\begin{align*}
\left\Vert \left(\mathcal{R}^{(\delta,\Omega,\gamma )} - \mathcal{R}^{(0,\Omega,\gamma)}\right)[\gm,\gp]\right\Vert_{L^2(\R)\times L^2(\R)} \lesssim \ep^{1 - \tau}\Vert g_{\pm 1} \Vert_{H^2(\R)}. 
\end{align*}
\item[(iii)] \label{reduced_equation:coercivity-3} The following inequality holds:
\begin{align*}
 \left\Vert (\mathcal{R}^{(0,\Omega,\gamma)})^{-1}\left(\mathcal{R}^{(\delta,\Omega,\gamma )} - \mathcal{R}^{(0,\Omega,\gamma)}\right)\right\Vert_{H^2(\R)\times H^2(\R)\to H^2(\R)\times H^2(\R)} <1, \quad \ep \in (0, \ep_0).
\end{align*}
 \end{enumerate} 
 \end{Lemma}
\begin{Proof}
The inequality in (iii) is a direct consequence of the results in (i) and (ii):
\begin{align*}
 &\left\Vert (\mathcal{R}^{(0,\Omega,\gamma)})^{-1}\left(\mathcal{R}^{(\delta,\Omega,\gamma )} - \mathcal{R}^{(0,\Omega,\gamma)}\right)\right\Vert_{H^2(\R)\times H^2(\R)\to H^2(\R)\times H^2(\R)}\\ 
 &\quad  \lesssim \left\Vert (\mathcal{R}^{(0,\Omega,\gamma)})^{-1}\right\Vert_{ {L^2(\R)\times L^2(\R)\to H^2(\R)\times H^2(\R)}}  \left\Vert\mathcal{R}^{(\delta,\Omega,\gamma )} - \mathcal{R}^{(0,\Omega,\gamma)}\right\Vert_{H^2(\R)\times H^2(\R)\to L^2(\R)\times L^2(\R)} \\
 & \quad \lesssim   \ep^{1- \tau} \Vert (\mathcal{R}^{(0,\Omega,\gamma)})^{-1}\Vert_{L^2(\R)\times L^2(\R)\to H^2(\R)\times H^2(\R)}.
\end{align*}
In order to prove (i) we first observe that the operator $v\mapsto \mathcal{R}^{(0,\Omega,\gamma)}[v]$ described in \eqref{Operator_formal_limit}  converges to the following endstates 
\begin{align*}
  v\mapsto \mathcal{R}_{ (x = -\infty)}^{(0,\Omega,\gamma )}[v] =  \left(\begin{array}{c}
                                                                     4 \partial_x^2 v_1  - \left( 3\pi   + 1\right)v_1 - 3\frac{\pi}{2}v_2 \\
                                                                     4 \partial_x^2 v_2   - \left( 3\pi   + 1\right)v_2 - 3\frac{\pi}{2}v_1 
                                                                    \end{array}\right), \qquad (\text{as} \, x \to - \infty);\\
  v\mapsto \mathcal{R}_{ (x = +\infty)}^{(0,\Omega,\gamma )}[v]=  \left(\begin{array}{c}
                                                                     4 \partial_x^2 v_1  - \left( 3\pi   - 1\right)v_1 - 3\frac{\pi}{2}v_2 \\
                                                                     4 \partial_x^2 v_2   - \left( 3\pi   - 1\right)v_2 - 3\frac{\pi}{2}v_1 
                                                                    \end{array}\right),\qquad (\text{as} \, x \to + \infty).
\end{align*}
We claim that both these constant coefficient operators are coercive. Indeed, whenever $f = (f_1,f_2)\in H^2(\R)\times H^2(\R),$
\begin{align*}
  \left\langle   \mathcal{R}_{ (x = -\infty)}^{(0,\Omega,\gamma )} f, f \right \rangle_{L^2(\R)} & \approx \left\langle \mathscr{F}\left[  \mathcal{R}_{ (x = -\infty)}^{(0,\Omega,\gamma )} f\right],\mathscr{F}\left[ f\right] \right \rangle_{L^2(\R)} \\
 & \gtrsim 4\left( \Vert \partial_x f_1\Vert_{L^2(\R)}^2 + \Vert \partial_x f_2\Vert_{L^2(\R)}^2\right)+ \left(3\pi^2\right)\left(\Vert f_1\Vert_{L^2(\R)}^2+\Vert f_2\Vert_{L^2(\R)}^2\right)\\
 & \quad + 3\pi\mathrm{Re}\left(\int_{\R}f_1(x)\overline{f_2}(x)\mathrm{d}x\right)\\ 
  & \gtrsim  \left( 3\pi^2- \frac{3\pi}{2} - 1\right)\Vert f(\cdot)\Vert_{H^1(\R)\times H^1(\R)}^2,
 \end{align*}
and the result follows. The same holds for  $\mathcal{R}_{ (x = +\infty)}^{(0,\Omega,\gamma )}$. Applying the Lax-Milgram Theorem (cf. \cite[Chapter 2, Section 9]{Lions_Magenes-I}) we assert the invertibility of the asymptotic operators $\displaystyle{\mathcal{R}_{ (x = \pm\infty)}^{(0,\Omega )}:H^2(\R) \to L^2(\R)  }$. Hence, arguing as in \cite[\S 3]{Salamon} (or similarly, as in \cite[Proposition 4.3]{Monteiro_Scheel-contact_angle}) we conclude that the operator $\displaystyle{ \mathcal{R}^{(0,\Omega,\gamma)}: H^2(\R) \to L^2(\R) }$ is a Fredholm operator with index 0. In order to prove the invertibility of this operator it suffices to show that its kernel and cokernel are trivial. In both cases these reasoning is the same:  a similar  calculation shows that the operator and its adjoint   are coercive operators, and we conclude that these are trivial spaces. Hence, $\mathrm{Ker}\left(\mathcal{R}^{(0,\Omega,\gamma)}\right)=\mathrm{coKer}\left(\mathcal{R}^{(0,\Omega,\gamma)}\right)= \{0\}$, and the operator  $v\mapsto \mathcal{R}^{(0,\Omega,\gamma)}[v]$ is invertible. 

With regards to property (ii),  we first note that
 \begin{align*}
  \left(\mathcal{R}^{(\delta,\Omega,\gamma )} - \mathcal{R}^{(0,\Omega,\gamma)}\right)\widehat{g_{\pm 1}}(\xi) =  \left(1 - \mathbbm{1}_{\left\{ \ep^{\tau-1}\mathcal{B}\right\}}(\xi)\right)\left(\begin{array}{cc}
										- \left[  \chi\left(\frac{x}{\ep}\right) + 3\pi^2\right]& -\frac{3\pi}{2}\\
										  -\frac{3\pi}{2}&  - \left[  \chi\left(\frac{x}{\ep}\right) + 3\pi^2\right]
										\end{array}\right)
										\left(\begin{array}{c}
		    \widehat{\gp}(\xi)\\
		    \widehat{\gm}(\xi)
		      \end{array}
\right)
 \end{align*}
 Now notice that 
 \begin{align*}
  \left\Vert \left(1 - \mathbbm{1}_{\left\{ \ep^{\tau-1}\mathcal{B}\right\}}(\xi)\right) \widehat{\gp}(\xi)\right\Vert_{L^2(\R)}^2 &\lesssim \int_{\vert\xi\vert > \ep^{\tau -1}} \vert \widehat{\gp}(\xi) \vert^2 \mathrm{d}\xi
   \lesssim \int_{\vert\xi\vert > \ep^{\tau -1}}\frac{1}{\vert\xi\vert^2} \vert\xi\vert^2 \vert \widehat{\gp}(\xi) \vert^2 \mathrm{d}\xi 
   \lesssim  \ep^{2(1-\tau)} \Vert g_{\pm 1} \Vert_{H^2(\R)}^2,
\end{align*}
which concludes the proof of (ii).
\end{Proof}
The following result validates our approach, showing that we can avoid the commutativity issues alluded to in Remark \ref{Rmk:digression}.
\begin{Lemma}[Persistence of band-limited properties under approximation]\label{Lemma:persistence} Any functions $\gm(\cdot), \gp(\cdot) \in H^{2}(\R)$ satisfying \eqref{Approximated_inverse:problem} are in $H_{\text{near}, \ep^{\tau-1}}^{2}(\R)$, that is,  
$$g_{\pm 1}(\cdot) \in H^2(\R), \quad \text{and further,} \quad  \mathrm{supp}\left(\widehat{g}_{\pm 1}\right) \subset \ep^{\tau -1}\mathcal{B}.$$
\end{Lemma}
\begin{Proof}
 It is clear that  both systems are equivalent. Hence, one can rewrite system \eqref{Approximated_inverse:problem} back to the form \eqref{second_main_eqt:simplified} and multiply it by $ 1 - \mathbbm{1}_{\left\{ \xi \,\in \,\ep^{\tau-1}\mathcal{B}\right\}}(\xi)$, obtaining
 \begin{align*}
4 \left( 1 - \mathbbm{1}_{\left\{ \xi \,\in \,\ep^{\tau-1}\mathcal{B}\right\}}(\xi)\right)\xi^2\mathscr{F} \left(\begin{array}{c}
 g\m(\cdot)\\
 g\p(\cdot)
\end{array}\right)(\xi) = \left(\begin{array}{c}
			    0\\
			    0
			  \end{array}\right)\end{align*}
			  which readily implies that $\mathscr{F}[g_{\pm1}](\cdot)= \widehat{g_{\pm 1}}(\cdot)$ are both supported in $\ep^{\tau -1}\mathcal{B}.$
\end{Proof}
\begin{Lemma}[Attaining small norms for  $h_{*}(\cdot)$ by choosing  $\chi(\cdot)$]\label{Lemma:almost_partition}
 Recall $h_*(\cdot) = (h_*^{(-)},h_*^{(+)})$ from \eqref{h_star}, where
 \begin{align*}
  h_*^{(\pm)}(\xi)= \frac{3 \pi }{4}e^{\pm i \gamma}\widehat{\Gamma_1}\left( \xi\right) - 4\pi e^{\pm i \gamma}\widehat{\Gamma_2}( \xi) - 2\frac{\delta^2}{\ep^2} \widehat{\Gamma_3^{(\mp)}}(\xi).
 \end{align*}
For any given number $c_0>0$, we can choose $\chi(\cdot)$  the latter satisfies the ODE
\begin{align}\label{ODE}
 -4\pi \partial_x^2 \chi(\cdot) +\frac{3 \pi}{4}\left( \chi^3(\cdot) - \chi(\cdot)\right) - c\partial_x\chi(x) =0, \quad \lim_{x\to -\infty}\chi(x) = 1, \quad \lim_{x\to -\infty}\chi(x) = 0,
\end{align}
for any  $0 <c <c_0$. Furthermore, this choice can be done in  such a way that following bounds hold,
 \begin{align*}
\Vert h_*^{(\pm)}(\cdot)\Vert_{L^2(\R)\times L^2(\R)}\leq c_0.
 \end{align*}
 In particular, we can choose $\chi(\cdot)$ so that it satisfies (\hyperref[H1]{H1}).
 \end{Lemma}
\begin{Proof} As shown in appendix \ref{existence_TW}, the existence of a solution $\chi(\cdot)$ to the ODE above satisfying the asserted spatial asymptotics in \eqref{ODE} is known.  In this manner, choosing a $\chi(\cdot)$ with this quality, we can exploit  the explicit formulas for $\Gamma_{1,2}$ (see \S \ref{sec:5_2}) to get
 \begin{align*}
  \frac{3 \pi}{4}\widehat{\Gamma_1}(\xi)  - 4\pi\widehat{\Gamma_2}( \xi)  = \mathscr{F}\left[-4\pi \partial_x^2 \chi(\cdot) +\frac{3 \pi}{4}\left( \chi^3(\cdot) - \chi(\cdot)\right)\right](\xi) = \mathscr{F}\left[c\partial_s\chi(\cdot)\right](\xi).
 \end{align*} 
Thus, one can write
\begin{align*}
h_*(\xi) =  \mathbbm{1}_{\left\{ \ep^{\tau-1}\mathcal{B}\right\}}(\xi)\mathscr{F}\left[c\partial_s\chi(\cdot)\right](\xi) - 2\frac{\delta^2}{\ep^2} \widehat{\Gamma_3}(\xi). 
\end{align*}
We shall estimate both terms on the right hand side. Beginning with the first one: it is shown in appendix \ref{existence_TW}, we have that $c\partial_x\chi(\cdot)\in L^2(\R)$, and $\displaystyle{\Vert c\partial_x\chi(\cdot)\Vert_{L^2(\R)}^2\lesssim \vert c\vert}$. Therefore, applying Plancherel's identity successively, we obtain
\begin{align*}
  \Vert \mathscr{F}\left[c\partial_x \chi(\cdot) \right]\Vert_{L^2(R)}^2  \lesssim \Vert c\partial_x \chi(\cdot) \Vert_{L^2(R)}^2 \lesssim \vert c\vert <c_0,
\end{align*}
which can be made as small as necessary by taking $0 < c_0 \ll 1$. 
With regards to the second term, first recall that
$$\Gamma_3^{(\pm)}(\cdot) = \mathbbm{1}_{[0,\infty)}(\cdot) \chi\left(\cdot \right)e^{-i\frac{(\cdot)}{\ep}}\cos\left(\frac{\cdot}{\ep}\pm \gamma\right),$$
We begin by invoking  (\hyperref[H2]{H2}) to bound $\displaystyle{\frac{\delta^2}{\ep^2}\leq 16}.$ Now,  we exploit  the fact that solutions to the ODE \eqref{ODE} are  translation invariant and exponentially decaying to 0 as $x\to +\infty$: considering solutions to the above ODE of the form $\chi(\cdot + \tau_x)$, for any given $c_0>0$ we can choose  $\tau_x^* >0$ sufficiently large so that
\begin{align*}
\frac{1}{2\pi}\Vert \widehat{\Gamma_3^{(\pm)}}\Vert_{L^2(\R)}^2 =   \Vert \Gamma_3^{(\pm)}\Vert_{L^2(\R)}^2 \leq  \int_{0}^{+\infty}\vert\chi(x+ \tau_x)\vert^2\mathrm{d}x \leq c_0^2, \quad \text{for all} \quad \tau_x >\tau_x^*;
\end{align*}
combining both estimates, the result follows. \end{Proof} 
\begin{Remark}[On the use of far/near (spatial) decompositions for numerical analysis purposes]\label{rmk:numerical}
 The fact that the bounds obtained in Lemma \ref{Lemma:reduced_equation:coercivity} are independent of $\ep$ is essential for choosing  $\chi(\cdot)$ in Lemma \ref{Lemma:almost_partition}. It further corroborates with the idea that  far/near (spatial) decomposition can be used as an interesting tool in numerical analysis of bifurcation in extended domains. For other approaches, see \cite{Morrissey,Lloyd-Sc} or \S \ref{open_problems:numerical}.
\end{Remark}
We are ready to put  things together, combining  the previous Lemmas to give a proof Proposition \ref{Prop:reduced}:
\begin{Proof}[of Proposition \ref{Prop:reduced}] Due to the equivalence between \eqref{second_main_eqt:simplified} and  \eqref{Approximated_inverse:problem}, we can work with the latter.  Recall  that we can write $\ep = \ep(\delta, v\n(\cdot))$, as a continuous function of its parameters; throughout the analysis, $\ep$ should be considered in this form. Moreover, thanks to observations in (\hyperref[H2]{H2}) and Corollary \ref{rolls:reparametrization}, we can replace $\delta$ for  $\ep$ in the analysis, for they are equivalent. In the following, thanks for property (\hyperref[H2]{H2}), we shall use $\Omega$ instead of $\om = \om(\delta,\Omega)= \sqrt{1 + \delta \Omega}.$ Therefore, throughout the proof, we shall use $(\delta,\Omega,\gamma)$ instead of the parametrization $(\delta,\omega,\gamma)$.

Now we handle the rest:  we use the operator $\left(\widetilde{\mathcal{R}^{(0,\Omega,\gamma)}}\right)^{-1}$ obtained in Lemma \ref{Lemma:reduced_equation:coercivity} to  act on \eqref{Approximated_inverse:problem}, the equation we aim to find a solution for; we rewrite the outcome as
 \begin{align}\label{Approximated_inverse:problem-equiv}
&\text{Id}[\widehat{\gm}, \widehat{\gp}] =\mathcal{M}[\widehat{\gm}, \widehat{\gp},\delta,\Omega,\gamma],
\end{align}
where we denote the identity mapping as $(\widehat{\gm}, \widehat{\gp}) \mapsto Id[\widehat{\gm}, \widehat{\gp}]$, and
\begin{align*}
 \mathcal{M}[\widehat{\gm}, \widehat{\gp},\delta,\Omega,\gamma]&:=  -\left(\widetilde{\mathcal{R}^{(0,\Omega,\gamma)}}\right)^{-1}\left[\left(\widetilde{\mathcal{R}^{(\delta,\Omega,\gamma)}}- \widetilde{\mathcal{R}^{(0,\Omega,\gamma)}}\right)[\widehat{\gm}, \widehat{\gp}] \right]\\
 &\quad +\left(\widetilde{\mathcal{R}^{(0,\Omega,\gamma)}}\right)^{-1} \left[\mathbbm{1}_{\left\{ \xi \,\in\, \ep^{\tau- 1}\mathcal{B}\right\}}(\xi)h_{*}(\xi)\right]\\
 &\quad +\left(\widetilde{\mathcal{R}^{(0,\Omega,\gamma)}}\right)^{-1}\left[\widetilde{\mathcal{Q}^{( \delta,\Omega,\gamma;\pm)}}[\widehat{\gm}, \widehat{\gp}]\right]\\
 &= C_{1} +C_{2}+C_{3}.
\end{align*}
Once more we rely on Lemma \ref{Lemma:reduced_equation:coercivity} in several ways: first, to assert that the term $C_{1}$ is $\mathcal{O}(\ep^{1-\tau})$; 
the second term $C_{2}$ can be made as small as we would want in virtue of both Lemma \ref{Lemma:almost_partition} and Lemma \ref{reduced_equation:coercivity-1}, for the bounds on the inverse mapping $\left(\widetilde{\mathcal{R}^{(0,\Omega,\gamma)}}\right)^{-1}$ are independent of $(\delta,\Omega,\gamma)$; for the last term $C_{3}$, we also rely on the same bounds, allied to the estimates \eqref{small_nonlinear_terms}. Thus, there exists a $r_*>0$ sufficiently small such that  the mapping  $(\widehat{\gm}, \widehat{\gp}) \mapsto \mathcal{M}[\widehat{\gm}, \widehat{\gp},\delta,\Omega,\gamma]$  is a (uniform) contraction and  maps the set
\begin{align}
\{\gm(\cdot), \gp(\cdot) \in H^2(\R) |\Vert \gm(\cdot)\Vert_{H^2(\R)} + \Vert \gm(\cdot)\Vert_{H^2(\R)} \leq r_{*} \}
\end{align}
into itself.  Thus, another application of the Contraction Mapping Theorem implies that there exists a solution to \eqref{Approximated_inverse:problem-equiv}, parametrized by $(\delta,\Omega,\gamma)$, that is, there exists a $0<\delta_*= \delta_{*}(\tau)\ll\frac{1}{3}$ for which
\begin{align}\label{Prop:reduced-continuity}
(0,\delta_{*})\times \left(-\delta_{*},\delta_{*}\right)\times \T \ni (\delta,\Omega,\gamma) \mapsto \left(\gm^{(\delta,\Omega,\gamma)}(\cdot),\gp^{(\delta,\Omega,\gamma)}(\cdot)\right), 
\end{align}
with
\begin{align*}
\text{Id}\left[\widehat{\gm^{(\delta,\Omega,\gamma)}}(\cdot),\widehat{\gp^{(\delta,\Omega,\gamma)}}(\cdot)\right] =\mathcal{M}\left[\widehat{\gm^{(\delta,\Omega,\gamma)}}(\cdot),\widehat{\gp^{(\delta,\Omega,\gamma)}}(\cdot),\delta\right],\quad \text{for all} \quad (\delta,\Omega,\gamma) \in (0,\delta_{*})\times \left(-\delta_*,\delta_*\right)\times \T.
\end{align*}
Finally, we invoke Lemma \ref{Lemma:persistence} to conclude that $\left(\gm^{(\delta,\Omega,\gamma)}(\cdot),\gp^{(\delta,\Omega,\gamma)}(\cdot)\right) \in H_{\text{near}, \ep^{\tau-1}}^{2}(\R)\times H_{\text{near}, \ep^{\tau-1}}^{2}(\R)$. 
Let us now study the regularity of these parametrizations: initially, notice that all the bounds obtained in Lemma \ref{Lemma:reduced_equation:coercivity} and Lemma \ref{Lemma:almost_partition} are uniform in $0<\ep\ll1$. Since we have pointwise convergence of these functions in $\ep$, using the Lebesgue Dominated Convergence Theorem  we conclude that
\begin{align*}
 (0,\delta_{*})\times \left(-\delta_*,\delta_*\right)\times \T \ni (\delta,\Omega,\gamma) \mapsto \mathcal{M}[\widehat{\gm}, \widehat{\gp},\delta,\Omega,\gamma] \in H^2(\R)\times H^2(\R)
\end{align*}
is a continuous  mapping for any fixed $\gm(\cdot)$, $\gp(\cdot)$. Another application of the Contraction Mapping Theorem implies the continuity of the mapping \eqref{Prop:reduced-continuity} (as in \cite[\S2, Theorem 2.2]{Hale}). We argue using the uniqueness of the fixed point and $2\pi$-periodicity of the mapping $\gamma \mapsto \mathcal{M}[\widehat{\gm}, \widehat{\gp},\delta,\Omega,\gamma]$ to derive the 2$\pi$-periodicity of $\widehat{\gm}, \widehat{\gp}$ on the same parameter,  and with this we conclude the proof.
\end{Proof}
%

\section{Wavenumber selection --  proof of Theorem \ref{Thm:main}}\label{section:wavenumber_selection}
According to the results in Corollary \ref{Prop_6-1-corollary} we write the Ansatz \eqref{Ansatz} as
\begin{align}\label{Ansatz-rewritten}
(x,\delta,\Omega,\gamma) \mapsto \mathcal{U}(x) = v\n^{(\delta,\Omega,\gamma)}(x) +v\f[\delta, \Omega,v\n^{(\delta,\Omega,\gamma)}](x)+ \chi(\ep x)\ur^{(\delta,\sqrt{1+ \delta\Omega},\gamma)}\left(x\right), 
\end{align}
which gives a solution to problem \eqref{SH-eq-bifurcation-1} with the qualities we were after, as stated in \ref{heterecolinic_def}; see Corollary \ref{Corollary:solution_to_the_problem} below.  It remains to be shown that we have found is indeed a solution in the sense of Definition \ref{heterecolinic_def}.
In this case, fixing $\gamma$, we have a full family of solutions parametrized by two parameters $(\delta,\Omega,\gamma)$. Now, we unfold all the reductions we performed, and plug $v(\cdot)$ back into problem \eqref{asymptotics_properties}. In this new  context, an interesting phenomenon happens: the conservation laws due to the  Hamiltonian structure of the problem impose a severe parameter restriction, namely, we obtain a selection mechanism, implying that only one parameter is necessary in the characterization of solutions to problem \eqref{SH-eq-bifurcation-1}. Moreover, these solutions have the qualities we are concerned with, as stated in Definition \ref{heterecolinic_def}.
\begin{Lemma}[Regularity]\label{Lemma:regularity} Writing $\ep = \ep(\delta,\Omega)$ and $v(\cdot) = v\n(\cdot) + v\f(\cdot)$ as before, the following regularity condition holds
$$x\mapsto \mathcal{U}(x) = v(x) + \chi(\ep x)\ur^{(\delta,\sqrt{1+ \delta\Omega},\gamma)}\left(x\right)\in\mathscr{C}^{(\infty)}\left(\R\setminus\{0\}\right)\cap \mathscr{C}^{(3)}\left(\R\right).$$ 
\end{Lemma}
\begin{Proof}
With regards to regularity, the embedding $H^4(\R) \hookrightarrow \mathscr{C}^3(\R)$ gives part of the result, while smoothness in $\R\setminus\{0\}$ (i.e., aways from the quenching-front) is proved using either elliptic regularity \cite[Chapter 2, \S 3.2]{Lions_Magenes-I}, or using  \cite[Corollary 3.1.6]{Hormander}. 
\end{Proof}
We can rewrite the constraint in \eqref{SH-eq-bifurcation-2} as
 \begin{align}\label{hamiltonian-know_constants}
- \delta^2u^2\Big\vert_{x=0} - \delta^2 = \mathcal{H}^{(l)}[\ur^{(\delta,\sqrt{1+ \delta\Omega},\gamma)}(\cdot)] - \mathcal{H}^{(r)}[0].
 \end{align}
 As we shall see next, this property is the main ingredient in the selection mechanism that relates $\delta$  to the wavenumber $\Omega$ of the rolls $\ur^{(\delta, \Omega)}$. But first we need to expand these functions in powers of $\delta$ in order to understand how they can be  approximated.
\begin{Lemma} \label{Lemma:hamiltonian_expansion} Assume (\hyperref[H2]{H2}), and choose  parameters  $\tau$ as in \eqref{pmtr_range} of Proposition \ref{Prop:near_far_domination}. Write the  Ansatz as in \eqref{Ansatz-rewritten}. 
Then,  the following approximations hold:
\begin{enumerate}[label=(\roman*), ref=\theTheorem(\roman*)]
\hangindent\leftmargin
\item \label{Lemma:hamiltonian_expansion:1}$\mathcal{U}^2(x)\Big|_{x=0} = o(\ep)= o(\delta).$
\item \label{Lemma:hamiltonian_expansion:2} $\mathcal{H}^{(l)}[\ur^{(\delta,\sqrt{1+ \delta\Omega},\gamma)}(\cdot)] -\mathcal{H}^{(r)}[0]=- \delta^2 +\frac{4}{3}\delta^3\left(\Omega + \Omega^3\right) + \mathcal{O}(\delta^3).$
\end{enumerate}
\end{Lemma}
\begin{Proof} To prove (i), we first note that
$$\mathcal{U}^2(0) \lesssim \max\left\{(\chi(0)\ur^{(\delta,\sqrt{1+ \delta\Omega},\gamma)}(0))^2, v\n^2(0), v\f^2(0)\right\}.$$
We must show that the right hand side in the above equation is $o(\ep)$ (hence, $o(\delta)$, thanks to (\hyperref[H2]{H2})). From  \eqref{rolls} we get $\displaystyle{\chi(0)\ur^{(\delta,\sqrt{1+ \delta\Omega},\gamma)}(0) = \mathcal{O}(\ep^2)}$. The second estimate comes from the fact that $v\n^p(\cdot)$ is a band limited function for any $p\in \mathbb{N}\setminus\{0\},$ hence $v\n^p(\cdot) \in H^s(\R)$ for all $s\geq1$. In particular, when $p=2$ we can make use of Remark \ref{embedding_band_limited} and apply the Sobolev embedding \eqref{Sobolev_embedding} to get
$$v\n^2(0) \leq \Vert v\n^2(\cdot)\Vert_{L^{\infty}(\R)}\lesssim \Vert v\n^2(\cdot)\Vert_{H^1(\R)} \lesssim \Vert v\n^2(\cdot)\Vert_{L^{2}(\R)}.$$
An application of \eqref{Nl_itr_inequ:1} of Lemma \ref{Nl_itr_Lm_1} with $\beta =1$ then gives $\displaystyle{v\n^2(0) \lesssim \mathcal{O}(\ep^{\frac{3}{2}}})= o(\ep) = o(\delta)$. For the last term, we use inequality \eqref{Prop:nr_eng_dmn:far_bounded_by_near_improved_II} with $\beta =1$, %
\begin{align*}
\vert v\f(0)\vert \leq \Vert v\f(\cdot)\Vert_{L^{\infty}(\R)} \lesssim \Vert v\f(\cdot) \Vert_{H^4(\R)} = \mathcal{O} \left(\ep^{\frac{5}{2} - 2\tau }\right).
\end{align*}
Choosing $\tau$ as in \eqref{pmtr_range},  the fact that $\delta \approx \ep$ allows us to conclude the result in (i).

To prove (ii) we expand $\mathcal{H}^{(l)}[\ur^{(\delta,\sqrt{1+ \delta\Omega},\gamma)}(\cdot)]$ using \eqref{rolls}, the fact that $\omega^2 = 1 + \delta \Omega$ from property (\hyperref[H2]{H2}), and the equivalence $\delta \approx \ep$. \end{Proof}
In principle, whenever $\delta =0$ we have that $\ur^{(\delta,\sqrt{1+ \delta\Omega},\gamma)}(\cdot) \equiv 0$, hence $\Omega $ would be allowed to take any value. It turns out that $\Omega $ can be chosen in an unique fashion if we extend it to the value it takes as $\delta \downarrow 0$; our approach is allusive to the technique used in \cite[\S 6.7]{Fefferman}.
\begin{Lemma}\label{Lemma:wavenumber_selec} Let  $\tau$ as in \eqref{pmtr_range} of Proposition \ref{Prop:near_far_domination}.  Recall the parametrization $\ep = \ep(\delta, \Omega)$, due to  (\hyperref[H2]{H2}). Let $\delta_*>0$ as in Proposition \ref{Prop:reduced}.  Consider the mapping
 \begin{align*}
  (0,\delta_*)\times \left(-\delta_*,\delta_*\right)\times \T \ni (\delta, \Omega,\gamma)\mapsto \mathcal{S}[\delta, \Omega,\gamma] : = \frac{- \delta^2\mathcal{U}^2\Big\vert_{x=0} - \delta^2 -  \mathcal{H}^{(l)}[\ur^{(\delta,\sqrt{1+ \delta\Omega},\gamma)}(\cdot)] + \mathcal{H}^{(r)}[0]}{\delta^2},
 \end{align*}
where $\mathcal{U}(x)$  is written as in \eqref{Ansatz-rewritten}.  Then, the following properties hold:
\begin{enumerate}[label=(\roman*), ref=\theTheorem(\roman*)]
\hangindent\leftmargin
\item \label{Prop:final_reduction-smoothness} The mapping $\mathcal{S}[\delta, \Omega,\gamma]:(0,\delta_*)\times \left(-\delta_*,\delta_*\right)\times \T \to \R $ is smooth;
\item \label{Prop:final_reduction-selection}(Selection mechanism) There exists a $\delta_{**}$ satisfying $0 <\delta_{**}<\delta_{*}<\frac{1}{3}$ and a mapping  
$$(\delta,\gamma)\mapsto \Omega^{(\delta,\gamma)}: (0,\delta_{**})\times \T \to  \left(-\frac{1}{3},\frac{1}{3}\right)$$
that is 2$\pi$-periodic in $\gamma$ and so that 
$\mathcal{S}[\delta, \Omega^{(\delta,\gamma)}]=0$ on $(0,\delta_{**})\times \T.$
\item \label{Prop:final_reduction-branching} (Branching) The mapping $(\delta,\gamma) \mapsto \mathcal{S}[\delta,\Omega^{(\delta,\gamma)}]$  can be extended continuously in an unique fashion  to a mapping
$$(\delta,\gamma)\mapsto\overline{\mathcal{S}}[\delta, \Omega^{(\delta,\gamma)}]=0,$$
on $\delta \in [0,\delta_{**})\times\T$.
Moreover, we must have
\begin{align*}
 \lim_{\delta\downarrow 0 }\Omega^{(\delta,\gamma)} =0,
\end{align*}
and further, we have  $\frac{\delta}{4} \leq \ep = \ep(\delta,\Omega) \leq 4\delta.$
\end{enumerate}
In particular, assumption (\hyperref[H2]{H2}) is fully satisfied for this parametrization.
\end{Lemma}
\begin{Proof} We know that $H^4(\R)$ is an algebra (cf. \cite[Corollary 8.10]{Brezis}), and from the fact that pointwise evaluation is continuous, thanks to \eqref{Sobolev_embedding}. Hence, using \eqref{Ansatz-rewritten} and 
Proposition \ref{Prop:reduced}, the continuity of the mapping 
$$ (\delta,\Omega,\gamma) \mapsto \mathcal{U}^2\Big\vert_{x=0} $$
is obtained. Since we know from  \eqref{rolls} that  the parametrization of the rolls is continuous in $H_{\text{per}}^4(\R)$, we can use the Sobolev embedding as in \eqref{Sobolev_embedding} to derive coninuity in the uniform norm. Thus, assertion (i) follows.
 Recall from (\hyperref[H2]{H2}) that $\om^2 = 1 + \delta \Omega$. Plugging the expansions derived in Lemma \ref{Lemma:hamiltonian_expansion}, we obtain
\begin{align*}
 0 = \frac{- \delta^2\mathcal{U}^2\Big\vert_{x=0} - \delta^2 -  \mathcal{H}^{(l)}[\ur^{(\delta,\sqrt{1+ \delta\Omega},\gamma)}(\cdot)] + \mathcal{H}^{(r)}[0]}{\delta^2} = - \mathcal{U}^2\Big\vert_{x=0} +\frac{4}{3} \, \Omega^{3} - \frac{4}{3} \, \Omega + \mathcal{O}\left(\delta^2\right),
\end{align*}
which we finally rewrite as
\begin{align}\label{last_contraction_mapping}
 \Omega = - \Omega^3 -  \frac{3}{4} \mathcal{U}^2\Big\vert_{x=0} + \mathcal{O}\left(\delta^2\right) =: \mathcal{G}[\Omega, \delta,\gamma].
\end{align}
 Hence,  $\Omega$ can be seen as a fixed point for the mapping $\Omega \mapsto \mathcal{G}[\Omega, \delta,\gamma]$, which is 2$\pi$-periodic in $\gamma$. An application of the Contraction Mapping Theorem (in the topology of the sup norm) then shows that 
$\Omega$ is parametrized by $(\delta,\gamma)$; in combination with (i) this implies that the parametrization holds in a continuous fashion  and that its periodicity in $\gamma$ persists, thus establishing (ii).  

Last, (iii) is derived from \eqref{last_contraction_mapping} once we take the limit $\delta \downarrow 0$: we conclude that 
$$\displaystyle{\lim_{\delta \downarrow 0}\left\{\Omega^{(\delta,\gamma)} + (\Omega^{(\delta,\gamma)})^3\right\}=0 }.$$
Since $\Omega^{(\delta,\gamma)} \in \left(-\delta_*,\delta_*\right)\subset \left(-\frac{1}{3},\frac{1}{3}\right)$ and the mapping y $\mapsto y + y^3$  is a diffeomorphism in a neighborhood of 0, we conclude that $\displaystyle{\lim_{\delta\downarrow0} \Omega^{(\delta,\gamma)} =\Omega_*} $ exists and this limit must satisfy
$$ \Omega_* = - (\Omega_*)^3$$
uniformly in $\gamma \in \T, $ for all the bounds only depend on $\delta$, namely, they are uniform in $\gamma$. As $\Omega_* \in \left[-\delta_*,\delta_*\right]\subset \left[-\frac{1}{3},\frac{1}{3}\right]$, we must have that $\Omega_* =0$. Last, since $\Omega \in \left(-\frac{1}{3}, \frac{1}{3}\right)$ we can the parametrization result referred to in \eqref{rolls} to have, $\frac{\delta}{4} \leq \ep = \ep(\delta,\Omega) \leq 4\delta$ (choosing $\delta_*>0$ smaller if necessary), and we are done.
\end{Proof}

Before putting things together, we go back to the stretching $\displaystyle{x\mapsto z:= \frac{x}{\om^{(\delta,\gamma)}}}$ that lead us to \eqref{asymptotics_properties}:

\begin{Corollary}[Existence of solutions to problem \eqref{SH-eq} and its reformulation \eqref{asymptotics_properties}]\label{Corollary:solution_to_the_problem} There exists an $\delta_{**}>0$ such that the mapping \eqref{Ansatz-rewritten}
\begin{align*}
 (x,\delta,\Omega,\gamma) \mapsto \mathcal{U}(x) = v\n^{(\delta,\Omega,\gamma)}(x) +v\f[\delta, \Omega,v\n^{(\delta,\Omega,\gamma)}](x)+ \chi(\ep x)\ur^{(\delta,\sqrt{1+ \delta\Omega},\gamma)}\left(x\right), 
\end{align*}
solves problem \eqref{asymptotics_properties} for all $(\delta,\gamma) \in [0,\delta_{**})\times \T$, and $\om =\om^{(\delta,\gamma)}$; moreover,  $\mathcal{U}^{(\delta,\gamma)}(\cdot)$ has the properties described in Definition \ref{heterecolinic_def}.
 
 In particular, changing back to the variable $x = \om^{(\delta,\gamma)}z$ we obtain
 \begin{align}\label{cont_loss_of_cont}
[0,\delta_{**})\times \T \ni (\delta,\gamma) \mapsto \mathcal{U}^{(\delta,\gamma)}(z) =  v\n(\om^{(\delta,\gamma)} z) + v\f(\om^{(\delta,\gamma)} z)+\chi(\ep \om^{(\delta,\gamma)} z)\ur^{(\delta,\sqrt{1+ \delta\Omega^{(\delta,\gamma)}},\gamma)}(\om^{(\delta,\gamma)} z),  
 \end{align}
 which solves problem \eqref{SH-eq}.
\end{Corollary}

\begin{Lemma}[Continuity and loss of continuity of \eqref{cont_loss_of_cont}  in the sup norm with respect to parameters $(\delta,\gamma)$]\label{Lem:to_be_not_to_be_continuous}
 Consider the mapping \eqref{cont_loss_of_cont}. Thus, the mapping $(\delta,\gamma) \mapsto \mathcal{U}^{(\delta,\gamma)}(\cdot) $ is continuous in the sup norm whenever  $\om^{(\delta,\gamma)} \equiv 1$ (equivalently, when $\Omega^{(\delta,\gamma)} \equiv 0$), that is,
 $$\left \Vert \mathcal{U}^{(\delta_1,\gamma_1)}(\cdot) - \mathcal{U}^{(\delta_2,\gamma_2)}(\cdot)\right\Vert_{L^{\infty}(\R)}\lesssim \vert\delta_1 -\delta_2\vert+ \vert\gamma_1 -\gamma_2\vert$$
 holds whenever $\om^{(\delta,\gamma)} \equiv 1$ for all $\gamma_1, \gamma_2 \in \T$ and sufficiently small $\delta_1,\delta_2 >0$. 
 
 On the other hand, this mapping is discontinuous in the sup norm when $\om^{(\delta,\gamma)} \not\equiv 1$.
\end{Lemma}
\begin{Proof}
 According to Proposition \ref{Prop:reduced} we already know that the mapping $(\delta,\Omega,\gamma) \mapsto v\n(\cdot) + v\f(\cdot)$ is continuous in $H^2(\R)$, which readily implies the result once we apply the Sobolev Embedding $H^2(\R) \hookrightarrow L^{\infty}(\R)$. 
 Thus, using the parametrization $(\delta,\gamma) \mapsto \Omega^{(\delta,\gamma)}$ derived in Lemma \ref{Lemma:wavenumber_selec}, in order to investigate the continuity of the mapping \eqref{cont_loss_of_cont} it suffices to investigate  the mapping
 $$ (\delta,\gamma) \mapsto \chi(\ep \om^{(\delta,\gamma)}x)\ur^{(\delta,\sqrt{1+ \delta\Omega^{(\delta,\gamma)}},\gamma)}(\om^{(\delta,\gamma)}x),$$
 First, notice that  $\ur^{(\delta,\sqrt{1+ \delta\Omega},\gamma)}(\cdot) \not \in L^2(\R)$, although  its $L^{\infty}(\R)$ norm scales in $\ep$, cf. Proposition \ref{scaling:u_roll}. 
  Next, take $(\delta_1,\gamma_1)$ and $(\delta_2,\gamma_2)$ in $\left[0,\delta_{**}\right)\times\T$ and write $\displaystyle{\ur^{(j)}(\cdot):= \ur^{(\delta_{j},\om^{(\delta_{j},\gamma_{j})},\gamma_{j})}(\cdot)}$, for $j \in \{1,2\}$, and $\om^{(\delta,\gamma)} = \sqrt{1 + \delta \Omega^{(\delta,\gamma)}}$. 
  
  The dichotomy with respect to $\om^{(\delta,\gamma)}$  is due to the following fact: when $\om^{(\delta,\gamma)}\equiv 1$ the period of $\ur^{(\delta, \gamma)}(\cdot)$ is fixed, legitimating the equality
 \begin{align}\label{continuity_L_infty}
  \Vert \ur^{(1)}(\cdot) - \ur^{(2)}(\cdot) \Vert_{L^{\infty}(\R)} = \Vert \ur^{(1)}(\cdot) - \ur^{(2)}(\cdot) \Vert_{L^{\infty}(\left[0,2\pi\right])}. 
 \end{align}
Since $\Vert \ur^{(1)}(\cdot) - \ur^{(2)}(\cdot) \Vert_{L^{\infty}(\left[0,2\pi\right])} \lesssim \left\vert \ep(\delta_1,\Omega^{(\delta_{1})}) - \ep(\delta_2,\Omega^{(\delta_{2})})\right\vert+ \vert\gamma_1 - \gamma_2 \vert,$ we immediately obtain the result. 
In contrast, if $\om^{(\delta,\gamma)}\not \equiv 1$, we only get
 \begin{align*}
  \left\Vert \ur^{(1)}(x) - \ur^{(2)}(x) \right\Vert_{L^{\infty}(\R)} \approx \left\vert \ep(\delta_1,\Omega^{(\delta_1)}) + \ep(\delta_2,\Omega^{(\delta_2)})\right\vert,
 \end{align*}
 in which case continuity does not hold.
\end{Proof}
After this result, the proof of Theorem \ref{Thm:main} is readily available:
\begin{Proof}[of Theorem \ref{Thm:main}] The proof unfolds as a successive derivation  of equivalent but reduced formulations of the same problem, so we just compile the result in the order it was constructed. From the very beginning we fix $\gamma\in \R$. Equation \eqref{SH-eq-bifurcation-1} sets an equation for the ``corrector'' $v(\cdot)$, which gets solved from  \S \ref{section:far_near} through \S \ref{sec:reduced_equation+approx+solv}; along  this resolution process, the following steps were taken:
\begin{enumerate}[label=(\roman*), ref=\theTheorem(\roman*)]

\hangindent\leftmargin
 \item In \S \ref{section:far_near} it was shown that $v(\cdot)$ could be written as $v(\cdot)= v\n(\cdot) + v\f(\cdot)$, and using a Lyapunov-Schmidtt reduction the equation \eqref{SH-eq-bifurcation-1} could be rewritten in an equivalent form as a system of coupled equations \eqref{mode_system_cutoff-a1}, \eqref{mode_system_cutoff-a2}. Then, in Proposition \ref{Prop:near_far_domination} of \S\ref{section:blow_up} it was proved that the parametrization%
 \begin{align*}
  (v\n(\cdot), \delta,\Omega,\gamma) \mapsto v\f[v\n, \delta,\Omega,\gamma](\cdot) 
 \end{align*}
holds in a continuous fashion. 
\item  Thanks to the previous results, the problem gets reduced to understanding $(v\n(\cdot), \delta,\Omega,\gamma)$ only. In \S \eqref{section:reduced_eq_simplifications} a blow-up in Fourier space is introduced in order to desingularize the limit $\ep\downarrow0$. In passing, we obtain the following characterization of $v\n(\cdot)$ given in \eqref{blow_up:near},
\begin{align*}
 v_{near}(x)& = \ep^{\beta} e^{+i  x} \gp(\ep^{\beta} x) + \ep^{\beta} e^{-i x} \gm(\ep^{\beta} x). 
\end{align*}
and, accordingly, many properties of  $\gm(\cdot)$ and $\gp(\cdot)$ (and, consequently, of $v\n(\cdot)$) are derived.
\item We show in Proposition \ref{Prop:decay} that $\displaystyle{\lim_{\vert x\vert \to +\infty}\partial_x^{\alpha}v(x) = 0}$ whenever $\alpha \in \{0,1,2,3\}$; this proves that the asymptotic spatial limits in stated in (i) hold.
\item In \S \ref{section:reduced_eq_simplifications},  a pair  of equivalent  equations \eqref{second_main_eqt:rewrite} on $\gm(\cdot)$ and $\gp(\cdot)$  is derived. After slightly adapting a result of \cite{Fefferman} to our purposes, in Proposition \eqref{beta_is_1:sufficiency} we are able to conclude that  $\beta =1$, which allow us to simplify our equations even further, reducing them to the form \eqref{second_main_eqt:simplified} presented in section   \S \ref{sec:5_4}.
\item In  \S \ref{sec:reduced_equation+approx+solv} we obtain (iii), showing in Lemma \ref{Lemma:almost_partition} that $\chi(\cdot)$ can be chosen in such a way that it satisfies the ODE
\begin{align*}
 -4\pi \partial_x^2 \chi(\cdot) +\frac{3 \pi}{4}\left( \chi^3(\cdot) - \chi(\cdot)\right) - c\partial_x\chi(x) =0, \quad \lim_{x\to -\infty}\chi(x) = 1, \quad \lim_{x\to -\infty}\chi(x) = 0.
\end{align*}
In passing, this provides nice conditions under which, once more,  the Contraction Mapping Theorem can be applied (Proposition \ref{Prop:reduced}) to show  the existence of a mapping
$$ (0,\delta_{*})\times \left(-\delta_*,\delta_*\right)\times \T\ni (\delta,\Omega,\gamma) \mapsto \left(\gm^{(\delta,\Omega,\gamma)}(\cdot),\gp^{(\delta,\Omega,\gamma)}(\cdot)\right) \in  H_{\text{near}, \ep^{\tau-1}}^{2}(\R)\times H_{\text{near}, \ep^{\tau-1}}^{2}(\R),$$
that is $2\pi$-periodic in $\gamma$ and solves problem \eqref{SH-eq-bifurcation-1} for 
$$v (\cdot) = v\n^{(\delta,\Omega,\gamma)} + v\f[v\n^{(\delta,\Omega,\gamma)},\delta,\Omega],$$ 
where $(0,\delta_{*})\times \left(-\delta_*,\delta_*\right)\times \T\ni (\delta,\Omega,\gamma)\mapsto  v_{near}(x) = \ep e^{+i  x} \gp^{(\delta,\Omega,\gamma)}(\ep x) + \ep e^{-i x} \gm^{(\delta,\Omega,\gamma)}(\ep x)$, which we extend as a mapping on $[0,\delta_{*})\times \left(-\delta_*,\delta_*\right)\times \T$, thanks to Proposition \ref{pmtr_range:extension_to_0}. 
\end{enumerate}
This long derivation is brought to full use in \S \ref{section:wavenumber_selection}, where we put  the Hamiltonian structure of the problem in the limelight: in Lemma \ref{Lemma:wavenumber_selec} we prove the selection mechanism asserted in (ii) using  again the Contraction Mapping Theorem: we have an implicit description of $\Omega$ in terms of $(\delta,\gamma)$ and continuity of the mapping
\begin{align*}
 (\delta,\gamma) \mapsto \Omega^{(\delta,\gamma)}, \quad \text{where} \quad \Omega^{(\delta,\gamma)}\Big\vert_{\delta =0} = 0 \quad \text{holds},\quad \forall \gamma \in \T;
\end{align*}
furthemore, this mapping is 2$\pi$ periodic in $\gamma$. When plugged back into $v\n^{(\delta,\Omega,\gamma)}(\cdot)$ we immediately obtain  (iv). 

Subsequently, we prove (i) as follows: first,  in Corollary \ref{Corollary:solution_to_the_problem} we obtain a solution to equation \eqref{SH-eq} of the form
\begin{align*}
 (\delta,\gamma) \mapsto \mathcal{U}^{(\delta,\gamma)}(\om^{(\delta,\gamma)} z) =  v\n(\om^{(\delta,\gamma)} z) + v\f(\om^{(\delta,\gamma)} z)+\chi(\ep \om^{(\delta,\gamma)} z)\ur^{(\delta,\sqrt{1+ \delta\Omega^{(\delta,\gamma)}},\gamma)}(\om^{(\delta,\gamma)} z),
\end{align*}
 Lemma \ref{Lemma:regularity} implies the regularity of this mapping in  $x$.  Lemma \ref{Lem:to_be_not_to_be_continuous} asserts  its regularity with respect to  $(\delta, \gamma)$ where the following dichotomy was proved:continuity  in the topology of uniform continuity holds if and only if $\om^{(\delta,\gamma)}\equiv 1$.  This establishes (i), and we are finally done with the proof of Theorem \ref{Thm:main}. \end{Proof}


\section{Open problems and further comments}\label{open_problems}
This project was highly inspired by the techniques and the perspective in the memoir \cite{Fefferman}. Our results are corroborated by the results in \cite{Weinburd}, which puts the ideas we advocate for in a safe ground for comparison with other mathematical tools. Some of the questions we address below are suggestively related to well established mathematical techniques (\S \ref{open_problems:techniques}-\S\ref{open_problems:multiplier}), while others (\S \ref{open_problems:invasion_fronts}-\S\ref{open:stability}) have a pure speculative nature; regardless of their plausibility, they should be read with caution. 
\subsection{The techniques in \cite{Fefferman}}\label{open_problems:techniques} The results in the memoir \cite{Fefferman} in many aspects seem to be related to the work of Schneider and touches upon interesting issues previously addressed in \cite{Newell, Newell-order_parameters}. It is also possible that the techniques presented in the referred memoir can provide existence of stationary solutions, as those  with different wavenumbers opposite sides of the  far-field, as considered in \cite{Doel-memoir}. On the other hand, in many cases an  invasion by a modulated front can be seen, as studied by Collet and Eckmann \cite{collet1986existence}. In that  case, a careful use of the Bloch-Floquet theory as done in \cite{Fefferman} indicates a first step in the construction of these objects using partial differential equations, multi-scale analysis and perturbation techniques. 
\subsection{Far/near reductions and  dynamical systems}\label{open_problems:far_near} As pointed out before, the shape of $v\n(\cdot)$ in  \eqref{blow_up:near} resembles the  initial steps in modulation theory that lead to a Ginzburg-Landau type of equation. In general this type of approximations are applied to pattern formation systems close to unstable states \cite[Part IV]{Schneider-book}; different approaches to this derivation using multiple scales analysis are also possible, cf. \cite{Harten}. With regards to the role of the reduced equation and approximation operator $\widetilde{\mathcal{R}^{(0,\Omega;\pm)}}[\widehat{\gm}, \widehat{\gp}]$ studied in \S \ref{sec:reduced_equation+approx+solv},  its relation to transversality theory, as used in \cite[Sections 2-(e)(f)(g)]{Weinburd}, is still to be clarified.
\subsection{On the role of multipliers}\label{open_problems:multiplier} In case of $n$ distinct singularities in the multiplier one can expect to obtain a system of ODEs in $n$ variables. More examples and  possibly a more general theory is still necessary  to elucidate how the location of the singularities, the stretching of the vicinities around them (imposed by the far/near decomposition) and the subsequent blow-up in Fourier space play a crucial role in the reduced equations obtained in the end.  Many other questions are also left behind: is it possible to say  that the reduced equations are always unique (or, somewhat, equivalent), up to some parametrization? 
It would also be interesting to see an example of this technique being applied to  non-local models  where the linearization has an associated multiplier directly obtained by convolution; this scenario is very interesting, because these cases are not directly amenable to dynamicl systems.  The techniques presented in \cite{ScheelTao} show an interesting direction of investigation, where the authors prove the existence of stationary solutions to a nonlocal (convolution-type) problem exploiting properties of the multiplier, resulting in a reduced type of equation \cite[\S 3]{ScheelTao}; in their case however  no far/field (spectral nor spatial) decomposition is used.

\subsection{Invasion fronts and the role of $\chi(\cdot)$}\label{open_problems:invasion_fronts} In the literature of pattern formation, whenever near/far (spatial) decompositions have been  applied, the functions $\chi(\cdot)$ are mostly introduced aiming localization of pattern properties in the far field (cf. \cite{Goh_Sch-CH,Goh_Sch-pattern,Lloyd-Sc, Morrissey,Monteiro_Scheel,Monteiro_Scheel-contact_angle,Mon-horizontal}), which happens mostly because  effects of bifurcation parameter variations to the the far field are known or predictable. In these cases, the role os the function $\chi(\cdot)$ is essentially that of a \emph{partition of unity}. 

In contrast, the type of the ODE satisfied by $\chi(\cdot)$ in  Lemma \ref{Lemma:almost_partition} bring to our mind the results in  \cite{collet1986existence,CE90}, based on  which several analogies can be made: one can say that the function $\chi(\cdot)$  plays a role of an envelope of the modulated (invasion) front. In the case \cite{collet1986existence} however, the interpretation is clearer: as the profile is positive,  one can easily discern  the invaded part from the wake of the front.

Last, we mention the interesting work \cite{Nepomnyashchy} (in particular, \S III), where multi-dimensional patterns are studied. Roughly speaking, the boundaries  of regions filled with rolls with different orientation are investigated using functions $r_1(\cdot)$ and $r_2(\cdot)$ that have a similar role to that of $\chi(\cdot)$; see also \S \ref{open_problems:grain_boundaries} below. 

\subsection{Nanopatterns and numerical aspects}\label{open_problems:numerical}  In case of patterns with characteristic wavelength smaller than computer floating numbers, our result seems to be useful in the computation of possible profiles that would be otherwise undetectable in numerical simulations. The papers \cite{Lions_micro,Lions_micro-bigger} exploit this among other questions in an interesting fashion, with techniques different from ours. 

\subsection{Grain boundaries, defects, and multi-dimensional patterns}\label{open_problems:grain_boundaries} In spite of its robustness,   dynamical systems techniques do not seem to be broad enough to capture unsurmountable difficulties in the study of multidimensional patterns.
Recently, different research avenues have been exploited: several studies have been done using rigorous numerical analysis \cite{Spots}, harmonic analysis techniques \cite{Jaramillo,Jaramillo-inhomogeneities,Lions_micro}, variational techniques  \cite{Rabinowitz}, or more functional-analytic based techniques \cite{Monteiro_Scheel-contact_angle}. Still, many classes of problems remain unsolved, as that of asymmetrical grain boundaries, a case that does not seem to be directly amenable to the spatial dynamics techniques as presented in \cite{Haragus,Qiliang}; in this scenario the far/near decompositions we presented might be relevant for analytical results. For numerical results which exploit far/near (spatial) decomposition, see \cite{Lloyd-Sc}).
 
\subsection{Directional quenching and wavenumber selection}\label{open_problems:dq_wavenumber_selec} In spite of the consistency between our result and that in \cite[Theorem 1.1]{Weinburd}, a  comparison between them shows that our selection mechanism result is weaker, for we only prove the existence of an implicit  parametrization $\Omega = \Omega^{(\delta,\gamma)}$, while therein the authors show a much stronger resul: an explicit formula for the lower order terms,
$$\Omega^{(\delta,\gamma)} = \delta\frac{\cos(2\gamma)}{16} + \mathcal{O}\left(\delta^{3}\right).$$
We must further highlight  that the shape of the control parameter $\mu(\cdot)$ was chosen in its simplest form as a parameter of jump type. The sharp discontinuity type was called here \emph{directional quenching}, and should be seen as a contrast to the case of  slow decay, physically closer to the process of \emph{annealing};  it is plausible that, in the latter scenario, no wavenumber selection happens. We highlight the interesting discussion in \cite[\S 4]{Weinburd} about wavenumber selection by ramp discontinuities; for a general overview of pattern selection mechanisms, the discussion in \cite[\S 3.1]{Nishiura} is also very illustrative. 

\subsection{Stability issues}\label{open:stability} Once we consider roll solutions embedded in a multidimensional space new types of instabilities are seen,  as it is the case of Zig-Zag instabilities (cf. \cite[\S 4]{Mielke}, \cite{Newell-order_parameters}); several mathematical and numerical issues related to the formation of these patterns are still not well understood (some steps to clarify them can be seen in \cite{Growing_stripes}). Nonlinear stability results have also been investigated, with some interesting recent results in one-dimension: the seminal paper \cite{SandsSch-defects}, the recent work \cite{2018-source-defects} and the memoir \cite{Doel-memoir}.

%
\appendix
%
\section{A few computations}

In this appendix we give a detailed computation of the results presented in \S \ref{sec:5_2} and \S\ref{sec:5_3}. We stress that part of the results in this appendix are computations necessary before we set $\beta =1$. This, in part \ref{appendix:comput-4} we still carry $\beta$  in computations, while in part \ref{appendix:comput-1} we already set $\beta =1$.

\subsection{Simplifying \texorpdfstring{$\mathcal{I}^{(4; \, \pm)}(\xi)$}{I [+-1](4)(xi)}: computations }\label{appendix:comput-4} 
Our goal is to estimate the terms in
\begin{align*}
\mathcal{I}^{(4; \, +)}(\xi)  &= \widetilde{\mathbb{P}\n^{(+)}}\circ\mathscr{F}\left(   \mathscr{N}^{(4)}[v\n,\ur^{(\delta,\om,\gamma)},\ep]\right)(\ep^{\beta}\xi)\nonumber\\
 & = \widetilde{\mathbb{P}\n^{(+)}}\circ\mathscr{F}\left(  \chi( \chi^2 -1)\left(\ur^{(\delta,\om,\gamma)}\right)^3 \right)(\ep^{\beta}\xi)  + \widetilde{\mathbb{P}\n^{(+)}}\circ\mathscr{F}\left( [(1+ \partial_x^2)^2,\chi]\ur \right)(\ep^{\beta}\xi)\\
 & \quad + \widetilde{\mathbb{P}\n^{(+)}}\circ\mathscr{F}\left(\delta^2(\mu -1)\chi \ur^{(\delta,\om,\gamma)} \right)(\ep^{\beta}\xi)\nonumber\\
 & = \mathcal{S}_{\text{I}}^{(4; \, +)}(\xi) + \mathcal{S}_{\text{II}}^{(4; \, +)}(\xi) + \mathcal{S}_{\text{III}}^{(4; \, +)}(\xi),
\end{align*} 
The term $\mathcal{S}_{\text{I}}^{(4; \, +)}(\cdot)$ has already been studied. For now, we estimate  $\mathcal{S}_{\text{II}}^{(4; \, +)}(\xi)$. In order to do so,  initially we study the term $[(1+ \partial_x^2)^2,\chi]\ur^{(\delta,\om,\gamma)}(x)$, which can be written, from lowest order terms to higher order terms,  as
\begin{align*}\label{expansion_commutant}
[(1+ \partial_x^2)^2,\chi]\ur^{(\delta,\om,\gamma)}(x) &= 4\partial_x(\chi(\ep^{\beta} x))(\partial_x^3 + \partial_x) \ur^{(\delta,\om,\gamma)}(x) +\partial_x^2(\chi_x(\ep x))[6\partial_x^2 +2]\ur^{(\delta,\om,\gamma)}(x)\nonumber\\
&\quad + \mathcal{O}\left((\partial_x^3 + \partial_x^4)\chi(\ep^{\beta}x )\right) \\
& = G_1 + G_2 + G_3,\nonumber
\end{align*}
Note that the relation $\partial_x^{j}\chi(\ep^{\beta}x) = \mathcal{O}(\ep^{j\beta})$ holds for all $j \in \N; $ with this in mind we begin the analysis of each term. We immediately have that $\vert G_3\vert = \mathcal{O}(\ep^{3\beta})$ and further, $\Vert G_3\Vert_{L^2(\R)}^2 = \mathcal{O}(\ep^{5\beta})$.

For the first term on the left hand side, we rely on  \eqref{rolls} to rewrite is as
\begin{align*}
G_1 &= 4 \ep^{\beta} \chi'(\ep^{\beta} x)(\partial_x^3 + \partial_x)\ur^{(\delta,\om,\gamma)}(x) = 4 \ep^{\beta} \chi'(\ep^{\beta} x)(\partial_x^3 + \partial_x) (\ep\cos(x + \gamma) + \mathcal{O}(\ep^3)),
\end{align*}
and since $(\partial_x^2 + 1)\cos(x + \gamma) =0$, we only need to keep the second term in the expansion of $\ur^{(\delta,\om,\gamma)}(\cdot)$, obtaining $\displaystyle{\vert G_1\vert = \mathcal{O}\left(\ep^{3+\beta}\right)}$, and $\Vert G_1\Vert_{L^2(\R)}^2= \mathcal{O}(e^{6+\beta})$. 

Last we  study $G_2$, once again expanding $\ur^{(\delta,\om,\gamma)}(\cdot)$:
\begin{align*}
 G_2& = \partial_x^2(\chi(\ep^{\beta} x))[6\partial_x^2 +2] (\ep\cos(x + \gamma) + \mathcal{O}(\ep^3)) = -4\ep^{1+2\beta}\cos(x+ \gamma)\chi^{''}(\ep^{\beta} x) +\widetilde{G_2},
\end{align*}
where $\vert \widetilde{G_2}\vert = \mathcal{O}(\ep^{3+2\beta})$ and $\Vert \widetilde{G_2}\Vert_{L^2(\R)}^2 = \mathcal{O}(e^{6+3\beta}).$ Altogether, we finally write $\mathcal{S}_{\text{II}}^{(4; \, +)}(\xi)$ as
\begin{align*}
 \mathcal{S}_{\text{II}}^{(4; \, +)}(\xi)  & =  -4\ep^{1+2\beta} \mathbbm{1}_{\left\{ \xi \,\in \,\ep^{\tau-\beta}\mathcal{B}\right\}}(\xi)\mathscr{F}\left[\cos(\cdot+ \gamma)\chi^{''}(\ep^{\beta} \cdot)\right](1 + \ep^{\beta}\xi)  + F_{\text{II}}\nonumber \\
 & = -4\ep^{1+2\beta}\mathbbm{1}_{\left\{ \xi \,\in \,\ep^{\tau-\beta}\mathcal{B}\right\}}(\xi)\int_{\R} \left( \cos(x+ \gamma)\chi^{''}(\ep^{\beta} x)  \right)e^{-i x}e^{-i \ep^{\beta} \xi x}\mathrm{d}x + F_{\text{II}}\nonumber \\
 & \stackrel{x= 2\pi z, }{=}-8\pi \ep^{1+2\beta}\mathbbm{1}_{\left\{ \xi \,\in \,\ep^{\tau-\beta}\mathcal{B}\right\}}(\xi)\int_{\R} \left( \cos(2\pi z+ \gamma)\chi^{''}(\ep^{\beta} 2\pi z)  \right)e^{-i 2\pi z}e^{-i \ep^{\beta} \xi 2\pi z}\mathrm{d}z + F_{\text{II}}\nonumber \\
 & \stackrel{\theta = 2\pi \ep^{\beta} }{=}-8\pi\left(\frac{\theta}{2\pi}\right)^{\frac{1+2\beta}{\beta}}\mathbbm{1}_{\left\{ \xi \,\in \,\ep^{\tau-\beta}\mathcal{B}\right\}}(\xi)\int_{\R} \left( \cos(2\pi z+ \gamma)\chi^{''}(\theta z)  \right)e^{-i 2\pi z}e^{-i \theta \xi z}\mathrm{d}z + F_{\text{II}},
\end{align*}
where $F_{\text{II}} = G_1 + \widetilde{G_2} + G_3$, and $\Vert F_{\text{II}} \Vert_{L^2(\R)}^2 =\mathcal{O}(\ep^{6+\beta}+\ep^{6+3\beta}+\ep^{5\beta})$. We rewrite it using  Lemma \ref{Lemma:fef_pois}: setting $Z = \theta z$, $f(z) = \cos(2\pi z+ \gamma)$, $g(z) = e^{-i 2\pi z}$ and $\Gamma_2(z, Z) = \Gamma_2(Z)= \chi^{''}(\theta z) =    \chi^{''}(Z)$, we get
\begin{align*}
 \mathcal{S}_{\text{II}}^{(4; \, +)}(\xi)  
 &= -8\pi\left(\frac{\theta}{2\pi}\right)^{\frac{1+\beta}{\beta}}\mathbbm{1}_{\left\{ \xi \,\in \,\ep^{\tau-\beta}\mathcal{B}\right\}}(\xi)\sum_{n\in \mathbb{Z}}\int_0^1 e^{2\pi i n x} \widehat{\Gamma_2}\left(\frac{2\pi n}{\theta} + \xi\right)\cos(2\pi z+ \gamma)e^{-i 2\pi z}\mathrm{d}z+F_{\text{II}}\nonumber \\
 &= -8\pi\left(\frac{\theta}{2\pi}\right)^{\frac{1+\beta}{\beta}}\mathbbm{1}_{\left\{ \xi \,\in \,\ep^{\tau-\beta}\mathcal{B}\right\}}(\xi)\times \\
 & \quad \times\left[  \widehat{\Gamma_2}( \xi)\int_0^1\cos(2\pi z+ \gamma)e^{-i 2\pi z}\mathrm{d}z + \sum_{\vert n \vert \geq 1} \widehat{\Gamma_2}\left(\frac{2\pi n}{\theta} + \xi\right)\int_0^1\cos(2\pi z+ \gamma) e^{2\pi i n x}e^{-i 2\pi z}\mathrm{d}z\right]+F_{\text{II}}\nonumber \\
 &= -8\pi\ep^{1+\beta}\mathbbm{1}_{\left\{ \xi \,\in \,\ep^{\tau-\beta}\mathcal{B}\right\}}(\xi) e^{i\gamma}\widehat{\Gamma_2}( \xi)\int_0^1\cos^2(2\pi z+ \gamma)\mathrm{d}z + \widetilde{F_{\text{II}}},
\end{align*}
where  estimates on $\Vert F_{\text{II}} \Vert_{L^2(\R)}$ and Lemma \ref{Lemma:fef_pois} imply that  $\Vert\widetilde{ F_{\text{II}}} \Vert_{L^2(\R)}^2 =\mathcal{O}(\ep^{6+\beta}+\ep^{6+3\beta}+\ep^{5\beta}+ \ep^{2(1 + 2\beta)})$. Note also that $\Gamma_2(\cdot) \in L^2(\R)$, thanks to (\hyperref[H1]{H1}). Hence, after further simplification using Lemma \ref{Lem:simplf_cos} we get
\begin{align*}
 \mathcal{S}_{\text{II}}^{(4; \, +)}(\xi)  &= -4\pi\ep^{1+\beta}\mathbbm{1}_{\left\{ \xi \,\in \,\ep^{\tau-\beta}\mathcal{B}\right\}}(\xi) e^{i\gamma}\widehat{\Gamma_2}( \xi) + \widetilde{F_{\text{II}}},
\end{align*}
Last, we use similar ideas to estimate $\mathcal{S}_{\text{III}}^{(4; \, +)}(\xi)$: 
\begin{align*}
 \mathcal{S}_{\text{III}}^{(4; \, +)}(\xi)  &=  \delta^2\mathbbm{1}_{\left\{ \xi \,\in \,\ep^{\tau-\beta}\mathcal{B}\right\}}(\xi)\mathscr{F}\left[(\mu(\cdot) -1)\chi(\ep^{\beta} \cdot)\ur^{(\delta,\om,\gamma)}(\cdot)\right](1 + \ep^{\beta}\xi) \nonumber \\
 &=  \ep\delta^2\mathbbm{1}_{\left\{ \xi \,\in \,\ep^{\tau-\beta}\mathcal{B}\right\}}(\xi)\mathscr{F}\left[(\mu(\cdot) -1)\chi(\ep^{\beta} \cdot)\widetilde{\ur^{(\delta,\om,\gamma)}}(\cdot)\right](1 + \ep^{\beta}\xi) \nonumber \\
 &=  -2\delta^2\ep\mathbbm{1}_{\left\{ \xi \,\in \,\ep^{\tau-\beta}\mathcal{B}\right\}}(\xi)\int_{0}^{\infty}\cos(x + \gamma)\chi(\ep^{\beta} x)e^{-i(1+\ep^{\beta} \xi )x} \mathrm{d} x + F_{\text{III}}\\
 & = -2\delta^2\ep^{1-\beta}\mathbbm{1}_{\left\{ \xi \,\in \,\ep^{\tau-\beta}\mathcal{B}\right\}}(\xi)\widehat{\Gamma_3^{(+)}}( \xi) + F_{\text{III}},
\end{align*}
where the last equality is due to property \eqref{Fourier_properties:dilation}. The remainder is estimated as $\Vert F_{\text{III}}  \Vert_{L^2(\R)}^2 = \mathcal{O}(\delta^4\ep^{2(3-\beta)})\stackrel{(\hyperref[H2]{H2})}{=} \mathcal{O}(\ep^{2(5-\beta)})$,  and  the quantity $\Gamma_3^{(\pm)}(\cdot)$ is written explicitly as
$$\Gamma_3^{(\pm)}(\cdot) =  \mathbbm{1}_{[0,\infty)}(\cdot) \chi\left(\cdot \right)e^{-i\frac{(\cdot)}{\ep}}\cos\left(\frac{\cdot}{\ep}\pm \gamma\right),$$
where $\Gamma_3^{(\pm)}(\cdot) \in L^2(\R)$ due to (\hyperref[H1]{H1}). At last, we combine these  estimates to get
\begin{align*}
 \mathcal{I}^{(4; \, +)}(\xi)  = \mathbbm{1}_{\left\{ \xi \,\in \,\ep^{\tau-\beta}\mathcal{B}\right\}}(\xi)\left\{\frac{3 \pi \ep^{3-\beta}}{4}e^{+ i \gamma}\widehat{\Gamma_1}\left( \xi\right) - 4\pi \ep^{1+\beta}e^{ + i \gamma}\widehat{\Gamma_2}( \xi) - 2\delta^2 \ep^{1-\beta}\widehat{\Gamma_3^{(+)}}(\xi) \right\} +\mathscr{M}^{(4; \, +)}(\xi).
\end{align*}
The following relations hold
\begin{align*}
 \mathscr{M}^{(4; \, +)}(\xi) = \mathbbm{1}_{\left\{ \xi \,\in \,\ep^{\tau-\beta}\mathcal{B}\right\}}(\xi)\mathscr{M}^{(4; \, +)}(\xi),\quad  \text{with} \quad \Vert \mathscr{M}^{(4; \, +)}\Vert_{L^2(\R)}^2=\mathcal{O}(\ep^{6+\beta}+\ep^{6+3\beta}+\ep^{5\beta}+ \ep^{2(1 + 2\beta)} + \ep^{2(5-\beta)})
\end{align*}
A similar analysis gives
\begin{align*}
 \mathcal{I}^{(4; \, -)}(\xi)  = \mathbbm{1}_{\left\{ \xi \,\in \,\ep^{\tau-\beta}\mathcal{B}\right\}}(\xi)\left\{\frac{3 \pi \ep^{3-\beta}}{4}e^{- i \gamma}\widehat{\Gamma_1}\left( \xi\right) - 4\pi \ep^{1+\beta}e^{ - i \gamma}\widehat{\Gamma_2}( \xi) - 2\delta^2 \ep^{1-\beta}\widehat{\Gamma_3^{(-)}}(\xi) \right\} +\mathscr{M}^{(4; \, -)}(\xi).
\end{align*}
where $\mathscr{M}^{(4; \, -)}(\xi) = \mathbbm{1}_{\left\{ \xi \,\in \,\ep^{\tau-\beta}\mathcal{B}\right\}}(\xi)\mathscr{M}^{(4; \, -)}(\xi)$,  and $\Vert \mathscr{M}^{(4; \, -)}\Vert_{L^2(\R)}=\mathcal{O}(\ep^{6+\beta}+\ep^{6+3\beta}+\ep^{5\beta}+ \ep^{2(1 + 2\beta)} + \ep^{2(5-\beta)})$.
\subsection{Simplifying \texorpdfstring{$\mathcal{I}^{(1; \, \pm)}(\xi)$}{I [+-1](1)(xi) }}\label{appendix:comput-1} 
\begin{SCfigure}[][htb]
 \includegraphics[width=6cm]{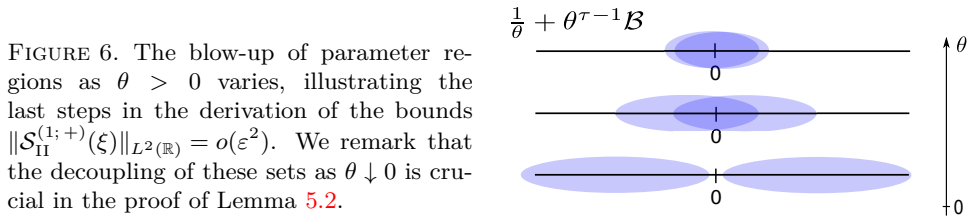}
 \caption{The blow-up of parameter regions as $\theta>0$ varies,  illustrating the last steps in the derivation of the bounds $\displaystyle{ \Vert \mathcal{S}_{\text{II}}^{(1; \, +)}(\xi)  \Vert_{L^2(\R)}= o(\ep^2)}$. We remark that the decoupling of these sets as $\theta \downarrow 0$  is crucial  in the proof of Lemma \ref{Lemma:fef_pois}. \label{no_overlap} }
\end{SCfigure}

Due to Proposition \ref{beta_is_1:sufficiency} in this section we set $\beta=1$. We aim to understand the terms
\begin{align*}
\mathcal{I}^{(1; \, \pm)}(\xi)& =-\mathbbm{1}_{\left\{ \ep^{\tau-1}\mathcal{B}\right\}}(\xi) \mathscr{F}\left[\left(\delta^2\mu(x) - 3 (\widetilde{\ur^{(\delta,\om,\gamma)}})^2\right) v\n(\cdot)\right](1 + \ep\xi).
\end{align*}
For the moment, let's focus on $\mathcal{I}^{(1; \, +)}(\xi)$. Using the relation \eqref{blow_up:near}, we can expand this term as
\begin{align*}
\mathcal{I}^{(1; \, +)}(\xi)& = -\delta^2\ep \mathbbm{1}_{\left\{ \ep^{\tau-1}\mathcal{B}\right\}}(\xi)\int_{\R} \mu(x) \gp(\ep x) e^{-i \ep \xi x}\mathrm{d}x\nonumber \\
& \quad -\delta^2\ep \mathbbm{1}_{\left\{ \ep^{\tau-1}\mathcal{B}\right\}}(\xi)\int_{\R} \mu(x)  \gm(\ep x)e^{-2ix} e^{-i \ep \xi x}\mathrm{d}x\nonumber \\
& \quad  +3\ep^3 \mathbbm{1}_{\left\{ \ep^{\tau-1}\mathcal{B}\right\}}(\xi)\int_{\R}  (\widetilde{\ur^{(\delta,\om,\gamma)}}(x))^2\gp(\ep x) e^{-i \ep \xi x}\mathrm{d}x\nonumber \\
& \quad  +3\ep^3 \mathbbm{1}_{\left\{ \ep^{\tau-1}\mathcal{B}\right\}}(\xi)\int_{\R}  (\widetilde{\ur^{(\delta,\om,\gamma)}}(x))^2 \gm(\ep x)e^{-2ix} e^{-i \ep \xi x}\mathrm{d}x\nonumber \\
& = \mathcal{S}_{\text{I}}^{(1; \, +)}(\xi) + \mathcal{S}_{\text{II}}^{(1; \, +)}(\xi) + \mathcal{S}_{\text{III}}^{(1; \, +)}(\xi) + \mathcal{S}_{\text{IV}}^{(1; \, +)}(\xi).
\end{align*}
In order to estimate $\mathcal{S}_{\text{I}}^{(1; \, +)}(\xi)$,
we do a  change variables $X = \ep x$,  and exploit the homogeneity of the parameter $\mu(\cdot)$, which enables us to  write $\mu\left(\frac{X}{\ep}\right) = \mu(X)$. Thus,
\begin{align*}
 \mathcal{S}_{\text{I}}^{(1; \, +)}(\xi) &=  -\delta^2 \mathbbm{1}_{\left\{ \ep^{\tau-1}\mathcal{B}\right\}}(\xi) \int_{\R} \mu(X)\left[ \gp(X) \right] e^{-i  \xi X}\mathrm{d}X.
\end{align*}
With regards to the case $\mathcal{S}_{\text{II}}^{(1; \, +)}(\xi) $, the analysis is easier. Indeed, 
\begin{align*}
 \mathcal{S}_{\text{II}}^{(1; \, +)}(\xi)  &= -\delta^2\ep \mathbbm{1}_{\left\{ \ep^{\tau-1}\mathcal{B}\right\}}(\xi)\int_{\R} \mu(x)\left[\gp(\ep x)e^{-i2x}\right] e^{-i \ep \xi x}\mathrm{d}x\nonumber \\
 &= -\delta^2 \mathbbm{1}_{\left\{ \ep^{\tau-1}\mathcal{B}\right\}}(\xi)\int_{\R} \mu(X)\left[\gp(X)\right] e^{-i\left(\frac{2}{\ep}+ \xi \right)X}\mathrm{d}X \\
  &= -\delta^2 \mathbbm{1}_{\left\{ \ep^{\tau-1}\mathcal{B}\right\}}(\xi)\mathscr{F}\left[\mu(\cdot)\gp(\cdot)\right]\left(\frac{2}{\ep}+ \xi \right).
\end{align*}
As $\frac{2}{\ep} + \ep^{\tau -1}\mathcal{B}\subset \{\vert \xi \vert \geq \frac{1}{\ep}\},$ we make use of \eqref{rolls} to get $\displaystyle{ \Vert \mathcal{S}_{\text{II}}^{(1; \, +)}(\xi)  \Vert_{L^2(\R)}^2 \leq \ep^4 \int_{\vert  \xi\vert  \geq \frac{1}{\ep}}\left\vert \mathscr{F}\left[\mu\left(\cdot\right)\gp(\cdot)\right](\xi) \right\vert^2 \mathrm{d} \xi,}$ thus obtaining that  $\displaystyle{ \Vert \mathcal{S}_{\text{II}}^{(1; \, +)}(\xi)  \Vert_{L^2(\R)}= o(\ep^2)};$ see Figure \ref{no_overlap}. 

\noindent
Both terms $\mathcal{S}_{\text{III}}^{(1; \, +)}(\xi)$ and $\mathcal{S}_{\text{IV}}^{(1; \, +)}(\xi)$ can be estimated using formula in Lemma \ref{Lemma:fef_pois}. Indeed, after change of variables, we set $\theta = 2\pi \ep$ and apply \eqref{poisson_summation}, obtaining the equivalent expression
\begin{equation}
 \begin{split}
  \mathcal{S}_{\text{III}}^{(1; \, +)}(\xi)  & = 3\ep^3 \mathbbm{1}_{\left\{ \ep^{\tau-1}\mathcal{B}\right\}}(\xi)\int_{\R}  (\widetilde{\ur^{(\delta,\om,\gamma)}}(x))^2\left[\gp(\ep x) \right] e^{-i \ep \xi x}\mathrm{d}x\nonumber \\
 & = 3\frac{\theta^2}{2\pi}\widehat{\gp}(\xi) \mathbbm{1}_{\left\{ \ep^{\tau-1}\mathcal{B}\right\}}(\xi)\int_0^1  (\widetilde{\ur^{(\delta,\om,\gamma)}}(2\pi z))^2 \mathrm{d}z + G_{\text{III}} 
 \end{split}
\end{equation}
with estimate $\Vert G_{\text{III}} \Vert_{L^2(\R) }= \mathcal{O}(\ep^3).$ Likewise, 
\begin{equation}
 \begin{split}
  \mathcal{S}_{\text{IV}}^{(1; \, +)}(\xi)   &=  3\ep^3 \mathbbm{1}_{\left\{ \ep^{\tau-1}\mathcal{B}\right\}}(\xi)\int_{\R}  (\widetilde{\ur^{(\delta,\om,\gamma)}}(x))^2\left[ \gm(\ep x)e^{-2ix}\right] e^{-i \ep \xi x}\mathrm{d}x\nonumber \\
 & =  3\frac{\theta^2}{2\pi}\widehat{\gm}(\xi) \mathbbm{1}_{\left\{ \ep^{\tau-1}\mathcal{B}\right\}}(\xi)\int_0^1  (\widetilde{\ur^{(\delta,\om,\gamma)}}(2\pi z))^2e^{-2i2\pi z}\mathrm{d}z+ G_{\text{IV}}.
 \end{split}
\end{equation}
where  $\Vert G_{\text{IV}} \Vert_{L^2(\R) }= \mathcal{O}(\ep^3)$. And using the expansion of $\ur^{(\delta,\om,\gamma)}(\cdot)$ given in Lemma \ref{scaling:u_roll} and Lemma \ref{Lem:simplf_cos} we get
\begin{align*}
 \mathcal{I}^{(1; \, +)}(\xi)  = \mathbbm{1}_{\left\{ \xi \,\in \,\ep^{\tau-\beta}\mathcal{B}\right\}}(\xi)\left\{-\delta^2 \int_{\R} \mu(X)\left[ \gp(X) \right] e^{-i  \xi X}\mathrm{d}X +  3\pi \ep^2\widehat{\gp}(\xi)  + 3\frac{\pi\ep^2}{2}\widehat{\gm}(\xi) \right\} +\mathscr{M}^{(1; \, +)}(\xi).
\end{align*}
where $\mathscr{M}^{(1; \, +)}(\xi) = \mathbbm{1}_{\left\{ \xi \,\in \,\ep^{\tau-\beta}\mathcal{B}\right\}}(\xi)\mathscr{M}^{(1; \, +)}(\xi)$ is so that 
$\Vert \mathscr{M}^{(1; \, +)} \Vert_{L^2(\R)} = \mathcal{O}\left(\ep^3\right) $. A similar analysis gives
\begin{align*}
 \mathcal{I}^{(1; \, -)}(\xi)  = \mathbbm{1}_{\left\{ \xi \,\in \,\ep^{\tau-\beta}\mathcal{B}\right\}}(\xi)\left\{-\delta^2 \int_{\R} \mu(X)\left[ \gm(X) \right] e^{-i  \xi X}\mathrm{d}X +  3\pi \ep^2\widehat{\gm}(\xi)  + 3\frac{\pi\ep^2}{2}\widehat{\gp}(\xi) \right\} +\mathscr{M}^{(1; \, -)}(\xi),
\end{align*}
where $\mathscr{M}^{(1; \, -)}(\xi) = \mathbbm{1}_{\left\{ \xi \,\in \,\ep^{\tau-\beta}\mathcal{B}\right\}}(\xi)\mathscr{M}^{(1; \, -)}(\xi)$ is so that $\Vert \mathscr{M}^{(1; \, -)} \Vert_{L^2(\R)} = \mathcal{O}\left(\ep^3\right).$

  \section{On the existence of traveling waves }\label{existence_TW}
\begin{SCfigure}[][h]
 \includegraphics[height=3cm]{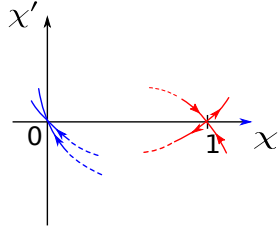}
 \caption{Phase plane sketch for the ODE $\partial_x^2\chi +  c\partial_x\chi +  f(\chi)=0$, which has a Heteroclinic orbit $\chi(\cdot)$ satisfying  $\chi(-\infty)=1$ and $\chi(+\infty)=0.$\label{fig_heteroclinic}}
\end{SCfigure}

Following the analysis  in \cite[\S 4]{Fife}, whenever  $c>0$ the ODE
$$\partial_x^2\chi +  c\partial_x\chi +  f(\chi)=0, $$
with $f(0) = f(1)=0$, $f'(0) >0$, $f'(1) <0$   and $f(x) >0$ whenever $x\in (0,1)$ admits a  heteroclinic orbit $\chi(\cdot) $  so that $\displaystyle{\lim_{x\to -\infty}\chi(x) = 1}$, and  $\displaystyle{\lim_{x\to +\infty}\chi(x) =0};$ see Figure \ref{fig_heteroclinic}. 

\noindent
It is also known that this convergence is exponential, that is,  for all $j \in \{0, \ldots, 4\}$ there exist constants $C_*>0$ and $S_*>0$ such that
\begin{align*}
 \vert \partial_x^j\left(\chi(x) - 1\right)\vert \lesssim e^{-C_* \vert x \vert},  \quad \text{whenever} \quad x\leq -S_*,\quad \text{and} \quad 
 \quad \vert  \partial_x^j\chi(x)\vert \lesssim e^{-C_* \vert x \vert}, \quad \text{whenever} \quad x\geq S_*.
\end{align*}
Furthermore, whenever $G(\cdot) = f'(\cdot)$ is known, it is possible to compute $\Vert\chi \Vert_{L^{2}(\R)}^2$: indeed, 
\begin{align*}
 c\int_{\R}\vert\chi'(s)\vert^2\mathrm{d}s = -\int_{\R}\chi'(s)\left(\chi^{''}(c) +G'(\chi) \right)\mathrm{d}s = -\left(\frac{\chi'(x)^2}{2} + G(\chi(x))\right)\Big \vert_{x= -\infty}^{x= +\infty}.
\end{align*}
Thus,
\begin{align*}
\int_{\R}\vert\chi'(s)\vert^2\mathrm{d}s \leq \left\vert \frac{G(\chi(+\infty)) - G(\chi(-\infty) )}{c} \right\vert = \left\vert \frac{G(1) - G(0)}{c} \right\vert.
\end{align*}
Moreover, 
\begin{align*}
\left\vert c^2\int_{\R}\vert\chi'(s)\vert^2\mathrm{d}s\right\vert \lesssim \vert c\vert.
\end{align*}
In the particular, when $x\mapsto f(x) := x - x^3$ we have $\displaystyle{G(x) = C - \frac{(x^{2 -1})^{2}}{4} }$ for $C\in \R$, yielding $$\left\vert c^2\int_{\R}\vert\chi'(s)\vert^2\mathrm{d}s\right\vert \leq \frac{\vert c\vert}{4}.$$
%

\newcommand{\etalchar}[1]{$^{#1}$}

\end{document}